\documentclass[a4paper,12pt]{amsart}
\input xypic

\usepackage{amssymb}

\def\hbar{\bar{h}}

\def\Rarr#1{\buildrel #1\over \longrightarrow}

\def\iso{\buildrel \sim\over\to}

\def\Gm{{\mathfrak{m}}}
\def\Gn{{\mathfrak{n}}}

\def\CA{{\mathcal{A}}}

\def\CC{{\mathcal{C}}}

\def\CE{{\mathcal{E}}}
\def\CF{{\mathcal{F}}}
\def\CG{{\mathcal{G}}}

\def\CI{{\mathcal{I}}}
\def\CJ{{\mathcal{J}}}
\def\CK{{\mathcal{K}}}
\def\CL{{\mathcal{L}}}

\def\CO{{\mathcal{O}}}

\def\CS{{\mathcal{S}}}
\def\CT{{\mathcal{T}}}

\def\BA{{\mathbf{A}}}

\def\BF{{\mathbf{F}}}

\def\BL{{\mathbf{L}}}

\def\BN{{\mathbf{N}}}

\def\BP{{\mathbf{P}}}

\def\BZ{{\mathbf{Z}}}

\def\eps{\varepsilon}

\def\add{\operatorname{add}\nolimits}

\def\amp{\operatorname{amp}\nolimits}

\def\can{{\mathrm{can}}}

\def\coh{\operatorname{coh}\nolimits}

\def\mcoh{\operatorname{\!-coh}\nolimits}

\def\coker{\operatorname{coker}\nolimits}
\def\colim{\operatorname{colim}\nolimits}

\def\cone{\operatorname{cone}\nolimits}

\def\diff{\operatorname{diff}\nolimits}

\def\en{{\mathrm{en}}}
\def\End{\operatorname{End}\nolimits}

\def\Ext{\operatorname{Ext}\nolimits}

\def\gldim{\operatorname{gldim}\nolimits}

\def\hocolim{\operatorname{hocolim}\nolimits}

\def\Hom{\operatorname{Hom}\nolimits}

\def\Id{\operatorname{Id}\nolimits}
\def\id{\operatorname{id}\nolimits}
\def\im{\operatorname{im}\nolimits}

\def\Ind{\operatorname{Ind}\nolimits}

\def\mMod{\operatorname{\!-mod}\nolimits}

\def\mMOD{\operatorname{\!-Mod}\nolimits}

\def\opp{{\operatorname{opp}\nolimits}}

\def\mperf{\operatorname{\!-perf}\nolimits}

\def\pdim{\operatorname{pdim}\nolimits}

\def\mproj{\operatorname{\!-proj}\nolimits}
\def\mProj{\operatorname{\!-Proj}\nolimits}

\def\mqcoh{\operatorname{\!-qcoh}\nolimits}

\def\Res{\operatorname{Res}\nolimits}

\def\Spec{\operatorname{Spec}\nolimits}

\def\mstab{\operatorname{\!-stab}\nolimits}

\def\Tor{\operatorname{Tor}\nolimits}

\def\ie{{\em i.e.}}
\def\eg{{\em e.g.}}

\def\tM{{\tilde{M}}}

\newtheorem{thm}{Theorem}[section]
\newtheorem{lemma}[thm]{Lemma}
\newtheorem{cor}[thm]{Corollary}
\newtheorem{prop}[thm]{Proposition}
\newtheorem{defi}[thm]{Definition}
\newtheorem{conj}[thm]{Conjecture}

\theoremstyle{definition}
\newtheorem{rem}[thm]{Remark}
\newtheorem{example}[thm]{Example}

\xyoption{all}

\def\en{{\operatorname{en}\nolimits}}
\def\Krulldim{\operatorname{Krulldim}\nolimits}
\def\repdim{\operatorname{repdim}\nolimits}
\def\wrepdim{\operatorname{wrepdim}\nolimits}
\def\lwrepdim{\operatorname{lwrepdim}\nolimits}
\def\rwrepdim{\operatorname{rwrepdim}\nolimits}
\def\Ll{\operatorname{ll}\nolimits}
\def\opp{\circ}

\title{Dimensions of triangulated categories}
\author{Rapha\"el Rouquier}

\begin{document}
\date{Septembre 2004}

\begin{abstract}
We define a dimension for a triangulated category. We prove a representability
Theorem for a class of functors on finite dimensional triangulated
categories. We study the dimension of the bounded
derived category of an algebra or a scheme and we show in particular
that the bounded derived category of coherent sheaves over a
variety has a finite dimension. For a self-injective algebra, a lower
bound for Auslander's representation dimension is given by the dimension
of the stable category. We use this to compute the representation
dimension of exterior algebras. This provides the first known examples
of representation dimension $>3$. We deduce that the
Loewy length of the group algebra over $\BF_2$ of a finite group
is strictly bounded below by the
$2$-rank of the group (a conjecture of Benson).
\end{abstract}

\maketitle

\tableofcontents

\section{Introduction}
In his 1971 Queen Mary College notes \cite{Au}, Auslander introduced an
invariant
of finite dimensional algebras, the representation dimension. It was meant
to measure how far an algebra is to having only finitely many classes of
indecomposable modules.
Whereas many upper bounds have been found for the representation dimension,
lower bounds were missing. In particular, it wasn't known whether
the representation dimension could be greater than $3$. 
A proof that all algebras have representation dimension at most $3$ would have
led for example to a solution of the finitistic dimension conjecture
\cite{IgTo}.

We prove here that the representation dimension of the exterior
algebra of a finite dimensional vector space is one plus the dimension of
that vector space --- in particular, the representation dimension can
be arbitrarily large.
Thus, the representation dimension is a useful invariant of finite
dimensional algebras of infinite representation type,
confirming the hope of Auslander. The case of algebras
with infinite global dimension is particularly interesting.

As a consequence of our results, we prove the characteristic $p=2$ case
of a conjecture of Benson asserting that the $p$-rank of a finite group
is less than the Loewy length of its group algebra over a field
of characteristic $p$.

\medskip
Our approach to these problems is to define and study
a ``dimension'' for triangulated categories. This is inspired by
Bondal and Van den Bergh's work \cite{BoVdB} and we generalize some of
their main results. This leads us to look more generally at finiteness
conditions for triangulated categories and their meaning in algebraic
and geometric examples.

\smallskip
Our results also shed some light on properties of dg algebras
related to geometry and might be viewed as requirements for non-commutative
geometry. Let us give two examples.
\begin{itemize}
\item
Given a projective scheme $X$ over a field, it is a classical fact that
there exists a dg algebra $A$ with $D(A)\simeq D(X)$~:
going to the dg world, $X$ ``becomes affine''. Given any such dg algebra $A$,
we show that dg-$A$-modules with finite dimensional
total cohomology admit ``resolutions'' (Remark \ref{resolproj}), a strong
condition on a dg algebra.

\item
Given a quasi-projective scheme $X$ over a perfect field, we show that
there is a dg algebra $A$ with $A\mperf\simeq D^b(X\mcoh)$. Furthermore,
for any such $A$, $A\mperf$ has finite dimension, a property
which might be viewed as some kind of homological regularity for $A$
(when $A=A_0$ is noetherian and the differential vanishes, then
$\dim A\mperf<\infty$ if and only if $\gldim A<\infty$)~:
going to the dg world, $X$ ``becomes regular''. Note that this 
notion is weaker that smoothness of a dg algebra $A$ =
perfection of $A$ as an $(A,A)$-dg bimodule.
\end{itemize}
More generally,
let $\CT$ be the bounded derived category of finitely generated modules
over an artinian ring or over
a noetherian ring with finite global dimension, or
the bounded derived category of coherent sheaves over a separated scheme of
finite type over a perfect field, or a quotient of any of these categories.
Then, $\dim\CT<\infty$, \ie, $\CT$ is equivalent to 
$A\mperf$, where $A$ is a dg algebra which is ``homologically regular''.

\medskip
Let us review the content of the chapters. Chapters
\S \ref{secdim}-\ref{seclocalization} deal with ``abstract'' triangulated
categories, whereas chapters \S \ref{subsecalg}-\ref{secfd} deal
with derived categories of rings and schemes (quasi-compact, quasi-separated),
and stable categories of self-injective finite dimensional algebras.

\smallskip
In a first part \S \ref{secdim}, we review various types of generation of
triangulated categories and we define a dimension for triangulated
categories. This is the minimum number of cones needed to build any
object (up to a summand) from finite sums of shifts of a given object.
Note that we introduce and use later the notion of
compactness of objects for triangulated categories that do not
admit arbitrary direct sums.

\smallskip
We consider various finiteness conditions for cohomological
functors on a triangulated category in \S\ref{subsecfin} and we derive
some stability properties of these classes of functors. We
define in particular locally finitely presented functors. On an
$\Ext$-finite triangulated category, they include locally finite functors.
The crucial property of locally finitely presented functors is that they
can be ``approximated'' by representable functors (\S\ref{seclfp} and
in particular Proposition \ref{approxlfp}). This leads, in
\S\ref{secrep}, to representability Theorems for locally finitely
presented functors, generalizing Brown-Neeman's representability Theorem
for ``big'' triangulated categories (cocomplete, generated by a set
of compact objects) as well as Bondal-Van den Bergh's Theorem
for ``small'' triangulated categories ($\Ext$-finite, with finite dimension).

In \S\ref{subsecfinob}, we consider properties of objects $C$ related
to properties of the functor $\Hom(-,C)$ restricted to compact objects.
Later (Corollary \ref{lfgnoetherian} and Proposition \ref{clfpcoh}),
we determine the corresponding categories
for derived categories of noetherian algebras or schemes~:
the cohomologically locally finitely presented objects are the complexes
with bounded and finite type cohomology.

\smallskip
Part \S\ref{seclocalization} develops a formalism for coverings of
triangulated categories mimicking the coverings of schemes by open
subschemes. More precisely, we consider Bousfield subcategories and
we introduce a notion of proper intersection of two Bousfield subcategories
(\S\ref{secproper}) and we study properties of families of Bousfield
subcategories intersecting properly. We obtain for example
Mayer-Vietoris triangles (Proposition \ref{mayervietoris}). The main 
part is \S \ref{seccoverings}, where we consider compactness. We
show that compactness is a ``local'' property (Corollary \ref{caractcompact})
--- this sheds some light
on the compact=perfect property for derived categories of schemes. We also
explain how to construct a generating set of compact objects from
local data (Theorem \ref{cocovering}). It is fairly quick to prove
the existence of a compact generator for the derived categories of schemes
from this (cf Theorem \ref{genscheme} for a version with supports).

\smallskip
Part \S \ref{subsecalg} is a study of various classes of objects in
derived categories of algebras and schemes. The main section
\S \ref{secschemes} considers complexes of $\CO$-modules with
quasi-coherent cohomology on a quasi-compact and quasi-separated scheme.
We show how the length of the cohomology (sheaves) of a complex
is related to the non-zero shifted groups of morphisms from a
fixed compact generator to the complex (Proposition \ref{amplitude}). We
give a characterization of pseudo-coherent complexes in triangulated
terms (Proposition \ref{pseudocoherent})~: they are the objects whose
cohomology can be ``approximated''
by compact objects --- such a result is classical in the presence of
an ample family of line bundles.

In \S \ref{seccompactbounded}, we show that
for noetherian rings or noetherian separated schemes, the compact objects
of the bounded derived category are the objects with finite type
cohomology. This gives a descent principle.

\smallskip
In \S \ref{secder}, we analyze the dimension of derived categories in
algebra and geometry. In \S \ref{secresoldiago}, we use resolution
of the diagonal methods. We show that for $A$ a finite dimensional
algebra or a commutative algebra essentially of finite type over a 
perfect field, then $\dim D^b(A\mMod)\le\gldim A$ (Proposition
\ref{boundgldimk}). For a smooth quasi-projective scheme $X$ over
a field, we have $\dim D^b(X\mcoh)\le 2\dim X$ (Proposition \ref{boundsmooth}).
We show (Proposition \ref{localregular}) that the
residue field of a commutative local noetherian algebra $A$ over a field
cannot be obtained by less than $\Krulldim(A)$ cones from sums of
$A$ and its shifts. This is the key result to get lower bounds~: we deduce
(Proposition \ref{inegvar}) that for $X$ a reduced separated scheme
of finite type over a field, then $\dim D^b(X\mcoh)\ge \dim X$ and
there is equality $\dim D^b(X\mcoh)=\dim X$ when $X$ is in addition
smooth and affine (Theorem \ref{smoothaffine}).

In \S \ref{secfingldim}, we investigate rings with finite global dimension
and regular schemes. As noted by Van den Bergh, a noetherian ring is regular if
and only $\dim A\mperf<\infty$ (Proposition \ref{regularring}). 
Analogously, the category of perfect complexes for a quasi-projective
scheme $X$ over a field has finite dimension if and only if $X$ is regular
(Proposition \ref{regularvar}).
For an artinian ring $A$, then $\dim D^b(A\mMod)$ is less than the Loewy length
(Proposition \ref{boundLoewy}).

The main result of \S \ref{schemes} is a proof
that the derived category of coherent
sheaves $D^b(X\mcoh)$
has finite dimension, for a separated scheme $X$ of finite type over a
perfect field (Theorem \ref{finitevar}). This is
rather surprising and it is a rare instance where $D^b(X\mcoh)$ is
better behaved than $X\mperf$. As a consequence, the stable derived
category $D^b(X\mcoh)/X\mperf$ has finite dimension as well.
In the smooth case (only), one has a stronger result about
the structure sheaf of the diagonal, due to Kontsevich.
There are very few cases where we can determine the exact dimension of
$D^b(X\mcoh)$ for a smooth $X$, and these are cases where it coincides with
$\dim X$ ($X$ affine or $X$ a projective space for example).
We conclude the chapter (Corollary \ref{swap})
by a determination of locally finite cohomological
functors on $X\mperf$ and $D^b(X\mcoh)^\opp$, for $X$ a projective
scheme over a perfect field $k$ (the first case is due to Bondal and
Van den Bergh)~: they are represented by an object of $D^b(X\mcoh)$ in
the first case and an object of $X\mperf$ in the second case --- this
exhibits some ``perfect pairing''
$\Hom(-,-):X\mperf\times D^b(X\mcoh)\to k\mMod$.

\smallskip
Finally, in \S \ref{secfd}, we study the dimension of the stable
category of a self-injective algebra, in relation with Auslander's
representation dimension. Via Koszul duality, we compute these
dimensions for the exterior algebra of a finite dimensional vector space~:
$\dim \Lambda(k^n)\mstab=n-1$ and the representation dimension of
$\Lambda(k^n)$ is $n+1$ (Theorem \ref{exterior}).
This enables us to settle the characteristic $2$ case of a conjecture of
Benson (Theorem \ref{Benson2}).

\medskip
Preliminary results have been obtained and exposed at the
conference ``Twenty years of tilting theory'' in Fraueninsel in November 2002.
I wish to thank the
organizers for giving me the opportunity to report on these early results and
the participants for many useful discussions,
particularly Thorsten Holm for introducing me to Auslander's work.

The geometric part of this work was motivated by lectures given by
A.~A.~Beilinson at the University of Chicago and by discussions with A.~Bondal.

\section{Notations and terminology}
For $\CC$ an additive category and $\CI$ a subcategory of $\CC$,
we denote by $\add(\CI)$ (resp. $\overline{\add}(\CI)$)
the smallest additive full subcategory of $\CC$ containing $\CI$ and
closed under taking direct summands (resp. and closed under direct
summands and direct sums).
We say that $\CI$ is {\em dense} if every object of $\CC$ is isomorphic to a
direct summand of an object of $\CI$.

We denote by $\CC^\opp$ the category opposite to $\CC$.
We identify a set of objects of $\CC$ with the full subcategory with
the corresponding set of objects.

\smallskip
Let $\CT$ be a triangulated category. A {\em thick} subcategory $\CI$
of $\CT$ is a full triangulated subcategory such that given
$M,N\in\CT$ with $M\oplus N\in\CI$, then $M,N\in\CI$. Whenever we consider
the quotient $\CT/\CI$, it will be assumed that this has small
$\Hom$-sets.

Given
$X\Rarr{f} Y\Rarr{g} Z\rightsquigarrow$ a distinguished triangle, then $Z$
is called a {\em cone} of $f$ and $X$ a {\em cocone} of $g$.

\smallskip
Given $\CA$ an abelian category, we denote by $D(\CA)$ the derived
category of $\CA$ and we denote by $D^{\le a}(\CA)$ the full subcategory
of objects with cohomology vanishing in degrees $>a$.

\smallskip
Let $A$ be a differential graded (=dg) algebra.
We denote by $D(A)$ the derived category of dg $A$-modules and
by $A\mperf$ the category of {\em perfect
complexes}, \ie, the smallest thick subcategory of $D(A)$ containing $A$.

\smallskip
Let $A$ be a ring. We denote by $A\mMOD$ the category of left
$A$-modules, by $A\mMod$ the category of
finitely generated left $A$-modules, by 
$A\mProj$ the category of projective $A$-modules and by
$A\mproj$ the category of finitely generated projective $A$-modules.
We denote by $\gldim A$
the global dimension of $A$. For $M$ an $A$-module, we denote by
$\pdim_A M$ the projective dimension of $M$. We denote by
$A^\opp$ the opposite ring to $A$. For $A$ an algebra over a 
commutative ring $k$,
we put $A^\en=A\otimes_k A^\opp$.

\smallskip
Let $X$ be a scheme. We denote by $X\mcoh$ (resp. $X\mqcoh$) the category of
coherent (resp. quasi-coherent) sheaves on $X$.
We denote by 
$D(X)$ the full subcategory of the derived category of sheaves of
$\CO_X$-modules consisting of complexes with quasi-coherent cohomology.
A complex of sheaves of $\CO_X$-modules is
{\em perfect} if it is locally quasi-isomorphic to a bounded complex of
vector bundles (=locally free sheaves of finite rank). We denote
by $X\mperf$ the full subcategory of perfect complexes of $D(X)$.
Given a complex of sheaves $C$, the notation
$H^i(C)$ will always refer to the cohomology sheaves, not to the
(hyper)cohomology groups.

\smallskip
Let $C$ be a complex of objects of an additive category and $i\in\BZ$.
We put $\sigma^{\le i}C=\cdots\to C^{i-1}\to C^i\to 0$ and
$\sigma^{\ge i}C=0\to C^i\to C^{i+1}\to \cdots$. 
Let now $C$ be a complex of objects of an abelian category.
We put $\tau^{\ge i}C=0\to C^i/\im d^{i-1}\to C^{i+1}\to C^{i+2}\to\cdots$
and $\tau^{\le i}C=\cdots\to C^{i-2}\to C^{i-1}\to \ker d^i\to 0$.

\section{Dimension}
\label{secdim}

\subsection{Dimension for triangulated categories}

\subsubsection{}
We review here various types of generation of triangulated categories,
including the crucial ``strong generation'' due to Bondal and Van den Bergh.

\smallskip
Let $\CT$ be a triangulated category.

Let $\CI_1$ and $\CI_2$ be two subcategories of $\CT$. We denote
by $\CI_1\ast\CI_2$ the full subcategory of $\CT$
consisting of objects $M$ such that there is a distinguished triangle
$M_1\to M\to M_2\rightsquigarrow$
with $M_i\in\CI_i$.

\medskip
Let $\CI$ be a subcategory of $\CT$.
We denote by
$\langle\CI\rangle$ 
the smallest full subcategory of $\CC$ containing $\CI$ and closed
under finite direct sums, direct summands and shifts.
We denote by $\overline{\CI}$ the smallest full subcategory of $\CC$
containing $\CI$ and closed under direct sums and shifts.

We put $\CI_1\diamond\CI_2=\langle \CI_1\ast\CI_2\rangle$.

We put $\langle\CI\rangle_0=0$ and
we define by induction
$\langle\CI\rangle_i=\langle\CI\rangle_{i-1}\diamond\langle\CI\rangle$
for $i\ge 1$. We put
$\langle\CI\rangle_\infty=\bigcup_{i\ge 0}\langle\CI\rangle_i$.
We define also $\CI^{\ast i}=\CI^{\ast (i-1)}\ast\CI$.

The objects of $\langle\CI\rangle_i$ are the direct summands of
the objects obtained by taking an $i$-fold extension of finite direct
sums of shifts of objects of $\CI$.

We will also write $\langle \CI\rangle_{\CT,i}$ when there is some ambiguity
about $\CT$.

\smallskip

We say that
\begin{itemize}
\item $\CI$ {\em generates} $\CT$ if
given $C\in\CT$ with $\Hom_\CC(D[i],C)=0$ for all $D\in\CI$ and all
$i\in\BZ$, then $C=0$
\item $\CI$ is a {\em $d$-step generator} of $\CT$ if
$\CT=\langle\CI\rangle_d$ (where $d\in\BN\cup\{\infty\}$)
\item $\CI$ is a {\em complete $d$-step generator} of $\CT$ if
$\CT=\langle\overline{\CI}\rangle_d$ (where $d\in\BN\cup\{\infty\}$).
\end{itemize}

We say that $\CT$ is
\begin{itemize}
\item {\em finitely generated} if there exists $C\in\CT$ which generates $\CT$
(such a $C$ is called a {\em generator})
\item {\em classically finitely (completely) generated}
if there exists $C\in\CT$
which is a (complete) $\infty$-step generator of $\CT$ (such a $C$ is called
a {\em classical (complete) generator})
\item {\em strongly finitely (completely) generated} if there exists $C\in\CT$
which is a (complete) $d$-step generator of $\CT$ for some $d\in\BN$
(such a $C$ is called a {\em strong (complete) generator}).
\end{itemize}

Note that $C$ is a classical generator of $\CT$ if and only if
$\CT$ is the smallest thick subcategory of $\CT$ containing $C$.
Note also that if $\CT$ is strongly finitely generated, then
every classical generator is a strong generator.


It will also be useful to allow only certain infinite direct sums.
We define $\widetilde{\CI}$ to be the smallest full subcategory of $\CT$ closed
under finite direct sums and shifts and containing multiples
of objects of $\CI$ (\ie, for $X\in\CI$ and $E$ a set such that
$X^{(E)}$ exists in $\CT$, then $X^{(E)}\in\widetilde{\CI}$).

\subsubsection{}
We now define a dimension for a triangulated category.

\begin{defi}
The dimension of $\CT$, denoted by
$\dim\CT$, is the minimal integer $d\ge 0$ such that there is
$M$ in $\CT$ with $\CT=\langle M\rangle_{d+1}$.

We define the dimension to be $\infty$ when there is no such $M$.
\end{defi}

\smallskip
The following Lemmas are clear.

\begin{lemma}
\label{densedim}
Let $\CT'$ be a dense full triangulated subcategory of $\CT$.
Then, $\dim\CT=\dim\CT'$.
\end{lemma}

\begin{lemma}
\label{quotient}
Let $F:\CT\to\CT'$ be a triangulated functor with dense image. If
$\CT=\langle \CI\rangle_d$, then $\CT'=\langle F(\CI)\rangle_d$. So,
$\dim\CT'\le\dim\CT$.

In particular,
let $\CI$ be a thick subcategory of $\CT$. Then,
$\dim \CT/\CI\le \dim\CT$.
\end{lemma}

\begin{lemma}
\label{decomposition}
Let $\CT_1$ and $\CT_2$ be two triangulated subcategories of $\CT$ such that
$\CT=\CT_1\diamond\CT_2$. Then, $\dim\CT\le 1+\dim\CT_1 + \dim\CT_2$.
\end{lemma}

\begin{lemma}
\label{generationopp}
The property of generation, strong generation, etc... for $\CT$ is
equivalent to the corresponding property for $\CT^\opp$.
We have $\dim\CT^\circ=\dim\CT$.
\end{lemma}

\subsection{Remarks on generation}
\subsubsection{}
\begin{rem}
One can strengthen the notion of generation of $\CT$ by $\CI$ by
requiring that $\CT$ is the smallest triangulated subcategory
containing $\CI$ and closed under direct sums.
Cf Theorem \ref{equivgen} for a case
where both notions coincide.
\end{rem}

\begin{rem}
Let $\langle\CI\rangle'$ be the smallest full subcategory of $\CT$ containing
$\CI$ and closed under finite direct sums and shifts. Define similarly
$\langle\CI\rangle'_d$.
Then, $\CI$ is a classical generator of $\CT$ if and only if the triangulated
subcategory $\langle\CI\rangle'_\infty$ of $\CT$ is dense.
By Thomason's characterization of dense subcategories (Theorem \ref{dense} below),
if $\CI$ classically generates $\CT$ and the classes of objects of
$\CI$ generate the abelian group $K_0(\CT)$, then 
$\CT=\langle\CI\rangle'_\infty$.

A similar statement does not hold in general for $d$-step generation,
$d\in\BN$~: take $\CT=D^b((k\times k)\mMod)$, where $k$ is a field.
Let $\CI$ be the full subcategory containing 
$k\times k$ and $k\times 0$ (viewed as complexes concentrated in degree $0$).
Then, $\CT=\langle\CI\rangle$ and
$K_0(\CT)=\BZ\times\BZ$ is generated by the classes of objects of $\CI$,
but $\langle\CI\rangle'$ is not a triangulated subcategory of $\CT$.

Note the necessity of allowing direct summands when $K_0(\CT)$ is not a finitely
generated group (\eg, when $\CT=D^b(X\mcoh)$ and $X$ is an elliptic curve).
\end{rem}

\begin{rem}
It would be interesting to study the ``Krull dimension'' as well.
We say that a
thick subcategory $\CI$ of $\CT$ is irreducible if given two thick
subcategories $\CI_1$ and $\CI_2$ of $\CI$ such that $\CI$ is classically
generated by $\CI_1\ast \CI_2$, then
$\CI_1=\CI$ or $\CI_2=\CI$.
We define
the Krull dimension of $\CT$ as the maximal integer $n$ such that there is
a chain of thick
irreducible subcategories $0\not=\CI_0\subset \CI_1\subset\cdots\subset
\CI_n=\CT$ with $\CI_i\not=\CI_{i+1}$.

By Hopkins-Neeman's Theorem \cite{Nee1}, given a commutative noetherian
ring $A$, the Krull dimension of the category of perfect complexes
of $A$-modules is the Krull dimension of $A$.

By \cite{BeCaRi}, given a finite $p$-group $P$, the Krull dimension of the
stable
category of finite dimensional representations of $P$ over a field of
characteristic $p$ is the $p$-rank of $P$ minus $1$.

\smallskip
Another approach would be to study the maximal possible value for the
transcendence degree of the field of fractions of the center of
$\bigoplus_{i\in\BZ} \Hom(\id_{\CT/\CI},
\id_{\CT/\CI}[i])$, where $\CI$ runs over finitely generated thick subcategories
of $\CT$.
\end{rem}

\begin{rem}
When $\CT$ has finite dimension, every classical generator is a strong
generator.
It would be interesting to study the supremum, over
all classical generators $M$ of $\CT$, of
$\min\{d|\CT=\langle M\rangle_{1+d}\}$.
\end{rem}

\begin{rem}
\label{otherdim}
One can study also, as a dimension,
the minimal integer $d\ge 0$ such that there is
$M$ in $\CT$ with $\CT=\langle\overline{M}\rangle_{d+1}$
or $\CT=\langle\widetilde{M}\rangle_{d+1}$
This is of interest for $D(A)$ and $D(X)$ or $D^b(A)$ and $D^b(X)$.
\end{rem}

\subsubsection{}
\label{devissage}
We often obtain d\'evissages of objects
in the following functorial way (yet another notion of dimension...)~:

Assume there are triangulated functors
$F_i:\CT\to\CT$ with image in $\langle\CI\rangle$ for $1\le i\le d$,
triangulated functors $G_i:\CT\to\CT$ for $0\le i\le d$ with
$G_0=\id$, $G_d=0$
and distinguished triangles
$F_i\to G_i\to G_{i-1}\rightsquigarrow$
for $1\le i\le d$.
Then, $\CT=\langle \CI\rangle_d$.

\subsection{Compact objects}
\subsubsection{}
Let $\CC$ be an additive category.
We say that $\CC$ is {\em cocomplete} if arbitrary direct sums exist in
$\CC$.

An object $C\in\CC$ is {\em compact} if for every set $\CF$ of objects of
$\CC$ such that $\bigoplus_{F\in\CF}F$ exists, then the canonical map
$\bigoplus_{F\in\CF}\Hom(C,F)\to \Hom(C,\bigoplus_{F\in\CF}F)$ is an isomorphism.
We denote by $\CC^c$ the set of compact objects of $\CC$.

A triangulated category  $\CT$ is {\em compactly generated} if it
generated by a set of compact objects. We say that a full triangulated subcategory
$\CI$ of $\CT$ is {\em compactly generated in} $\CT$ if it is generated by a
set of objects of $\CI\cap\CT^c$.

\subsubsection{}
\label{Homcolim}
Let $\CT$ be a triangulated category. Then, $\CT^c$ is a thick subcategory
of $\CT$.

Let $X_0\Rarr{s_0} X_1\Rarr{s_1}\cdots$
be a sequence of objects and maps of $\CT$. If $\bigoplus_{i\ge 0}X_i$ 
exists, then the {\em homotopy colimit} of the sequence,
denoted by $\hocolim X_i$, is a cone of the morphism
$\sum_i \id_{X_i}-s_i:\bigoplus_{i\ge 0}X_i\to\bigoplus_{i\ge 0}X_i$.

We have a canonical map
$$\colim \Hom_\CT(Y,X_i)\to \Hom_\CT(Y,\hocolim X_i)$$
that makes the following diagram commutative
$$\tiny\xymatrix{
&\Hom(Y,\bigoplus X_i)\ar[r] & \Hom(Y,\bigoplus X_i)\ar[r] &
 \Hom(Y,\hocolim X_i)\ar[r] & 
\Hom(Y,\bigoplus X_i[1])\ar[r] & \Hom(Y,\bigoplus X_i[1]) \\
0\ar[r] &\bigoplus \Hom(Y,X_i)\ar[r]\ar[u] & \bigoplus \Hom(Y,X_i)\ar[r]\ar[u] &
 \colim \Hom(Y,X_i)\ar[r]^0\ar@{.>}[u] & 
\bigoplus\Hom(Y,X_i[1])\ar[r]\ar[u] &\bigoplus \Hom(Y,X_i[1])\ar[u] 
}$$

Since the horizontal sequences of the diagram above are exact, we deduce
(cf \eg\ \cite[Lemma 1.5]{Nee2})~:

\begin{lemma}
\label{compactcolim}
The canonical map
$\colim \Hom_\CT(Y,X_i)\to \Hom_\CT(Y,\hocolim X_i)$
is an isomorphism if $Y$ is compact.
\end{lemma}

We now combine the commutation of $\Hom(Y,-)$ with colimits and with
direct sums in the following result (making more precise a classical 
result \cite[Lemma 2.3]{Nee2})~:

\begin{prop}
\label{colimsum}
Let $0=X_0\to X_1\to X_2\to\cdots$ be a directed system in $\CT$,
let $\CF_i$ be a set of compact objects such that $\bigoplus_{C\in \CF_i}C$
exists and let $X_{i-1}\to X_i\to\bigoplus_{C\in\CF_i}C\rightsquigarrow$
be a distinguished triangle, for $i\ge 1$.

Let $Y$ be a compact object and $f:Y\to\hocolim X_i$.
Then, there is an integer $d\ge 1$, a finite subset $\CF_i'$ of $\CF_i$
for $1\le i\le d$ and a commutative diagram
$$\xymatrix{
0=X_0\ar[rr] && X_1\ar[rr]\ar[dl] && X_2\ar[rr]\ar[dl] && X_3\ar[r]\ar[dl]
 & \cdots\ar[r] & X_d\\
& \bigoplus_{\CF_1}C\ar@{~>}[ul] && \bigoplus_{\CF_2}C\ar@{~>}[ul] 
 && \bigoplus_{\CF_3}C\ar@{~>}[ul]\\
& \bigoplus_{\CF'_1}C\ar@{~>}[dl]\ar[u] && \bigoplus_{\CF'_2}C\ar@{~>}[dl]
 \ar[u] && \bigoplus_{\CF'_3}C\ar@{~>}[dl] \ar[u]\\
0=X'_0\ar[rr]\ar[uuu] && X'_1\ar[rr]\ar[ul]\ar[uuu] && X'_2\ar[rr]\ar[ul]
 \ar[uuu] && X'_3\ar[r]\ar[ul]\ar[uuu] & \cdots \ar[r] & X'_d\ar[uuu]_h
}$$
such that $f$ factors through $X'_d\Rarr{h} X_d\Rarr{\can}\hocolim X_i$.
\end{prop}

\begin{proof}
By Lemma \ref{compactcolim}, there is $d\ge 1$ such that $f$ factors through
the canonical map $X_d\to \hocolim X_i$.
We proceed now by induction on $d$.
The composite map $Y\to X_d\to \bigoplus_{C\in\CF_d}C$ factors through
the sum indexed by a finite subset $\CF'_d$ of $\CF_d$. 
Let $Z$ be the cocone of the corresponding map $Y\to \bigoplus_{C\in\CF'_d}C$
and $X''_d$ the cocone of the composite map
$\bigoplus_{C\in\CF'_d}C\to \bigoplus_{C\in\CF_d}C\to X_{d-1}[1]$. The composite
map $X''_d\to \bigoplus_{C\in\CF'_d}C\to \bigoplus_{C\in\CF_d}C$ factors
through $X_d$. The map $Y\to X_d$ factors through $X''_d$ and
the composite map $Z\to Y\to X''_d$ factors through $X_{d-1}$.
Summarizing, we have a commutative diagram
$$\xymatrix{
X_{d-1}\ar[r] & X_d\ar[r] & \bigoplus_{C\in\CF_d}C \ar@{~>}[r] & \\
X_{d-1}\ar@{=}[u]\ar[r]   & X''_d\ar[r]\ar@{.>}[u] & \bigoplus_{C\in\CF'_d}C 
  \ar@{~>}[r]\ar[u] & \\
Z        \ar@{.>}[u]\ar[r] & Y\ar@/^2pc/[uu]\ar@{.>}[r]\ar@{.>}[u]_a &
 \bigoplus_{C\in\CF'_d}C \ar@{~>}[r]\ar@{=}[u] & 
}$$
By induction, we have already a commutative diagram as in the
proposition for the corresponding map $Z\to X_{d-1}$. We define
now $X'_d$ to be the cocone of the composite map $\bigoplus_{C\in\CF'_d}C
\to Z[1]\to X'_{d-1}[1]$. There is a commutative diagram
$$\xymatrix{
X_{d-1}\ar[r] & X''_d\ar[r] & \bigoplus_{C\in\CF'_d}C\ar[rr]\ar[dr] &&
 X_{d-1}[1]  \\
&&&Z[1]\ar[ur]\ar[dr] \\
X'_{d-1}\ar[r]\ar[uu] & X'_d\ar[r] & \bigoplus_{C\in\CF'_d}C\ar[rr]\ar@{=}[uu]
 && X'_{d-1}[1]\ar[uu]\\
Z\ar[r] \ar[u] & Y\ar@/^2pc/[uuu]^a\ar@{.>}[u] \ar[r] &
 \bigoplus_{C\in\CF'_d}C\ar[rr]\ar@{=}[u] && Z[1]\ar[u]
}$$
The composite map
$Z\to Y\to X''_d$ factors through $X_{d-1}$, hence through $X'_{d-1}$.
It follows that $a$ factors through $X'_d$ and we are done.
\end{proof}

We deduce the following descent result \cite[Proposition 2.2.4]{BoVdB}~:

\begin{cor}
\label{intercompact}
Let $\CI$ be a subcategory of $\CT^c$ and let $d\in\BN\cup\{\infty\}$. Then,
$\CT^c\cap \langle\overline{\CI}\rangle_d=\langle\CI\rangle_d$.
\end{cor}

\begin{proof}
Let $Y$ be a compact object and $f:Y\to X_d$ be a split injection where $X_d$
is obtained by taking a $d$-fold extension of objects of
$\langle\overline{\CI}\rangle$.
Proposition \ref{colimsum} shows that $f$ factors through an object
$X'_d\in\langle\CI\rangle_d$ and we obtain a split injection $Y\to X'_d$.
\end{proof}

\subsection{Relation with dg algebras}
\label{dgend}
Following Keller, we say that a triangulated category $\CT$ is {\em algebraic}
if it is the stable category of
a Frobenius exact category \cite[Chapter 5, \S 2.6]{GeMa}
(for example, $\CT$ can be
the derived category of an abelian category).

\smallskip
Recall the construction of \cite[\S 4.3]{Ke}.
Let $\CT=\CE\mstab$ be the stable category of a Frobenius exact category $\CE$.
Let $\CE'$ be the category of acyclic complexes of
projective objects of $\CE$ and $Z^0:\CE'\to\CE\mstab$ be the functor
that sends $C$ to $\coker d_C^{-1}$.

Given $X$ and $Y$ two complexes of objects of $\CE$, we denote by
$\Hom^\bullet(X,Y)$ the total $\Hom$ complex (\ie,
$\Hom^\bullet(X,Y)^i=\prod_{j\in\BZ} \Hom_{\CE}(X^j,Y^{i+j})$).

Let $M\in\CE\mstab$ and $M'\in\CE'$ with $Z^0(M')\iso M$.
Let $A=\End^\bullet(M')$ be the dg algebra
of endomorphisms of $M'$.
The functor $\Hom^\bullet(M',-):\CE'\to D(A)$
factors through $Z^0$ and induces a triangulated functor
$R\Hom^\bullet(M,-):\CE\mstab\to D(A)$.
That functor restricts to an equivalence $\langle M\rangle_\infty\iso A\mperf$.
In particular, if
$M$ is a classical generator of $\CT$, then we get the equivalence
$\CT\iso A\mperf$.

So,

\begin{prop}
Let $\CT$ be an algebraic triangulated category. Then, $\CT$ is classically finitely
generated if and only if it is equivalent to the category of perfect
complexes over a dg algebra.
\end{prop}

This should be compared with the following result.

Assume now $\CE$ is a cocomplete Frobenius category (\ie, all
direct sums exist and are exact).
If $M$ is compact, then $R\Hom^\bullet(M,-)$ restricts to
an equivalence between the smallest full triangulated subcategory of
$\CT$ containing $M$ and closed under direct sums and $D(A)$
(cf Theorem \ref{equivgen} (2) and Corollary \ref{genDG} below).
So, using Theorem \ref{equivgen} (2) below, we deduce \cite[Theorem 4.3]{Ke}~:

\begin{thm}
Let $\CE$ be a cocomplete Frobenius category and $\CT=\CE\mstab$.
Then, $\CT$ has a compact generator if and only if it is equivalent
to the derived category of a dg algebra.
\end{thm}

\section{Finiteness conditions and representability}
\label{secfin}
\subsection{Finiteness for cohomological functors}
\label{subsecfin}
We introduce a class of ``locally finitely presented'' cohomological
functors that includes the representable functors, inspired by Brown's
representability Theorem. It extends the class of locally finite functors,
of interest only for $\Ext$-finite triangulated categories.

\subsubsection{}
Let $k$ be a commutative ring.

Let $\CT$ be a $k$-linear triangulated category.
Let
$H:\CT^\opp\to k\mMOD$ be a 
($k$-linear) functor. We say that $H$ is
{\em cohomological} if for every distinguished triangle
$X\Rarr{f} Y\Rarr{g} Z\rightsquigarrow$, then the associated sequence
$H(Z)\Rarr{H(g)}H(Y)\Rarr{H(f)}H(X)$ is exact.

For $C\in \CT$, we denote by $h_C$ the cohomological
functor $\Hom_{\CT}(-,C):\CT^\opp\to k\mMOD$.

\smallskip
We will repeatedly use Yoneda's Lemma~:

\begin{lemma}
Let $X\in \CT$ and $H:\CT^\opp\to k\mMOD$ a functor. Then, the canonical map
$\Hom(h_C,H)\to H(C), f\mapsto f(C)(\id_C)$ is an isomorphism.
\end{lemma}

\smallskip
Let $H:\CT^\opp\to k\mMOD$ be a functor.
We say that $H$ is
\begin{itemize}
\item
{\em locally bounded} (resp. {\em bounded above}, resp. {\em bounded below}) if for
every $X\in\CT$, we have
$H(X[i])=0$ for $|i|\gg 0$ (resp. for $i\ll 0$, resp. for $i\gg 0$)
\item
{\em locally finitely generated} if for every $X\in\CT$, there is 
$D\in\CT$ and $\alpha:h_D\to H$ such that $\alpha(X[i])$ is surjective
for all $i$.
\item
{\em locally finitely presented}
if it is locally finitely generated and
the kernel of any map $h_E\to H$ is locally finitely generated.
\end{itemize}

Let $X\in\CT$. We introduce two conditions~:
\begin{itemize}
\item[(a)] there is 
$D\in\CT$ and $\alpha:h_D\to H$ such that $\alpha(X[i])$ is surjective
for all $i$
\item[(b)] for every $\beta:h_E\to H$, there is $f:F\to E$ such that
$\beta h_f=0$ and
$h_F(X[i])\Rarr{h_f} h_E(X[i])\Rarr{\beta} H(X[i])$ is an exact sequence
for all $i$.
\end{itemize}

Note that $H$ is locally finitely presented if and only if
for every $X\in\CT$, then conditions (a) and (b) are fulfilled.

\begin{lemma}
\label{presislocfipres}
For $C\in\CT$, then $h_C$ is locally finitely presented.
\end{lemma}

\begin{proof}
We take $D=C$ and $\alpha=\id$ for condition (a). For (b), a map
$\beta:h_E\to h_C$ comes from a map $g:E\to C$ and
we pick a distinguished triangle
$F\Rarr{f}E\Rarr{g}C\rightsquigarrow$.
\end{proof}

\begin{prop}
\label{lfpstable}
Let $H_0\to H_1\to H\to H_2\to H_3$ be an exact sequence of
functors $\CT^\opp\to k\mMOD$. 

If $H_1$ and $H_2$ are locally finitely generated and $H_3$ is
locally finitely presented, then $H$ is locally finitely generated.

If $H_0$ is locally finitely generated and
$H_1$, $H_2$ and $H_3$ are locally finitely presented,
then $H$ is locally finitely presented.
\end{prop}

\begin{proof}
Let us name the maps~:
$H_0\Rarr{t_0}H_1\Rarr{t_1}H\Rarr{t_2}H_2\Rarr{t_3}H_3$.
Let $X\in\CT$.

Let $\alpha_2:h_{D_2}\to H_2$ as in (a).
Let $\beta_3=t_3\alpha_2:h_{D_2}\to H_3$.
Let $f_3:E\to D_2$ as in (b).
Since $H(E)\to H_2(E)\to H_3(E)$ is exact, the composite
map $\alpha_2 h_{f_3}:h_E\to H_2$ factors as
$h_E\Rarr{\gamma} H\Rarr{t_2}H_2$.
Let $\alpha_1:h_{D_1}\to H_1$ as in (a).

Let $a:h_X\to H$. The composite $t_2a:h_X\to H_2$ factors as
$t_2a:h_X\Rarr{b} h_{D_2}\Rarr{\alpha_2}H_2$.
The composition $t_3(t_2a):h_X\to H_3$ is zero,
hence $b$ factors as $b:h_X\Rarr{c} h_E\Rarr{h_{f_3}}h_{D_2}$.
Now, we have
$t_2\gamma c=\alpha_2 h_{f_3} c=\alpha_2 b=t_2a$. Since
the composite $t_2(a-\gamma c):h_X\to H_2$ is zero,  it follows that
$a-\gamma c$ factors as $h_X\Rarr{a_1} H_1\Rarr{t_1} H$. Now, $a_1$ factors
through $\alpha_1$.
So, we have shown that $a$ factors through
$\gamma+t_1\alpha_1:h_E\oplus h_{D_1}\to H$, hence $H$ satisfies (a).

\smallskip
Let $\alpha_0:h_{D_0}\to H_0$ as in (a).
Let $\beta':h_E\to \ker t_2$.
Then, there is $\beta_1:h_E\to H_1$ such that $\beta'=t_1\beta_1$.
Since $H_1$ is locally finitely presented, there are
$u:h_F\to h_{D_0}$ and $v:h_F\to h_E$ such that
$(\beta_1+t_0\alpha_0)(u-v)=0$ and
$h_F(X[i])\Rarr{u-v} h_E(X[i])\oplus h_{D_0}(X[i])
\xrightarrow{\beta_1+t_0\alpha_0} H_1(X[i])$ is
exact for every $i$. Summarizing, we have a commutative diagram
$$\xymatrix{
h_F\ar[d]_u\ar[dr]^v \\
h_{D_0}\ar[d]_{\alpha_0}
\ar[dr] & h_E\ar[d]^{\beta_1}\ar[dr]^{\beta'} \\
H_0\ar[r]_{t_0} & H_1 \ar[r]_{t_1} & \ker t_2
}$$
It follows that $\beta v=0$ and
$h_F(X[i])\Rarr{v} h_E(X[i])\Rarr{\beta}(\ker t_2)(X[i])$ is exact
for every $i$, hence
$\ker t_2$ satisfies (b).

Let now $\beta:h_E\to H$. Let $G=\ker\beta$ and $G_2=\ker(t_2\beta)$.
Now, we have exact sequences
$0\to G\to G_2\to \ker t_2$ and
$0\to G_2\to h_E\to H_2$. The first part of the Proposition together
with Lemma \ref{presislocfipres} shows
that $G$ is finitely generated. Consequently, $H$ is locally finitely
generated.
\end{proof}

\subsubsection{}
We will now study conditions (a) and (b) in the definition of
locally finitely presented functors.

\begin{lemma}
\label{surj}
Let $H:\CT^\opp\to k\mMOD$ be a $k$-linear functor and $X\in\CT$.
\begin{itemize}
\item
Let $\beta_r:h_{E_r}\to H$ for $r\in\{1,2\}$ such that (b) holds for
$\beta=\beta_1+\beta_2:h_{E_1\oplus E_2}\to H$. Then, (b)
holds for $\beta_1$ and $\beta_2$.
\item
Assume (a) holds. If
(b) holds for those $\beta:h_E\to H$ such that $\beta(X[i])$ is
surjective for all $i$, then (b) holds for all $\beta$.
\end{itemize}
\end{lemma}

\begin{proof}
Let $E=E_1\oplus E_2$. Denote by $i_r:E_r\to E$ and $p_r:E\to E_r$ the 
injections and projections.
There is $f:F\to E$ such that $\beta h_f=0$ and
$h_F(X[i])\Rarr{h_f} h_E(X[i])\Rarr{\beta} H(X[i])$ is an
exact sequence for all $i$.

Fix a distinguished triangle
$F_1\Rarr{f'_1}F\Rarr{p_2f}E_2\rightsquigarrow$ and let
$f_1=p_1ff'_1:F_1\to E_1$. We have $\beta_1 h_{f_1}=0$ since
$\beta_1h_{p_1f}=-\beta_2h_{p_2f}$.

For all $i$, the horizontal sequences
and the middle vertical sequence in the following commutative
diagram are exact
$$\xymatrix{
&h_{F_1}(X[i])\ar[r]^{h_{f'_1}}\ar[d]_{h_{f_1}} & h_{F}(X[i])\ar[r]^{h_{p_2f}}
 \ar[d]^{h_f} & h_{E_2}(X[i])\ar@{=}[d] \\
0\ar[r] & h_{E_1}(X[i])\ar[r]^{h_{i_1}}\ar[d]_{\beta_1} &
 h_{E}(X[i])\ar[r]^{h_{p_2}}\ar[d]^{\beta} & h_{E_2}(X[i]) \ar[r] & 0 \\
&H(X[i]) \ar@{=}[r] & H(X[i])\ar[d] \\
&&0
}$$
hence the left vertical sequence is exact as well.

\smallskip
Let us now prove the second part of the Lemma.
Let $\beta:h_E\to H$. Since (a) holds, there is $D\in\CT$ and
$\alpha:h_D\to H$ such that $\alpha(X[i])$ is surjective
for all $i$. Let $E'=D\oplus E$ and $\beta'=\alpha+\beta:h_{E'}\to H$.
Then, (b) holds for $\beta'$, hence it holds for $\beta$ by the first part
of the Lemma.
\end{proof}

\begin{rem}
For the representability Theorem (cf Lemma \ref{approx}),
only the surjective case of (b) is needed,
but the previous Lemma shows that this implies that (b) holds in general.
\end{rem}

\begin{lemma}
\label{Xthick}
Let $H:\CT^\opp\to k\mMOD$ be a cohomological functor.

The full subcategory of $X$ in $\CT$ such that (a) and (b) hold
is a thick triangulated subcategory of $\CT$.

In particular, if $X$ is a classical generator for $\CT$ and (a), (b)
hold, then $H$ is locally finitely presented.
\end{lemma}

\begin{proof}
Let $\CI$ be the full subcategory of those $X$ such that (a) and (b) hold.
It is clear that $\CI$ is closed under shifts and under taking direct summands.
So, we are left with proving that $\CI$ is stable under extensions.

%
%
Let $X_1\Rarr{u} X \to X_2\rightsquigarrow$ be a distinguished triangle in
$\CT$ with $X_1,X_2\in\CI$.
Pick $D_r\in\CT$ and 
$\alpha_r:h_{D_r}\to H$ such that $\alpha_r(X_r[i])$ is surjective for all $i$.
Put $E=D_1\oplus D_2$ and $\beta=\alpha_1+\alpha_2:h_E\to H$.
There is $F_r\in\CT$ and $f_r:F_r\to E$ such that
$\beta h_{f_r}=0$ and
$h_{F_r}(X_r[i])\Rarr{h_f} h_E(X_r[i])\Rarr{\beta} H(X_r[i])$ is an exact
sequence for all $i$.
Put $F=F_1\oplus F_2$ and $f=f_1+f_2:F\to E$. Let $F\Rarr{f}E\Rarr{t} E'
\rightsquigarrow$ be a distinguished triangle.
We have an exact sequence $H(E')\to H(E)\to H(F)$.
The image in $H(F)$ of the element of $H(E)$ corresponding to $\beta$
is $0$, since $\beta h_f=0$. Hence, $\beta$ factors as
$h_E\Rarr{h_t}h_{E'}\Rarr{\gamma}H$.
Let $D=E\oplus E'$ and $\alpha=\beta+\gamma:h_D\to H$.

Let $a:h_X\to H$. Then, there is a commutative diagram where the top
horizontal sequence is exact
$$\xymatrix{
h_{X_2[-1]} \ar[r]\ar@{.>}[d] & h_{X_1}\ar[r]^{h_u}\ar@{.>}[d]_c & h_X \ar[d]^a\\
h_F\ar[r]_{h_f} & h_E\ar[r]_\beta & H
}$$
The composite $h_{X_2[-1]}\to h_{X_1}\Rarr{c} h_E \Rarr{h_t} h_{E'}$ is zero, hence
$h_t c:h_{X_1}\to h_{E'}$ factors as $h_{X_1}\Rarr{h_u}h_X\Rarr{b} h_{E'}$.
We have $ah_u=\beta c=\gamma h_tc=\gamma bh_u$, hence 
the composite $h_{X_1}\to h_X\Rarr{a'}H$ is zero, where
$a'=a-\gamma b$.
So, $a'$ factors through a map $h_{X_2}\to H$. Such a map
factors through $\beta$, hence $a'$ factors through $\beta$ and $a$ factors
through $\alpha$.
The same conclusion holds for $a$ replaced by any map
$h_{X[i]}\to H$ for some $i\in\BZ$.
So, every map
$h_{X[i]}\to H$ factors through $\alpha$, \ie, (a) holds for $X$.

\smallskip
Consider now a map $\beta':h_{E'}\to H$. Let
$\beta'':h_{E''}\to H$ such that
$\beta''(X_1[i])$ is surjective for all $i$. 
Let $\beta=\beta'+\beta'':h_E\to H$, where $E=E'\oplus E''$.
In order to prove that $\beta'$ satisfies (b), it suffices to prove
that $\beta$ satisfies (b), thanks to Lemma \ref{surj}.

There is $F_1\in\CT$ and
$f_1:F_1\to E$ such that $\beta h_{f_1}=0$ and
$h_{F_1}(X_1[i])\Rarr{h_{f_1}} h_E(X_1[i])\Rarr{\beta} H(X_1[i])$ is an exact
sequence for all $i$.
Let $E_1$ be the cone of $f_1$. As in the discussion above, $\beta$
factors through a map $\gamma:h_{E_1}\to H$.
Let $F_2\in\CT$ and $f_2:F_2\to E_1$ such that $\gamma h_{f_2}=0$ and
$h_{F_2}(X_2[i])\Rarr{h_{f_2}} h_{E_1}(X_2[i])\Rarr{\gamma} H(X_2[i])$
is an exact sequence for all $i$.
Let $F$ be the cocone of the sum map $E\oplus F_2\to E_1$.
The composition $h_F\to h_{E}\Rarr{\beta} H$ is zero.
We have a commutative diagram
$$\xymatrix{
&& h_F\ar[rr]\ar[d] && h_{F_2}\ar[d]^{h_{f_2}} \\
h_{F_1}\ar[rr]^{h_{f_1}} && h_E\ar[rd]_\beta \ar[rr] && h_{E_1}\ar[dl]^\gamma \\
&&&H
}$$
In the diagram, the square is homotopy cartesian, \ie, given $Y\in\CT$ and
$u:Y\to E$, $v:Y\to F_2$ such that the compositions $Y\Rarr{u}E\to E_1$ and
$Y\Rarr{v}F_2\to E_1$ are equal, then there is $w:Y\to F$ such that
$u$ is the composition $Y\Rarr{w}F\to E$ and $v$ the composition
$Y\Rarr{w}F\to F_2$.

Let $a:h_X\to h_E$ such that $\beta a=0$.
The composite $h_{X_1}\to h_X\Rarr{a}h_E$ factors through $h_{F_1}$.
It follows that the composition $h_{X_1}\to h_X\Rarr{a}h_E\to h_{E_1}$ is
zero. Hence, the composite $h_X\Rarr{a}h_E\to h_{E_1}$ factors
through a map $b:h_{X_2}\to h_{E_1}$.
The composite $b':h_{X_2}\Rarr{b}h_{E_1}\to H$ factors through a map
$c:h_{X_1[1]}\to H$, since $h_X\to h_{X_2}\Rarr{b'}H$ is zero.
Now, $c$ factors as $h_{X_1[1]}\Rarr{d}h_{E_1}\Rarr{\gamma} H$.
Summarizing, we have a diagram all of whose squares and triangles but
the one marked ``$\not=$'' are commutative and where the 
horizontal sequences are exact
$$\xymatrix{
h_{X_1} \ar[r]\ar@{.>}[d] & h_X\ar[r]\ar[d]_a & h_{X_2}\ar@{.>}[d]_b\ar[r] &
h_{X_1[1]}\ar@{.>}[dl]^d_{\not=} \ar@{.>}@/^2pc/[ddll]^c \\
h_{F_1} \ar[r] & h_E \ar[r]\ar[d]_\beta & h_{E_1}\ar[dl]^\gamma \\
& H
}$$
Let $d'$ be the composition $h_{X_2}\to h_{X_1[1]}\Rarr{d}h_{E_1}$.
Then, the composition $h_{X_2}\Rarr{b-d'} h_{E_1}\Rarr{\gamma}H$ is zero.
The map $b-d'$ factors as $h_{X_2}\Rarr{d''}h_{F_2}\Rarr{h_{f_2}}h_{E_1}$.
It follows that $h_X\Rarr{a}h_E\to h_{E_1}$ factors as
$h_X\to h_{X_2}\Rarr{d''}h_{F_2}\Rarr{h_{f_2}}h_{E_1}$. Using the homotopy
cartesian square above, we deduce that
$a$ factors through $a':h_X\to h_F$. So, the sequence
$h_F(X)\to h_E(X)\to H(X)$ is exact. The same holds for all $i$, hence
(b) holds for $X$.
\end{proof}

\begin{rem}
\label{inequalities}
All the results concerning locally finitely generated and presented 
functors above remain valid if we replace the conditions ``given $X\in\CT$,
a certain statement is true for all $i\in\BZ$''
by ``given $X\in\CT$ and $a\in \BZ$, a certain statement is true for 
$i\ge a$'' (or ``$i\le a$'' or ``$i=a$'') in
(a) and (b). Cf for example Proposition \ref{pseudocoherent}.
\end{rem}

\subsubsection{}

\begin{prop}
Let $H:\CT^\opp\to k\mMOD$ be a cohomological functor
and $X$ be a classical generator for $\CT$. Then, $H$ is locally finitely generated
if and only if $\bigoplus_i H(X[i])$ is a finitely generated $\End^*(X)$-module.
\end{prop}

\begin{proof}
Assume first $\bigoplus_i H(X[i])$ is a finitely generated $\End^*(X)$-module.
Let $f_r\in \Hom(h_{X[n_r]},H)$ be a finite set of elements such that
the $f_r(X[n_r])(\id_{X[n_r]})$ generate $\bigoplus_i H(X[i])$ as
an $\End^*(X)$-module.
Let $D=\bigoplus_r X[n_r]$ and $f=\sum_r f_r:h_D\to H$.
Then, $f(X[i])$ is surjective for every $i$, \ie, condition (a) is satisfied.

Conversely, assume (a) is satisfied.
\end{proof}

\subsubsection{}
Assume $k$ is noetherian.
We say that $\CT$ is {\em $\Ext$-finite} if
$\bigoplus_i \Hom(X,Y[i])$ is a finitely generated $k$-module, for
every $X,Y\in\CT$.

Assume now $\CT$ is $\Ext$-finite and let $H:\CT^\opp\to k\mMOD$ be
a functor. We say that $H$ is 
{\em locally finite} if for every $X\in\CT$, the $k$-module
$\bigoplus_i H(X[i])$ is finitely generated.

\begin{prop}
\label{locfin}
Let $H$ be a locally finite functor. Then,
$H$ is locally bounded and locally finitely presented.
\end{prop}

\begin{proof}
It is clear that $H$ is locally bounded.
Let $X\in\CT$. Let $I_i$ be a minimal (finite) family of
generators of $H(X[i])$ as a $k$-module. We have
$I_i=\emptyset$ for almost all $i$, since $H$ is locally bounded. Put
$D=\bigoplus_i X[i]\otimes_k k^{I_i}$ and let
$\alpha:h_D\to H$ be the canonical map. The map
$\alpha(X[i])$ is surjective for all $i$. So,
every locally finite functor is locally finitely generated.

Let now $\beta:h_E\to H$.
Let $G=\ker\beta$. Since $\CT$ is $\Ext$-finite, $G$ is again locally finite,
hence locally finitely generated.
\end{proof}

The results on finitely generated and presented functors discussed above
have counterparts for locally bounded functors, the proofs being trivial
in this case.

\begin{prop}
\label{lb}
Let $H_1\to H\to H_2$ be an exact sequence of functors $\CT^\circ\to k\mMOD$.
If $H_1$ and $H_2$ are locally bounded (resp. bounded above, resp. bounded
below), then $H$ is locally bounded (resp. bounded above, resp. bounded
below).

Let $H:\CT^\circ\to k\mMOD$ be a cohomological functor. Then,
the full subcategory of $X\in\CT$ such that $H(X[i])=0$ for
$|i|\gg 0$ (resp. $i\ll 0$, resp. $i\gg 0$) is a thick subcategory.
\end{prop}

\subsection{Locally finitely presented functors}
\label{seclfp}
\subsubsection{}
Let us start with some remarks on cohomological functors.

\smallskip
Given $0\to H_1\to H_2\to H_3\to 0$ an exact sequence of functors
$\CT^\opp\to k\mMOD$, if two of the functors amongst the $H_i$'s are
cohomological, then the third one is cohomological as well.
The category of cohomological functors $\CT^\opp\to k\mMOD$ is
closed under direct sums.

Given $H_1\to H_2\to\cdots$ a directed system of cohomological
functors $\CT^\opp\to k\mMOD$, we have an exact sequence
$0\to \bigoplus H_i\to \bigoplus H_i\to \colim H_i\to 0$. This shows that
$\colim H_i$ is a cohomological functor.

\begin{lemma}
\label{compo0}
Let $H_1,\ldots,H_{n+1}$ be cohomological functors on $\CT$ and
$f_i:H_i\to H_{i+1}$ for $1\le i\le n$.
Let $\CI_i$ be a subcategory of $\CT$ on which
$f_i$ vanishes. Then,
$f_n\cdots f_1$ vanishes on $\CI_1\diamond\cdots\diamond\CI_n$.
\end{lemma}

\begin{proof}
Note first that if a morphism between cohomological functors vanishes
on a subcategory $\CI$, then it vanishes on $\langle\CI\rangle$.

By induction, it is enough to prove the Lemma for $n=2$.
Let $X_1\to X\to X_2\rightsquigarrow$ be a distinguished triangle with
$X_i\in \CI_i$. The map $f_1(X)$ factors through $H_2(X_2)$, \ie, 
we have a commutative diagram with exact horizontal sequences
$$\xymatrix{
H_1(X_2)\ar[r]\ar[d] & H_1(X)\ar[r]\ar[d]\ar@{.>}[dl] & H_1(X_1) \ar[d]^0 \\
H_2(X_2)\ar[r]\ar[d]_0 & H_2(X)\ar[r]\ar[d] & H_2(X_1)   \ar[d] \\
H_3(X_2)\ar[r]       & H_3(X)\ar[r]       & H_3(X_1)    
}$$
This shows that $f_2f_1(X)=0$.
\end{proof}

\begin{rem}
Let $M\in\CT$ be a complete classical generator.
Let $f:\bigoplus_i
\Hom(\id_\CT,\id_\CT[i])\to \bigoplus_i \Hom(M,M[i])$ be the canonical map.
Let $\zeta\in\ker f$. It follows from Lemma \ref{compo0} that
$\zeta$ is locally nilpotent.
If $\CT=\langle \overline{M}\rangle_d$, then $(\ker f)^d=0$.
\end{rem}

\subsubsection{}
In this part, we study convergence conditions on directed systems. This
builds on \cite[\S 2.3]{BoVdB}.

Let $V_1\Rarr{f_1}V_2\Rarr{f_2}\cdots$ be a system of abelian groups.
We say that the system $(V_i)$
 is {\em almost constant} if one of the
following equivalent conditions is satisfied~:
\begin{itemize}
\item
$V_i=\im f_{i-1}\cdots f_2 f_1+\ker f_i$ and
$\ker f_{i+r}\cdots f_i=\ker f_i$ for any $r\ge 0$ and $i\ge 1$.
\item
Denote by $\alpha_i:V_i\to V=\colim V_i$ the canonical map. Then,
$\alpha_i$ induces an isomorphism $V_i/\ker f_i\iso V$.
\end{itemize}

\medskip

Let $\CT$ be a triangulated category and $\CI$ a subcategory of $\CT$.
Let $H_1\to H_2\to\cdots$ be a directed system of functors
$\CT^\opp\to k\mMOD$.
We say that $(H_i)_{i\ge 1}$ is
{\em almost constant} on $\CI$ if for every $X\in\CI$, the system
$H_1(X)\to H_2(X)\to\cdots$ is almost constant.

\smallskip
Given $1\le r_1<r_2<\cdots$, we denote by $(H_{r_i})$ the system
$H_{r_1}\xrightarrow{f_{r_2-1}\cdots f_{r_1+1}f_{r_1}}H_{r_2}
\xrightarrow{f_{r_3-1}\cdots f_{r_2+1}f_{r_2}}H_{r_3}\to\cdots$.

\begin{prop}
\label{propconv}
Let $(H_i)_{i\ge 1}$ be a directed system of cohomological
functors on $\CT$.
\begin{itemize}
\item[(i)]
If $(H_i)_{i\ge 1}$ is almost constant on $\CI_1,\CI_2,\ldots,\CI_n$,
then, for any $r>0$, the system $(H_{ni+r})_{i\ge 0}$ is
almost constant on $\CI_1\diamond\cdots\diamond\CI_n$.
\end{itemize}
Assume now 
$(H_i)_{i\ge 1}$ is almost constant on $\CI$. Then,
\begin{itemize}
\item[(ii)]
$(H_i)_{i\ge 1}$ is almost constant on $\add(\CI)$. If in addition the
functors $H_i$ commute with products, then 
$(H_i)_{i\ge 1}$ is almost constant on $\overline{\add}(\CI)$.
\item[(iii)]
$(H_{ir+s})_{i\ge 0}$ is almost constant on $\CI$ for any $r,s>0$.
\item[(iv)]
the canonical map
$H_{n+1}\to \colim H_i$ is a split surjection, when the functors are restricted
to $\langle \CI\rangle_n$.
\end{itemize}
\end{prop}

\begin{proof}
Let $H=\colim H_i$ and let $K_i=\ker (H_i\to H)$.
Take $\CI$ and $\CI'$ such that $(H_i)$ is almost constant on $\CI$ and
$\CI'$.
Let $I\to J\to I'$ be a distinguished triangle with $I\in\CI$ and $I'\in\CI'$.

Given $i\ge 1$, we have a commutative diagram with exact rows and columns
$$\xymatrix{
& & H_i(I') \ar[d]\ar[r] & H(I') \ar[d] \ar[r] & 0 \\
& & H_i(J) \ar[d]\ar[r] & H(J) \ar[d] & \\
0\ar[r] & K_i(I) \ar[r]\ar[d] &  H_i(I)\ar[d]\ar[r] & H(I) \ar[d] \ar[r] & 0 \\
0\ar[r] & K_i(I'[-1]) \ar[r]\ar[d] &  H_i(I'[-1])\ar[d]\ar[r] & H(I'[-1]) \ar[r] & 0 \\
0\ar[r] & K_i(J[-1]) \ar[r] &  H_i(J[-1])
}$$
This shows that $H_i(J)\to H(J)$ is onto.
By induction, we deduce that $H_i(X)\to H(X)$ is onto for any $i\ge 1$ and
any $X\in \CI_1\diamond\cdots\diamond\CI_n$.
It follows from Lemma \ref{compo0} that
the composition
$K_i\Rarr{f_i} K_{i+1}\to\cdots\to K_{i+n}$ vanishes on
$\CI_1\diamond\cdots\diamond\CI_n$. We deduce that (i) holds.

\smallskip
The assertions (ii) and (iii) are clear.

\smallskip
By (i), it is enough to prove (iv) for $n=1$.
The map $f_1:H_1\to H_2$ factors through $H_1/K_1$ as
$\bar{f}_1:H_1/K_1\to H_2$. We have a commutative diagram
$$\xymatrix{
H_1/K_1\ar[rr]^{\bar{f}_1}\ar[dr] && H_2 \ar[dl] \\
& H}$$
When restricted to $\CI$, the canonical map
$H_1/K_1\to H$ is an isomorphism, hence the canonical map
$H_2\to H$ is a split surjection. This proves (iv).
\end{proof}

We say that a direct system $(A_1\Rarr{f_1} A_2\Rarr{f_2}\cdots)$
of objects of $\CT$ is
{\em almost constant} on $\CI$ if the system
$(h_{A_i})$ is almost constant on $\CI$.

\subsubsection{}

We study now approximations of locally finitely presented functors.

\begin{lemma}
\label{approx}
Let $\CT$ be a triangulated category and $G\in\CT$. Let $H$ be
a locally finitely presented cohomological functor.
Then, there is a directed system
$A_1\Rarr{f_1} A_2\Rarr{f_2}\cdots$ in $\CT$ that is almost constant
on $\{G[i]\}_{i\in\BZ}$
and a map
$\colim h_{A_i}\to H$ that is an isomorphism on $\langle G\rangle_\infty$.
\end{lemma}

\begin{proof}
%
Since $H$ is locally finitely presented, there is $A_1\in\CT$ and
$\alpha_1:h_{A_1}\to H$ such that $\alpha_1(G[r])$ is onto for
all $r$.

We now construct the system by induction on $i$.
Assume
$A_1\Rarr{f_1}A_2\Rarr{f_2}\cdots \Rarr{f_{i-1}}A_i$ and
$\alpha_1,\ldots,\alpha_i$
have been constructed.

Since $H$ is locally finitely presented, there is $g:B\to A_i$ with
$\im h_g(G[r])=\ker \alpha_i(G[r])$ for all $r$ and with
$h_g\alpha_i=0$. Let $B\Rarr{g} A_i\Rarr{f_i} A_{i+1}\rightsquigarrow$
be a distinguished triangle. We have an exact sequence
$h_B\Rarr{h_g} h_{A_i}\Rarr{h_{f_i}} h_{A_{i+1}}$, hence, there is
$\alpha_{i+1}:h_{A_{i+1}}\to H$
with $\alpha_i=\alpha_{i+1}f_i$. We have a surjection
$h_g(G[r]):h_B(G[i])\to \ker\alpha_i(G[r])$, hence
$\ker\alpha_i(G[r])\subseteq \ker f_i(G[r])$.
So, the system is almost constant on $\{G[i]\}_{i\in\BZ}$.
It follows from Proposition \ref{propconv} (iv) that
the canonical map $H\to\colim h_{A_i}$ is an isomorphism on
$\langle G\rangle_\infty$.
\end{proof}

\begin{prop}
\label{approxlfp}
Let $\CT$ be a triangulated category classically generated by an object $G$.
Let $H$ be a cohomological functor.
Then, $H$ is locally finitely presented if and only if there is a 
 directed system
$A_1\Rarr{f_1} A_2\Rarr{f_2}\cdots$ in $\CT$ that is almost
constant on $\{G[i]\}_{i\in\BZ}$ and an isomorphism
$\colim h_{A_i}\iso H$.
\end{prop}

\begin{proof}
The first implication is given by Lemma \ref{approx}.
Let us now show the converse.

Since $\CT$ is classically generated by $G$, it is enough to
show conditions (a) and (b) for $X=G$ (cf Lemma \ref{Xthick}).
Condition (a) is obtained with $\alpha_1:h_{A_1}\to H$.
Fix now $\beta:h_E\to H$. There is an integer $i$ such that
$E\in\langle G\rangle_i$. By Proposition \ref{propconv} (iii) and (iv),
the restriction of
$\alpha_{i+1}$ to $\langle G\rangle_i$ has a right inverse $\rho$.
We obtain a map $\rho\beta$ between the functors $h_E$ and $h_{A_{i+1}}$ 
restricted to $\langle G\rangle_i$.
It comes from a map
$f:E\to A_{i+1}$. Let $F$ be the cocone of $f$.
The kernel of 
$h_f(G[r]):h_E(G[r])\to h_{A_{i+1}}(G[r])$ is the same as the kernel of
$\beta(G[r])$. So, the exact sequence
$h_F(G[r])\to h_E(G[r])\to h_{A_{i+1}}(G[r])$ induces an exact sequence
$h_F(G[r])\to h_E(G[r])\to H(G[r])$ and (b) is satisfied.
\end{proof}

\subsection{Representability}
\label{secrep}
\subsubsection{}
We can now state a representability Theorem for 
strongly finitely generated triangulated categories.

\begin{thm}
\label{Brown}
Let $\CT$ be a strongly finitely generated triangulated category and
$H$ be a cohomological functor.

Then, $H$ is locally finitely presented if and only if it 
is a direct summand of a representable functor.
\end{thm}

\begin{proof}
Let $G$ be a $d$-step generator of $\CT$ for some $d\in\BN$.
Let $(A_i)$ be a directed system as in Lemma \ref{approx}.
Then, $\alpha_{d+1}:h_{A_{d+1}}\to H$ is a split surjection
by Proposition \ref{propconv} (iv).
The converse follows from Lemmas \ref{presislocfipres} and \ref{Xthick}.
\end{proof}

Recall that an additive category is {\em Karoubian} if for every 
object $X$ and every idempotent $e\in\End(X)$, there is an object
$Y$ and maps $i:Y\to X$ and $p:X\to Y$ such that $pi=\id_Y$ and $ip=e$.

\begin{cor}
\label{strongKaroubian}
Let $\CT$ be a strongly finitely generated Karoubian triangulated category.
Then, every locally finitely presented
cohomological functor is representable.
\end{cor}

Via Proposition \ref{locfin}, Theorem \ref{Brown} generalizes the following result of
Bondal and Van den Bergh \cite[Theorem 1.3]{BoVdB}.

\begin{cor}
\label{repExtfinite}
Let $\CT$ be an $\Ext$-finite strongly finitely generated
Karoubian triangulated category.
A cohomological functor $H:\CT^\opp\to k\mMOD$ is representable if and
only if it is locally finite.
\end{cor}

The following Lemma is classical:

\begin{lemma}
\label{idemp}
Let $\CT$ be a triangulated category closed under countable multiples. Then,
$\CT$ is Karoubian.
\end{lemma}

\begin{proof}
Given $X\in\CT$ and $e\in\End(X)$ an idempotent, then
$\hocolim(X\xrightarrow{e}X\xrightarrow{e} X\to\cdots)$ is the image of $e$.
\end{proof}

We have a variant of Theorem \ref{Brown}, with a similar proof~:

\begin{thm}
\label{Brown2}
Let $\CT$ be a triangulated category that has a strong complete generator
and $H$ be a cohomological functor that commutes with products.

Then, $H$ is locally finitely presented if and only if it 
is a direct summand of a representable functor.

If $\CT$ is closed under countable multiples, then
$H$ is locally finitely presented if and only if it is
representable.
\end{thm}

\subsubsection{}
Let us now consider cocomplete and compactly generated triangulated
categories --- the ``classical'' setting.

\begin{lemma}
\label{alllfp}
Assume $\CT$ is cocomplete.
Then, every functor is locally finitely presented.
\end{lemma}

\begin{proof}
Let $H$ be a functor and $X\in\CT$.
Let $D=\bigoplus_i X[i]^{|H(X[i])|}$
and $\alpha:h_D\to H$ the canonical map. Then, $\alpha(X[i])$ is
surjective for every $i$. It follows that $H$ is locally finitely generated.

Now, the kernel of  a map $h_E\to H$ will also be locally finitely
generated, hence $H$ is locally finitely presented.
\end{proof}

So, we can derive the classical representability Theorem
(\cite[Theorem 3.1]{Nee3}, \cite[Theorem 5.2]{Ke}, \cite[Lemma 2.2]{Nee2})~:

\begin{thm}
\label{equivgen}
Let $\CT$ be a cocomplete triangulated category generated by a set
$\CS$ of compact objects.
Then, 
\begin{enumerate}
\item
a cohomological functor $\CT^\opp\to k\mMOD$ is representable if and only
if it commutes with products
\item
every object of $\CT$ is a homotopy colimit
of a system $A_1\Rarr{f_1}A_2\Rarr{f_2}\cdots$ almost constant on 
$\langle\overline{\CS}\rangle$ and such that $A_1$ and
the cone of $f_i$ for all $i$ are in $\overline{\CS}$. In particular,
$\CT$ is the smallest full triangulated subcategory containing $\CS$ and
closed under direct sums.
\item
$\CS$ classically generates $\CT^c$.
\end{enumerate}
\end{thm}

\begin{proof}
Let $G=\bigoplus_{S\in\CS} S$.
Let $H:\CT^\opp\to k\mMOD$ be a cohomological functor that commutes with
products.
Let $(A_i,f_i)$ be a directed system constructed as in Lemma \ref{approx} and
$C=\hocolim A_i$. Note that we can assume that
$A_1$ and the cone of $f_i$ are direct sums
of shifts of $G$ (cf Lemmas \ref{approx} and \ref{alllfp}). By
Proposition \ref{propconv} (ii), the system is almost constant on
$\langle\overline{\CS}\rangle$.

The distinguished triangle
$\bigoplus A_i\to \bigoplus A_i \to C\rightsquigarrow$
induces an exact sequence
$H(C)\to\prod H(A_i)\to \prod H(A_i)$, since
$H$ takes direct sums in $\CT$ to products. Consequently, there
is a map $f:h_C\to H$ that makes the following diagram commutative
$$\xymatrix{
& h_C \ar@{.>}[dr]^f \\
\colim h_{A_i}\ar[rr]\ar[ur] && H
}$$
where the canonical maps from $\colim h_{A_i}$ are isomorphisms when
the functors are restricted to $\langle \CS\rangle$ (cf Lemma
\ref{compactcolim}).
So, the restriction of $f$ to $\langle \CS\rangle$ is an
isomorphism. Consequently, $f$ is an isomorphism on the smallest full 
triangulated subcategory $\CT'$ of $\CT$ containing $\CS$ and closed under
direct sums. To conclude, it is enough to show that $\CT'=\CT$ and
we will prove the more precise assertion (2) of the Theorem.

We take $X\in\CT$ and $H=h_X$. Then, $f$ comes from a map
$g:C\to X$. The cone $Y$ of $g$ is zero, since 
$\Hom(S[i],Y)=0$ for all $S\in\CS$ and $i\in\BZ$. Hence,
$g$ is an isomorphism, so (2) holds.

Assume finally that $X\in\CT^c$. Then, $g^{-1}:X\iso C$ factors through
some object of $\langle \CS\rangle_i$ by Proposition \ref{colimsum}, hence
$X\in\langle \CS\rangle_i$.
\end{proof}

\subsubsection{}
We deduce a general duality property for compact objects.
\begin{cor}
\label{Serre}
Let $\CT$ be a cocomplete compactly generated triangulated category over
a field $k$.
Then, there is a faithful functor $S:\CT^c\to \CT$ and bifunctorial
isomorphisms
$$\Hom(C,D)^*\iso \Hom(D,S(C))$$
for $C\in T^c$ and $D\in \CT$.
If $\Hom(C,D)$ is finite dimensional for all $C,D\in\CT^c$, then
$S$ is fully faithful.
\end{cor}

\begin{proof}
Let $C\in \CT^c$. 
The cohomological functor $\Hom(C,-)^*:\CT^\opp\to k\mMOD$ 
commutes with products, hence it is representable by
an object $S(C)\in \CT$ by Theorem \ref{equivgen}.
By Yoneda's Lemma,
this defines a functor $S:\CT^c\to \CT$.
Now, if $D\in \CT^c$, then $S$ is equal to the composition
$$\Hom(C,D)\xrightarrow{\can} \Hom(C,D)^{**}
\iso \Hom(D,S(C))^*\iso \Hom(S(C),S(D)).$$
\end{proof}

Whenever $\CT^c$ admits a Serre functor, it must be the restriction of
the $S$ above.
\begin{cor}
Let $\CT$ be a cocomplete compactly generated triangulated category over
a field $k$. Assume there is a self-equivalence $S'$ of $\CT^c$ together
with bifunctorial isomorphisms 
$$\Hom(C,D)^*\iso \Hom(D,S'(C))$$
for $C,D\in T^c$.
Then, $S$ takes values in $\CT^c$ and there is a unique
isomorphism $S'\iso S$ making the following diagram commutative
for any $C,D\in\CT^c$
$$\xymatrix{
 &\Hom(D,S'(C))\ar[dd]\\
\Hom(C,D)^*\ar[ur]\ar[dr] \\
 &\Hom(D,S(C))
}$$
\end{cor}

\begin{proof}
A bifunctorial isomorphism $\Hom(D,S'(C))\iso\Hom(D,S(C))$ comes
from a unique functorial 
map $S'(C)\to S(C)$. Its cone is right orthogonal to $\CT^c$,
hence it is zero, since $\CT$ is generated by $\CT^c$.
\end{proof}

\subsection{Finiteness for objects}
\label{subsecfinob}
\subsubsection{}
We say that $C$ is {\em cohomologically
locally bounded} (resp. {\em bounded above}, resp. {\em bounded below}, resp.
{\em finitely generated}, resp. {\em finitely presented}, resp. {\em finite})
if the restriction of $h_C$ to $\CT^c$ has that property.

From Lemma \ref{presislocfipres}, we deduce

\begin{lemma}
Let $C\in\CT^c$. Then, $C$ is cohomologically locally finitely presented.
\end{lemma}

\begin{lemma}
\label{carclfp}
Let $C\in\CT$ be cohomologically locally finitely generated. Then, $C$ is
cohomologically locally finitely presented if and only if given $X\in\CT^c$,
$E\in\CT^c$ and $\beta:E\to C$ such that $\Hom(X[i],\beta)$ is surjective for every
$i\in\BZ$, then the cocone of $\beta$ is 
cohomologically locally finitely generated.
\end{lemma}

\begin{proof}
Let $F$ be the cocone of $\beta$. We have an exact sequence
$$0\to \Hom(X[i],F)\to \Hom(X[i],E)\xrightarrow{\Hom(X[i],\beta)}\Hom(X[i],C)\to 0.$$
The Lemma follows now from Lemma \ref{surj}.
\end{proof}

From Lemma \ref{approx} and Proposition \ref{propconv} (ii), we obtain

\begin{lemma}
\label{clfpfactor}
Assume $\CT$ is cocomplete and generated by a compact object $G$. Let
$C\in\CT$. Let $C$ be a cohomologically locally finitely presented object
of $\CT$. Then,
there is a system $A_1\to A_2\to\cdots$ in $\CT^c$ which is almost
constant for $\langle\overline{\{G[i]\}_{i\in\BZ}}\rangle$ and an
isomorphism $\hocolim A_i\iso C$.

In particular, given $d\ge 0$, there is $D\in\CT^c$ and $f:D\to C$ such that
every map from an object of $\langle\overline{\{G[i]\}_{i\in\BZ}}\rangle_d$
to $C$ factors through $f$.
\end{lemma}

From Propositions \ref{lfpstable} and \ref{lb}, we deduce

\begin{prop}
\label{clfpthick}
The full subcategory of $\CT$ of cohomologically 
locally finitely presented (resp. bounded) objects is a thick subcategory.
\end{prop}

Note that the full subcategory of cohomologically locally 
bounded (resp. bounded above, resp. bounded below) objects is also a thick
subcategory.

From Theorem \ref{Brown}, we deduce

\begin{cor}
\label{lfpiscompact}
Let $\CT$ be a triangulated category such that
$\CT^c$ is strongly finitely generated. Then, $C\in\CT^c$ if and only if
$C$ is cohomologically locally finitely presented.
\end{cor}

\begin{rem}
\label{CKN}
Not all cohomological functors
on $\CT^c$ are isomorphic to the restriction of $h_C$, for some $C\in\CT$.
This question has been studied for example in \cite{Nee4,Bel,ChKeNee}.
Let us mention the following result \cite[Lemma 2.13]{ChKeNee}~:
let $\CT$ be a cocomplete and compactly generated triangulated category.
Assume $k$ is a field. Let $H$ be a cohomological functor on $\CT^c$
with value in the category $k\mMod$ of finite dimensional vector spaces. Then
there is $C\in\CT$ such that $H$ is isomorphic to the restriction of 
$h_C$ to $\CT^c$.
\end{rem}

\section{Localization}
\label{seclocalization}
\subsection{Compact objects}
\subsubsection{}
Let us recall Thomason's classification of dense subcategories
\cite[Theorem 2.1]{Th}~:

\begin{thm}
\label{dense}
Let $\CT$ be a triangulated category and $\CI$ a dense full
triangulated subcategory. Then, an object of $\CT$ is isomorphic to
an object of $\CI$ if and only if its class is in the image of
the canonical map $K_0(\CI)\to K_0(\CT)$.
\end{thm}

The following Lemma is proved in \cite[Lemma 1.5]{BoNee}.
\begin{lemma}
\label{sumquotient}
Let $\CT$ be a cocomplete triangulated category and $\CI$ be a thick
subcategory closed under direct sums. Then, $\CT/\CI$ is cocomplete
and the quotient functor $\CT\to\CT/\CI$ commutes with direct sums.
\end{lemma}

The following is a version of Thomason-Trobaugh-Neeman's
Theorem \cite[Theorem 2.1]{Nee2}.

\begin{thm}
\label{locNeeman}
Let $\CT$ be a cocomplete and compactly generated triangulated category.
Let $\CI$ a full triangulated subcategory closed under direct sums and
compactly generated in $\CT$.
Denote by $F:\CT\to\CT/\CI$ the quotient functor. Then,

\begin{itemize}
\item[(i)]
$\CI$ is a cocomplete compactly generated triangulated category and
$\CI^c=\CI\cap\CT^c$.
\item[(ii)]
Given $X\in\CT^c$ and $Y\in\CT$,
the canonical map
$$\lim \Hom_\CT(X',Y)\iso \Hom_{\CT/\CI}(FX,FY)$$
is an isomorphism,
where the limit is taken over the maps $X'\to X$ whose cone is in
$\CI^c$.
Also, if $FY$ is in $F(\CT^c)$, then, there
is $C\in\CT^c$ and $f:C\to Y$ such that $F(f)$ is an isomorphism.
\item[(iii)]
$F$ commutes with direct sums and
the canonical functor 
$\CT^c/\CI^c\to \CT/\CI$
factors through a fully faithful
functor $G:\CT^c/\CI^c\to (\CT/\CI)^c$.

\item[(iv)]
An object of $(\CT/\CI)^c$ is isomorphic to an object in the
image of $G$ if and only if its class is in the image of $K_0(G)$.
\end{itemize}
\end{thm}

\begin{proof}
It is clear that $\CI$ is cocomplete and that $\CT^c\cap\CI\subset\CI^c$.
Let $\CS_\CI$ be a set of
objects of $\CT^c\cap\CI$ that generates $\CI$. It follows from
Theorem \ref{equivgen} (3) that $\CS_\CI$ classically generates $\CI^c$.
Since $\CT^c\cap\CI$ is a thick subcategory of $\CI$, it follows that
$\CI^c=\CT\cap\CI^c$ and (i) is proven.

Let $X\in\CT^c$ and $Y\in\CT$. Let $\phi:W\to X$ and $\psi:W\to Y$ with
$W\in\CT$.
Let $Z$ be a cone of $\phi$ and assume $Z\in\CI$.
By Theorem \ref{equivgen} (2) and Proposition \ref{colimsum},
$X\to Z$ factors through a map $\alpha:X\to Z'$ for some $Z'\in \CI\cap\CT^c$.
Let $X'$ be the cocone of $\alpha$.
The map $X'\to X$ factors as a composition $\phi\zeta$. This shows
(ii).
$$
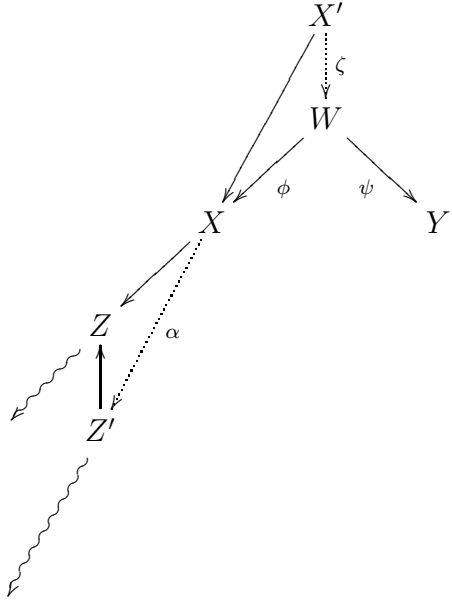
\xymatrix{
& &   & X'\ar@{.>}[d]^\zeta\ar[ddl] \\
& &   & W\ar[dr]_\psi\ar[dl]^\phi \\
& & X\ar[dl]\ar@{.>}[ddl]^\alpha & & Y \\
& Z \ar@{~>}[dl] \\
& Z'\ar[u]\ar@{~>}[ddl]\\
&\\
&\\
}$$
\smallskip
Since $\CT$ is cocomplete and the direct sum in $\CT$ of
objects of $\CI$ is in $\CI$, it follows from
Lemma \ref{sumquotient} that $F$ commutes with direct sums.

Let now $X\in\CT^c$ and
$\{Z_i\}$ be a family of elements of $\CT$. Let $f:F(X)\to 
\bigoplus_i F(Z_i)=F(\bigoplus_i Z_i)$. There is $\phi:X'\to X$ and
$\psi:X'\to \bigoplus_i Z_i$ with the cone of $\phi$ in $\CI\cap\CT^c$ and
$f=F(\psi)F(\phi)^{-1}$. Since $X'$ is compact, $\psi$ factors through
a finite sum of $Z_i$'s, hence $f$ factors through
a finite sum of $F(Z_i)$'s. Consequently, $F(X)$ is compact. The
fully faithfulness of $G$ comes from (ii).

\smallskip
Let us now prove (iv).
By Theorem \ref{equivgen} (3), 
$(\CT/\CI)^c$ is classically generated by $F(\CT^c)$.
Since $F(\CT^c)$ is a full triangulated subcategory of
$(\CT/\CI)^c$, it is dense.
The result follows now from Theorem \ref{dense}.
\end{proof}

\begin{cor}
Let $\CT$ be a cocomplete and compactly generated triangulated category.
Let $\CI$ be a full triangulated subcategory closed under direct sums and
generated by an object $G\in \CT^c\cap\CI$ such that for all $C\in\CT^c$, then 
$\Hom(C,G[i])=0$ for $|i|\gg 0$.

If $\CT^c$ is strongly finitely generated, then $(\CT/\CI)^c$ is strongly
finitely generated.
\end{cor}

\begin{rem}
Let $\CT$ be a cocomplete triangulated category generated by a set
$\CE$ of compact objects and let
$\CI$ a thick subcategory closed under direct sums. If the inclusion
functor $\CI\to\CT$ has a left adjoint $G$, then 
$G(\CE)$ is a generating set for $\CI$ and it consists of compact objects
of $\CT$.
\end{rem}

\subsection{Proper intersections of Bousfield subcategories}
\label{seccov}
\subsubsection{}
Let $\CT$ be a triangulated category and $\CI$ be a thick subcategory.
We have a canonical fully faithful functor $i_*:\CI\to \CT$ and
a canonical essentially surjective quotient functor $j^*:\CT\to\CT/\CI$.
We say that there is an {\em exact sequence} of triangulated categories
$$0\to \CI\xrightarrow{i_*} \CT\xrightarrow{j^*} \CT/\CI\to 0$$
We say that $C\in\CT$ is $\CI$-{\em local} if
$\Hom(M,C)=0$ for all $M\in\CI$. Note that given $C,D\in\CT$ with
$D$ an $\CI$-local object, then
$\Hom(C,D)\iso \Hom(j^*C,j^*D)$.

\smallskip
Let $\CI'$ be a thick subcategory of $\CI$. Then, we have a commutative
diagram of exact sequences of triangulated categories
$$\xymatrix{
&&0\ar[d] & 0\ar[ld] \\
&& \CI'\ar[dl]\ar[d] && 0\\
0\ar[r] & \CI\ar[r]\ar[d] & \CT\ar[r]\ar[d] & \CT/\CI\ar[r]\ar[ur] & 0 \\
0\ar[r] & \CI/\CI'\ar[r]\ar[d] & \CT/\CI'\ar[ur]\ar[d] \\
& 0 & 0
}$$

\subsubsection{}
Let us recall the construction of Bousfield localization (cf \eg\
\cite[\S 9.1]{Nee5}).

We say that $\CI$ is a {\em Bousfield} subcategory if the
quotient functor $j^*:\CT\to\CT/\CI$ has a right adjoint $j_*$.
We then denote by $\eta:\id_\CT\to j_*j^*$ the corresponding unit.

Assume $\CI$ is a Bousfield subcategory.
Note that $C$ is $\CI$-local if and only
if $\eta(C):C\to j_*j^*C$ is an isomorphism if and only
if $C\simeq j_* C'$ for some $C'\in\CT/\CI$.

We denote by $i_*:\CI\to\CT$ the inclusion functor.
Let $C\in\CT$ and $C'$ be the cocone of $\eta(C)$. We have
$j^*C'=0$, hence $C'\in\CI$. Since $j_*j^*C[-1]$ is $\CI$-local, the object
$C'$ is well defined up to unique isomorphism. So, there is a functor
$i^!:\CT\to\CI$ and a map $\eps:i_*i^!\to \id_\CT$ such that the
following triangle is distinguished
\begin{equation}
\label{openclosed}
i_*i^!\xrightarrow{\eps} \id_\CT \xrightarrow{\eta} j_*j^*\rightsquigarrow.
\end{equation}
Furthermore, $\eps$ provides $(i_*,i^!)$ with the structure of an adjoint
pair.

Since $i_*$ and $j^*$ have right adjoints, they commute with direct sums.
Also, $\CI$ is closed under direct sums (taken in $\CT$) and we
have $i^!j_*=j^*i_*=0$. The unit of adjunction $\id_\CI\iso i^!i_*$ is an 
isomorphism, as well as the counit $j^*j_*\iso \id_{\CT/\CI}$.

\smallskip
Let $\CI$ be a thick subcategory of $\CT$. Then, the following
conditions are equivalent
\begin{itemize}
\item $\CI$ is a Bousfield subcategory
\item
for any $C\in\CT$,
there is a distinguished triangle 
$C_1\to C\to C_2\rightsquigarrow$ with $C_1\in\CI$ and $C_2$ an
$\CI$-local object.
\item the restriction of $j^*$ to the full subcategory of $\CI$-local
objects is an equivalence.
\end{itemize}

\smallskip
Let $\CI'$ be a Bousfield subcategory of $\CT$ containing $\CI$.
Then, $\CI$ is a Bousfield subcategory of $\CI'$. The right adjoint to
the inclusion of $\CI$ in $\CI'$ is $i^!i'_*$. Also, 
$\CI'/\CI$ is a Bousfield subcategory of $\CT/\CI$ and the left adjoint
to the quotient $\CT/\CI\to\CT/\CI'$ is $j^*j'_*$.

\smallskip
Assume $\CT$ is cocomplete and compactly generated and $\CI$ is
a full triangulated subcategory closed under direct sums. Then,
$\CI$ is a Bousfield subcategory \cite[Example 8.4.5]{Nee5}.
Indeed, given $D\in\CT/\CI$, the functor
$\Hom(j^*(-),D):\CT^\opp\to k\mMOD$ is cohomological and commutes
with products (Theorem \ref{locNeeman}),
hence is representable by Theorem \ref{equivgen}.
The thickness follows from Lemma \ref{idemp}.

\begin{rem}
Let $\CT$ be a cocomplete compactly generated triangulated category and $\CI$ a
Bousfield subcategory. If $C\in\CT$ is cohomologically locally bounded,
then $i^!C$ is cohomologically bounded. An object $C'\in\CT/\CI$ is
cohomologically locally bounded if and only if $j_*C'$ is 
cohomologically locally bounded.
\end{rem}

\subsubsection{}
\label{secproper}

Let $\CI_1$ and $\CI_2$ be two Bousfield subcategories of $\CT$.

\begin{lemma}
\label{properly}
The following assertions are equivalent
\begin{enumerate}
\item $i_{1*}i_1^!(\CI_2)\subset\CI_2$ and $i_{2*}i_2^!(\CI_1)\subset\CI_1$
\item $j_{1*}j_1^*(\CI_2)\subset\CI_2$ and $j_{2*}j_2^*(\CI_1)\subset\CI_1$
\item the canonical functor
$\CI_1/(\CI_1\cap\CI_2)\oplus \CI_2/(\CI_1\cap\CI_2)\to\CT/(\CI_1\cap\CI_2)$ is
fully faithful
\item
given $M_1\in\CI_1$ and $M_2\in\CI_2$, every map $M_1\to M_2$ and every map
$M_2\to M_1$ factors through an object of $\CI_1\cap\CI_2$.
\end{enumerate}
\end{lemma}

\begin{proof}
Given $N\in\CI_2$, we have a distinguished triangle
$i_{1*}i_1^!N\to N\to j_{1*}j_1^*N\rightsquigarrow$.
This shows immediately the equivalence between (1) and (2).

Let $f:M\to N$ with $M\in \CI_1$. Then, there is $g:M\to i_{1*}i_1^!N$ such that
$f=\eta_1(N)g$.
It is now clear that (1)$\Rightarrow$(4).  Assume (4). Then, there is
$L\in\CI_1\cap\CI_2$ and $\phi:i_{1*}i_1^!N\to L$ and $\psi:L\to N$
such that $\eps(N)=\psi\phi$. Now, there is $\phi':L\to i_{1*}i_1^!N$
such that $\psi=\eps(N)\phi'$. So, $\eps(N)(1-\phi'\phi)=0$. Since
the canonical map $\End(i_{1*}i_1^!N)\to \Hom(i_{1*}i_1^!N,N),\ h\mapsto
\eps(N)h$ is injective, it follows that $i_{1*}i_1^!N$ is a direct of
$L$, hence $i_{1*}i_1^!N\in\CI_2$. So, (4)$\Rightarrow$(1).

A map in $\CT$ factors through an object of $\CI_1\cap\CI_2$ if and only
if it becomes $0$ in $\CT/(\CI_1\cap\CI_2)$. This shows the equivalence of
(3) and (4).
\end{proof}

We say that $\CI_1$ and $\CI_2$ {\em intersect properly} if
the assertions of Lemma \ref{properly} are satisfied.
This property passes to intersections, unions, quotients... A collection
of Bousfield subcategories any two of which intersect properly behaves like
a collection of closed subsets.

We will identify $\CI_1/(\CI_1\cap\CI_2)$ with its essential image in
$\CT/\CI_2$.

There are commutative diagrams of inclusions of subcategories and
of quotients of categories
$$\xymatrix{
& \CI_1\ar[dr]^-{i_{1*}} &&& 
 & \CT/\CI_1\ar[dr]^-{j_{1\cup}^*} \\
\CI_1\cap\CI_2\ar[rr]^-{i_{\cap*}}\ar[ur]^-{i_{1\cap*}}\ar[dr]_-{i_{2\cap*}}&&
 \CT && 
 \CT\ar[ur]^-{j_1^*}\ar[dr]_-{j_2^*}\ar[rr]^-{j_\cup^*} &&
 \CT/\langle\CI_1\cup\CI_2\rangle_\infty \\
& \CI_2\ar[ur]_-{i_{2*}} &&&
 & \CT/\CI_2\ar[ur]_-{j_{2\cup}^*}
}$$

\begin{lemma}
\label{inter1}
Assume $\CI_1$ and $\CI_2$ intersect properly. Let $\{a,b\}=\{1,2\}$. Then,
\begin{itemize}
\item $\CI_1\cap\CI_2$ and $\langle\CI_1\cup\CI_2\rangle_\infty$
are Bousfield subcategories of $\CT$.
\item We have
$i_{\cap *}i_\cap^!\simeq i_{a*}i_a^!i_{b*}i_b^!$ and
$j_{\cup *}j_\cup^*\simeq j_{a*}j_a^*j_{b*}j_b^*$.
\item
There are commutative diagrams
$$\xymatrix{
\CI_a\ar[r]^-{i_{a*}}\ar[d]_{i_{a\cap}^!} & \CT\ar[d]^{i_b^!} & &
  \CT/\CI_a \ar[r]^-{j_{a*}}\ar[d]_-{j_{a\cup}^*}&\CT\ar[d]^-{j_b^*} \\
\CI_1\cap\CI_2\ar[r]_-{i_{b\cap*}} & \CI_b &&
  \CT/\langle\CI_1\cup\CI_2\rangle_\infty\ar[r]_-{j_{b\cup*}}&\CT/\CI_b
}$$
\item
The canonical functor 
$\CI_a/(\CI_1\cap\CI_2)\iso \langle\CI_1\cup\CI_2\rangle_\infty/\CI_b$ is
an equivalence
and we have a commutative diagram of exact sequences of triangulated categories
$$\xymatrix{
& & 0\ar[dr] && 0 \\
& 0 \ar[dr] & & \CI_1/(\CI_1\cap\CI_2)\ar[ur]\ar[dr] && 0\\
0\ar[dr] && \CI_1\ar[ur]\ar[dr] & & \CT/\CI_2\ar[ur]\ar[dr] && 0\\
&\CI_1\cap\CI_2\ar[ur]\ar[dr] & & \CT\ar[ur]\ar[dr] & &
 \CT/\langle\CI_1\cup\CI_2\rangle_\infty\ar[dr]\ar[ur] \\
0\ar[ur] & & \CI_2\ar[ur]\ar[dr] & & \CT/\CI_1\ar[ur]\ar[dr] && 0\\
& 0 \ar[ur] && \CI_2/(\CI_1\cap\CI_2)\ar[ur]\ar[dr] && 0 \\
&&0\ar[ur] &&0
}$$
\end{itemize}
\end{lemma}

\begin{proof}
Let $C\in\CT$. We have distinguished triangles
$$i_{1*}i_1^!C\to C\to j_{1*}j_1^*C\rightsquigarrow \text{ and }
i_{2*}i_2^!i_{1*}i_1^!C\to i_{1*}i_1^!C\to j_{2*}j_2^*i_{1*}i_1^!C\rightsquigarrow$$
The octahedral axiom shows that there are $C'\in\CT$ and distinguished triangles
$$i_{2*}i_2^!i_{1*}i_1^!C\to C\to C'\rightsquigarrow \text{ and }
j_{2*}j_2^!i_{1*}i_1^!C \to C'\to j_{1*}j_1^*C\rightsquigarrow$$
Since $C'$ is $(\CI_1\cap\CI_2)$-local and
$i_{2*}i_2^!i_{1*}i_1^!C\in\CI_1\cap\CI_2$,
we deduce that $\CI_1\cap\CI_2$ is a Bousfield subcategory of $\CT$.
The map $i_{2*}i_2^!i_{1*}i_1^!C\to C$ factors uniquely through
the canonical map $i_{\cap *}i_\cap^!C\to C$ and similarly
the map $i_{\cap *}i_\cap^!C\to C$ factors uniquely through
$i_{2*}i_2^!i_{1*}i_1^!C\to C$ and this provides functorial inverse morphisms
between  $i_{\cap *}i_\cap^!C$ and $i_{2*}i_2^!i_{1*}i_1^!C$.

The case of $\langle\CI_1\cup\CI_2\rangle_\infty$ is similar, using
$i_{2*}i_2^!j_{1*}j_1^!C\to j_{1*}j_1^*C\to j_{2*}j_2^*j_{1*}j_1^*C\rightsquigarrow$
as a second distinguished triangle.

\smallskip
We have $i_1^!i_{2*}(\CI_2)\subset \CI_1\cap\CI_2$, hence the canonical
map $i_{1\cap*}i_{1\cap}^!i_1^!i_{2*}\iso i_1^!i_{2*}$ is an isomorphism.
Now, we have canonical isomorphisms
$$i_{1\cap*}i_{2\cap}^!\iso
i_{1\cap*}i_{2\cap}^!i_2^!i_{2*}\iso
i_{1\cap*}i_{1\cap}^!i_1^!i_{2*}$$
and we get the first commutative square. The proof of the commutativity
of the second square is similar.

The last assertion is clear.
\end{proof}

\begin{lemma}
Let $\CF$ be a finite family of Bousfield subcategories of $\CT$ any two of
which intersect properly.

Given $\CF'$ a subset of $\CF$, then $\bigcap_{\CI\in\CF'}\CI$ 
(resp. $\langle \bigcup_{\CI\in\CF'}\CI\rangle_\infty$) is a Bousfield
subcategory of $\CT$ that intersects properly any subcategory in $\CF$.

Given $\CI,\CI_1,\CI_2\in\CF$, then $\CI_1/(\CI\cap\CI_1)$ and
$\CI_2/(\CI\cap\CI_2)$ are Bousfield subcategories
of $\CT/\CI$ that intersect properly.
\end{lemma}

\begin{proof}
By induction, it is enough to prove the first assertion when $\CF'$ has
two elements, $\CF'=\{\CI_2,\CI_3\}$ and the result is then given by Lemma
\ref{inter1}.

Let $M\in\CI_1$, $N\in\CI_2$, $L\in\CT$ and $f:L\to M$ and $g:L\to N$
such that $f$ becomes an isomorphism in $\CT/\CI$. Then, the cone of $f$ is
in $\CI$, so $L\in \langle \CI_1\cup\CI\rangle_\infty$.
The first part of the Lemma shows that $\langle \CI_1\cup\CI\rangle_\infty$
and $\CI_2$ intersect properly. It follows that $g$ factors through
an object of $\langle \CI_1\cup\CI\rangle_\infty\cap\CI_2$.
Consequently, the image of $g$ in $\CT/\CI$ factors through an object
of $\left(\CI_1/(\CI\cap\CI_1)\right)\cap \left(\CI_2/(\CI\cap\CI_2)\right)$.
We have shown that every
map in $\CT/\CI$ between $M$ and $N$ factors through an object of
$\left(\CI_1/(\CI\cap\CI_1)\right)\cap \left(\CI_2/(\CI\cap\CI_2)\right)$
and we deduce the proper intersection property.
\end{proof}

\subsubsection{}
We have two Mayer-Vietoris triangles ("open" and "closed" cases).
\begin{prop}
\label{mayervietoris}
Assume $\CI_1$ and $\CI_2$ intersect properly.

\begin{itemize}
\item[(1)] If $\CT=\langle\CI_1,\CI_2\rangle_\infty$, then,
there are isomorphisms of functors
$i_{\cap}^!\iso i_{1\cap}^!i_1^!$ and
$i_{\cap}^!\iso i_{2\cap}^!i_2^!$ giving
a distinguished triangle of functors
$$i_{\cap*}i_\cap^!\xrightarrow{i_{1*}\eps_{1\cap}i_1^!+
i_{2*}\eps_{2\cap}i_2^!} i_{1*}i_1^!\oplus i_{2*}i_2^!
\xrightarrow{\eps_1-\eps_2}\id_\CT \rightsquigarrow.$$
\item[(2)]
If $\CI_1\cap \CI_2=0$, then there are isomorphisms of functors
$j_{1*}j_{1\cup*}\iso j_{\cup*}$ and
$j_{2*}j_{2\cup*}\iso j_{\cup*}$
giving a distinguished triangle of functors
$$\id_\CT\xrightarrow{\eta_1-\eta_2}j_{1*}j_1^*+j_{2*}j_2^*
\xrightarrow{j_{1*}\eta_{1\cup}j_1^*+j_{2*}\eta_{2\cup}j_2^*}
 j_{\cup*}j_\cup^*\rightsquigarrow.$$
\end{itemize}
\end{prop}

\begin{proof}
It is an easy general categorical fact that there is
an isomorphism of functors
$i_{\cap}^!\iso i_{a\cap}^!i_a^!$ such that
$\eps_\cap=\eps_a\circ (i_{a*}\eps_{a\cap}i_a^!)$. Then,
$(\eps_1-\eps_2)\circ (i_{1*}\eps_{2\cap}i_1^!+
i_{2*}\eps_{1\cap}i_2^!)=0$.

Given $M\in\CI_2$ an $(\CI_1\cap\CI_2)$-local object, then
$i_{2*}M$ is $\CI_1$-local. Since the canonical functor
$\CI_2/(\CI_1\cap\CI_2)\iso \CT/\CI_1$ is an equivalence, it follows
that the $\CI_1$-local objects of $\CT$ are contained in $\CI_2$, hence
$j_{1*}(\CT/\CI_1)\subset\CI_2$. As a consequence, given $N\in\CT/\CI_1$
such that $i_2^!j_{1*}N=0$, we have $N=0$.
Consider now $C\in\CT$ such that $i_1^!C=i_2^!C=0$.
Then, $C\iso j_{1*}j_1^*C$. Since $i_2^!C=0$, it follows that
$j_1^*C=0$, hence $C=0$.
We deduce that in order to prove that the triangle of the Lemma is
distinguished, it is sufficient to prove so after applying the functor
$i_1^!$ and after applying the functor $i_2^!$.

The map $i_1^!i_{2*}\eps_{2\cap}i_2^!:
i_1^!i_{2*}i_{2\cap*}i_{2\cap}^!i_2^!\to
i_1^!i_{2*}i_2^!$ is an isomorphism since
$i_1^!i_{2*}\simeq i_{1\cap*}i_{2\cap}^!$ (Lemma \ref{inter1}).
As the map $i_1^!\eps_1$ is an isomorphism, we deduce that after applying
$i_1^!$, the triangle is a split distinguished triangle.

The second assertion has a similar proof.
\end{proof}

We say that two subcategories $\CC_1$ and $\CC_2$ of a category $\CC$
are {\em orthogonal} if $\Hom(C_1,C_2)=\Hom(C_2,C_1)=0$ for all
$C_1\in\CC_1$ and $C_2\in\CC_2$. Note that this is equivalent to requiring
that $\CI_1\cap\CI_2=0$ and $\CI_1$ and $\CI_2$ intersect properly.

\subsection{Coverings}
\label{seccoverings}
\subsubsection{}

The following proposition shows that compactness is a local property,
in a suitable sense~:

\begin{prop}
\label{localcompact}
Let $\CI_1$ and $\CI_2$ be two orthogonal Bousfield subcategories of $\CT$.
Let $C\in\CT$. If $j_1^*C$, $j_2^*C$ and $j_\cup^*C$ are
compact, then $C$ is compact.
\end{prop}

\begin{proof}
Let $\CF$ be a set of objects of $\CT$ whose direct sum exists. Let
$a\in\{1,2,\cup\}$.
We have canonical isomorphisms
\begin{align*}
\bigoplus_D \Hom(C,j_{a*}j_a^*D)&\iso
\bigoplus_D \Hom(j_a^*C,j_a^*D)\iso
\Hom(j_a^*C,\bigoplus_D j_a^*D)\iso
\Hom(j_a^*C,j_a^*\bigoplus_D D)\\
&\iso\Hom(C,j_{a*}j_a^*\bigoplus_D D).
\end{align*}

We have a commutative diagram
$$\tiny\xymatrix{
\cdots \ar[r] & \bigoplus_D \Hom(C,D)\ar[d]\ar[r] &
 \bigoplus_D \Hom(C,j_{1*}j_1^*D)\oplus
\bigoplus_D \Hom(C,j_{2*}j_2^*D)\ar[r]\ar[d]_\sim & 
\bigoplus_D \Hom(C,j_{\cup*}j_\cup^*D)\ar[d]_\sim\ar[r] & \cdots\\
\cdots\ar[r] &\Hom(C,\bigoplus_D D)\ar[r] &
 \Hom(C,j_{1*}j_1^*\bigoplus_D D)\oplus
\Hom(C,j_{2*}j_2^*\bigoplus_D D)\ar[r] & 
\Hom(C,j_{\cup*}j_\cup^*\bigoplus_D D)\ar[r] & \cdots
}$$
where the exact horizontal rows come from the Mayer-Vietoris triangles
(Proposition \ref{mayervietoris} (2)).
It follows that the canonical map
$\bigoplus_D\Hom(C,D)\iso \Hom(C,\bigoplus_D D)$
is an isomorphism.
\end{proof}

Combining Theorem \ref{locNeeman} and Proposition \ref{localcompact}, we get
\begin{cor}
\label{caractcompact}
Let $\CT$ be a compactly generated cocomplete triangulated category
and let $\CI_1$ and $\CI_2$ be two orthogonal Bousfield subcategories of $\CT$.
Assume $\CI_a$ is compactly generated in $\CT$ for $a\in\{1,2\}$.

Let $C\in\CT$. Then, $C$ is compact if and only if $j_1^*C$ and
$j_2^*C$ are compact.
\end{cor}

\begin{proof}
The only new part is that the compactness of $j_\cup^* C$
follows from that of $j_1^*C$.
Since compact objects of $\CT$ remain compact in $\CT/\CI_1$
(Theorem \ref{locNeeman}), it follows
that $\CI_2$ is compactly generated in
$\CT/\CI_1$. So, if $j_1^*C$ is
compact, then $j_\cup^*C$ is compact (Theorem \ref{locNeeman} again).
\end{proof}

\subsubsection{}
We have now a converse to the localization Theorem \ref{locNeeman}:

\begin{prop}
\label{converseloc}
Let $\CT$ be a triangulated category and $\CI$
be a Bousfield subcategory of $\CT$.

Let $\CE$ be a set of objects of $\CI\cap\CT^c$ generating $\CI$ and
$\CE'$ be a set of objects of $\CT^c$ which generates $\CT/\CI$.
Then $\CT$ is generated by the set $\CE\cup\CE'$.
\end{prop}

\begin{proof}
Let $C\in\CT$ such that $\Hom(D[n],C)=0$ for all $D\in\CE$ and $n\in\BZ$. Then,
using the distinguished triangle (\ref{openclosed}), we get
$\Hom(D[n],i_*i^!C)=0$, hence $i^!C=0$. If follows that $C$ is
$\CI$-local.

Assume now in addition $\Hom(D'[n],C)=0$ for all $D'\in\CE'$ and $n\in\BZ$.
We have $C\iso j_*j^*C$, hence
$\Hom(j^*D'[n],j^*C)\iso \Hom(D'[n],j_*j^*C)=0$. So, $j^*C=0$ and finally $C=0$.
\end{proof}

\begin{prop}
\label{constructiongenerator}
Let $\CT$ be a cocomplete
triangulated category and $\CI_1$, $\CI_2$ be two orthogonal
Bousfield subcategories.
Assume 
\begin{itemize}
\item
$\CT/\CI_a$ is compactly generated and
\item
$\CI_b$ is compactly generated in $\CT/\CI_a$
\end{itemize}
for $\{a,b\}=\{1,2\}$.

Then, $\CT$ is compactly generated.

More precisely, let $\CE$ be a generating set of objects of
$\CI_2$ which are compact in $\CT/\CI_1$
and let $\CE'$ be a set of objects of $(\CT/\CI_2)^c$
generating $\CT/\CI_2$. Then, 
\begin{itemize}
\item $\CE\subset \CT^c$
\item given $M\in\CE'$, there is $\tM\in\CT^c$
such that $j_2^*\tM\simeq M\oplus M[1]$
\item $\CE\cup\{\tM\}_{M\in\CE'}$ generates $\CT$.
\end{itemize}
Let $\CJ$ be a Bousfield subcategory of $\CT$ intersecting properly
$\CI_1$ and $\CI_2$. Assume
\begin{itemize}
\item
$\CJ/(\CI_a\cap\CJ)$ is compactly generated in $\CT/\CI_a$ and
\item
$\CI_b\cap\CJ$ is compactly generated in $\CT/\CI_a$
\end{itemize}
for $\{a,b\}=\{1,2\}$.

Then, $\CJ$ is compactly generated in $\CT$.
\end{prop}

\begin{proof}
Since $\CT$ is cocomplete and $\CI_a$ is Bousfield,
it follows that $\CT/\CI_a$ is cocomplete.

Let $\CE$ be a generating set of objects of $\CI_2$
which are compact in $\CT/\CI_1$.
Given $C\in\CE$, we have $j_2^*C=j_\cup^*C=0$ and
$j_1^*C$ is compact. It follows from Proposition \ref{localcompact} that
$C$ is a compact object of $\CT$. In particular $\CI_2$ is compactly
generated.

Let $\CE'$ be a set of compact objects generating $\CT/\CI_2$.
Let $M\in\CE'$ and $D_2=M\oplus M[1]$. By Theorem \ref{locNeeman},
$D_\cup=j_{2\cup}^*D_2$ is compact and
there is $D_1\in(\CT/\CI_1)^c$ with an isomorphism $j_{1\cup}^*D_1\iso D_\cup$.
Let now $\tM$ be the cocone of the sum of canonical maps
$j_{2*}D_2\oplus j_{1*}D_1\to j_{\cup*}D_\cup$. We have
$j_a^*\tM\simeq D_a$ for $a\in\{1,2,\cup\}$. It follows from Proposition
\ref{localcompact} that $\tM$ is compact. Let $\CE'_2=\{\tM\}_{\tM\in\CE'}$.
Now, Proposition \ref{converseloc} shows that $\CE\cup \CE'_2$ generates $\CT$.

For the case of $\CJ$, we apply the first part of the Proposition to
the cocomplete triangulated category $\CJ$ with its orthogonal
Bousfield categories $\CI_1\cap\CJ$ and $\CI_2\cap\CJ$. We obtain a generating
set $\CE_\CJ$ of objects of $\CJ$ with the property that their images in
$\CJ/(\CI_1\cap\CJ)$ and $\CJ/(\CI_2\cap\CJ)$
are compact. These objects are thus compact in
$\CT/\CI_1$ and $\CT/\CI_2$ by Theorem \ref{locNeeman}. Since
$\CI_1$ and $\CI_2$ are compactly generated in $\CT$, it follows from
Corollary \ref{caractcompact} that $\CE_\CJ\subset\CT^c$.
\end{proof}

\subsubsection{}

A {\em cocovering} of $\CT$ is
a finite set $\CF$ of Bousfield subcategories of $\CT$
any two of which intersect properly and such that $\cap_{\CI\in\CF}\CI=0$.

The following result gives a construction of a compact generating set from
(relative) compact generating sets for the quotients $\CT/\CI$.

\begin{thm}
\label{cocovering}
Let $\CF$ be a cocovering of $\CT$.
\begin{itemize}
\item Let $C$ be an object of $\CT$ which is compact in
$\CT/\langle\bigcup_{\CI\in\CF'}\CI\rangle$ for all non empty
$\CF'\subset\CF$. Then, $C$ is compact in $\CT$.
\end{itemize}
Assume from now on that
for all $\CI\in\CF$ and $\CF'\subset\CF-\{\CI\}$, then
$\bigcap_{\CI'\in\CF'}\CI'/\bigcap_{\CI'\in\CF'\cup\{\CI\}}\CI'$ is
compactly generated in $\CT/\CI$.

\begin{itemize}
\item
Then, $\CT$ is compactly generated
and an object of $\CT$ is compact if and only if it is compact in $\CT/\CI$
for all $\CI\in\CF$.

\item
Let $\CJ$ be a Bousfield subcategory of $\CT$ intersecting properly every 
element of $\CF$ and such that for all $\CI\in\CF$ and
$\CF'\subset\CF-\{\CI\}$, then
$\CJ\cap\bigcap_{\CI'\in\CF'}\CI'/\CJ\cap\bigcap_{\CI'\in\CF'\cup\{\CI\}}\CI'$
is compactly generated in $\CT/\CI$.
Then, $\CJ$ is compactly generated in $\CT$.
\end{itemize}
\end{thm}

\begin{proof}
We prove each assertion of the Theorem by induction on the cardinality of $\CF$.

Let $\CI_1\in\CF$. We put $\CI_2=\bigcap_{\CI\in\CF-\{\CI_1\}}\CI$ and
$\bar{\CT}=\CT/\CI_2$. Given $\CI\in\CF$, we put $\bar{\CI}=\CI/(\CI\cap\CI_2)$,
viewed as a full subcategory of $\bar{\CT}$. Let
$\bar{\CF}=\{\bar{\CI}\}_{\CI\in\CF-\{\CI_1\}}$.
We have canonical equivalences $\CT/\CI\iso \bar{\CT}/\bar{\CI}$ and
$\CT/(\CI\cap\CI')\iso \bar{\CT}/(\bar{\CI}\cap\bar{\CI'})$. 
This shows that $\bar{\CF}$ is a cocovering of $\bar{\CT}$.
Let $\tilde{\CT}=\CT/\langle\CI_1\cup\CI_2\rangle_\infty$.
Given $\CI\in\CF-\{\CI_1\}$, let $\tilde{\CI}=
\CI/\langle\CI_2\cup(\CI\cap\CI_1)\rangle_\infty$.
Let $\tilde{\CF}=\{\tilde{\CI}\}_{\CI\in\CF-\{\CI_1\}}$. This
is a cocovering of $\tilde{\CT}$.

\smallskip
Let $\CF'$ be a non-empty subset of $\CF-\{\CI_1\}$. We have equivalences
$\bar{\CT}/\langle\bigcup_{\CI\in\CF'}\bar{\CI}\rangle_\infty\simeq
\CT/\langle\bigcup_{\CI\in\CF'}\CI\rangle_\infty$ and
$\tilde{\CT}/\langle\bigcup_{\CI\in\CF'}\tilde{\CI}\rangle_\infty\simeq
\CT/\langle\bigcup_{\CI\in\CF'}\CI\cup\CI_1\rangle_\infty$.
Let $C\in\CT$ such that $C$ is compact in
$\CT/\langle\bigcup_{\CI\in\CF'}\CI\rangle_\infty$ for all non empty
$\CF'\subset\CF$. By induction, $C$ is compact in $\bar{\CT}$ and in
$\tilde{\CT}$. Since it is also compact in $\CT/\CI_1$, it follows from
Proposition \ref{localcompact} that $C$ is compact.

\smallskip
Given $\CF'\subset\CF-\{\CI_1\}$ and $\CI\in\CF-(\CF'\cup\{\CI_1\})$,
then we have a canonical equivalence
$\bigcap_{\CI'\in\CF'}\CI'/\bigcap_{\CI'\in\CF'\cup\{\CI\}}\CI'\iso
\bigcap_{\CI'\in\CF'}\bar{\CI}'/\bigcap_{\CI'\in\CF'\cup\{\CI\}}\bar{\CI}'$.
This shows that 
$\bigcap_{\CI'\in\CF'}\CI'/\bigcap_{\CI'\in\CF'\cup\{\CI\}}\CI'$ is
compactly generated in $\bar{\CT}/\bar{\CI}$.
By induction, we deduce that $\bar{\CT}$ is compactly generated.

The induction hypothesis shows
that $\CI_1$ is compactly generated in $\bar{\CT}$.
Now, by assumption, $\CT/\CI_1$ is compactly generated and $\CI_2$
is compactly generated in $\CT/\CI_1$. So, Proposition \ref{constructiongenerator}
shows that $\CT$ is compactly generated.

\smallskip
Consider now $\bar{\CJ}=\CJ/(\CJ\cap\CI_2)$. Then, $\bar{\CJ}$
intersects properly any $\bar{\CI}\in\bar{\CF}$. Also, 
$\bar{\CJ}/(\bar{\CJ}\cap\bar{\CI})$ is compactly generated in
$\bar{\CT}/\bar{\CI}$. By induction, we deduce that $\bar{\CJ}$ is compactly
generated in $\bar{\CT}$. Also,
$\CJ\cap\CI_1$ is compactly generated in $\bar{\CT}$. By assumption,
$\CJ/(\CJ\cap\CI_1)$ and $\CJ\cap\CI_2$ are compactly generated in
$\CT/\CI_1$.
It follows from Proposition \ref{constructiongenerator}
shows that $\CJ$ is compactly generated in $\CT$.

\smallskip
Let $C\in\CT$.
By induction, the image $\bar{C}$ of $C$ in $\bar{\CT}$ is compact if and only
if $C$ is compact in $\CT/\CI$ for $\CI\in\CF-\{\CI_1\}$. Now, Corollary
\ref{caractcompact} shows that $C$ is compact in $\CT$ if and only if it is compact
in $\bar{\CT}$ and in $\CT/\CI_1$ and we are done.
\end{proof}

Note that the proof of Theorem \ref{cocovering}
actually provides a construction of a generating set.
For example, if the generating sets in the hypotheses of the Theorem
are all finite, then $\CT$ is generated by a finite set of compact objects,
hence by a single compact object (and the same holds for $\CJ$).


\begin{prop}
Let $\CF$ be a cocovering of $\CT$.
Then,
$\dim\CT< \sum_{\CI\in\CF}(1+\dim\CT/\CI)$.
\end{prop}

\begin{proof}
As in the proof of Theorem \ref{cocovering}, we proceed by induction on the
cardinality of
$\CF$. We take $\CI_1\in\CF$ and put $\CI_2=\bigcap_{\CI\in\CF-\{\CI_1\}}\CI$.
By induction, we have
$\dim\CT/\CI_2< \sum_{\CI\in\CF-\{\CI_1\}}(1+\dim\CT/\CI)$.
On the other hand, we have an essentially surjective functor
$\CT/\CI_1\to\CI_2$, hence $\dim\CI_2\le \dim \CT/\CI_1$
(Lemma \ref{quotient}).
\end{proof}

Note that this holds as well for the other two definitions of dimension
of Remark \ref{otherdim} when the functor
$j_{\CI*}$ commutes with direct sums (then, $i_{\CI}^!$ commutes
with direct sums as well) --- the corresponding result is certainly
more interesting. This holds in the geometric setting of \S\ref{subsecschemes}.

\section{Derived categories of algebras and schemes}
\label{subsecalg}
We study here the concepts of \S \ref{secfin} for derived categories
of algebras and schemes.

\subsection{Algebras}
\subsubsection{}

From Theorem \ref{equivgen} (3), we deduce the following result
\cite[\S 5.3]{Ke}~:

\begin{cor}
\label{genDG}
Let $A$ be a dg algebra. Then,
$D(A)^c=\langle A\rangle_\infty$.
\end{cor}


\begin{prop}
\label{dgbounded}
Let $A$ be a dg algebra and $C\in D(A)$. Then,
$C$ is cohomologically locally bounded (resp. bounded above, resp. bounded below)
if and only if $H^i(C)=0$ for $|i|\gg 0$ (resp. for $i\gg 0$, resp. for $i\ll 0$).
In particular, if $C$ is cohomologically locally finitely generated, then
$C\in D^b(A)$.
\end{prop}

\begin{proof}
We have $D(A)^c=\langle A\rangle_\infty$ (Corollary \ref{genDG}). Hence,
$C$ is cohomologically locally bounded (resp. bounded above, resp. bounded
below) if and only if $h_C(A[i])=0$ for $|i|\gg 0$ (resp. for $i\ll 0$,
resp. for $i\gg 0$). Since $h_C(A[i])\iso H^{-i}(C)$, the result follows.
\end{proof}

\subsubsection{}
For $A$ an algebra, we denote by $K^{-,b}(A\mproj)$ (resp.
$K^{-,b}(A\mProj)$) the 
homotopy category of right bounded complexes of finitely
generated projective $A$-modules (resp. projective $A$-modules)
with bounded cohomology.

\begin{prop}
\label{catalg}
Let $A$ be an algebra. The canonical functors
induce equivalences between
\begin{itemize}
\item $K^b(A\mproj)$ and $D(A)^c$
\item
$K^{-,b}(A\mproj)$ and the full subcategory of $D(A)$ of cohomologically
locally finitely presented objects
\end{itemize}
\end{prop}

\begin{proof}
The first assertion is an immediate consequence of Corollary \ref{genDG}.

Recall that the canonical functor $K^{-,b}(A\mProj)\to D^b(A)$ is
an equivalence.

\smallskip
We now prove the second assertion.
Let $C\in D(A)$. By Corollary \ref{genDG} and Lemma \ref{Xthick},
$C$ is cohomologically locally finitely
presented if and only if conditions (a) and (b) hold for
$X=A$.

\smallskip
Let $C$ be a right bounded complex of finitely generated projective
$A$-modules with bounded cohomology. Consider $r$ such that $H^i(C)=0$ for
$i\le r$. 
The canonical
map from the stupid truncation $\sigma^{\ge r}C$ to $C$ is surjective on
cohomology, so $C$ satisfies (a), hence $C$ is cohomologically locally
finitely generated.
Now, Lemma \ref{carclfp} shows that $C$ is cohomologically locally
finitely presented.


\smallskip
Let $C$ be a cohomologically locally finitely generated object. Then,
$C$ has bounded cohomology (Proposition \ref{dgbounded}). 
Let $i$ be maximal such that
$H^i(C)\not=0$. Up to isomorphism, we can assume $C^j=0$ for $j>i$.
By assumption, there is a bounded complex $D$ of finitely generated projective
$A$-modules and $f:D\to C$ a morphism of complexes
such that $H(f)$ is onto. In particular, we have
a surjection $D^i\to C^i\to H^i(C)$, hence
$H^i(C)$ is finitely generated.

\smallskip
Let $C$ be  cohomologically locally finitely
presented.

Assume first $C=M$ is a complex concentrated in degree $0$.
Let $f:D^0\to M$ be a surjection, with $D^0$ finitely generated projective.
Then, $\ker f$ is cohomologically locally finitely presented
(Proposition \ref{clfpthick}), hence is
the quotient of a finitely generated projective module. By induction,
it follows that $M$ has a left resolution by finitely generated projective
$A$-modules.

%
%

\smallskip
We take now for $C$ an arbitrary  cohomologically locally finitely
presented object.
We know that $C$ has bounded cohomology and
we now prove by induction on
$\sup\{i|H^i(C)\not=0\}-\min\{i|H^i(C)\not=0\}$ that $C$ is isomorphic
to an object of $K^{-,b}(A\mproj)$.

Let $i$ be maximal such that
$H^i(C)\not=0$. As proven above, there is a finitely generated projective
$A$-module $P$ and a 
morphism of complexes $f:P[-i]\to C$ such that $H^i(f)$ is surjective.
Let $C'$ be the cone of $f$.
By Proposition \ref{clfpthick},
$C'$ is again cohomologically
locally finitely presented. By induction, $C'$ is isomorphic to an
object of $K^{-,b}(A\mproj)$ and we are done.
\end{proof}

\begin{cor}
\label{lfgnoetherian}
Let $A$ be a noetherian algebra. Then, the full subcategory of
cohomologically locally finitely presented objects of $\CT=D(A)$
is equivalent to $D^b(A\mMod)$.
\end{cor}

\begin{rem}
For a dg algebra, there might be no non-zero
cohomologically locally bounded objects (\eg, for $k[x,x^{-1}]$ with
$x$ in degree $1$ and differential zero).
The notion of cohomologically locally finitely presented
objects is more interesting for our purposes.
\end{rem}

\subsection{Schemes}
\label{secschemes}
\subsubsection{}
\label{subsecschemes}
Recall that a scheme is quasi-compact and quasi-separated
if it has a finite covering $\CC$
by open affine subschemes such that given $U,U'\in\CC$, then $U-(U\cap U')$ is 
a closed subscheme of $U$ defined by a finite number of equations.

Let $X$ be a   quasi-compact and quasi-separated scheme. The category
$D(X)$ is a cocomplete triangulated
category. The perfect complexes have bounded cohomology.
If $X$ is in addition separated, then the canonical functor
$D(X\mqcoh)\to D(X)$ is an equivalence \cite[Corollary 5.5]{BoNee}.
If $X=\Spec R$, then $D(X)\simeq D(R\mMOD)$. If $X$ is a noetherian
scheme, then it is quasi-compact and quasi-separated and we denote by
$D_{\coh}(X)$ the full subcategory of $D(X)$ of complexes with
coherent cohomology sheaves.

\smallskip
Let $U$ be a quasi-compact open subscheme of $X$ (\ie, a finite union of affine
open subschemes). We denote by $D_{X-U}(X)$ the full subcategory of
$D(X)$ of complexes with cohomology supported by $X-U$. 
We denote by $j:U\to X$ the open immersion and $i:X-U\to X$ the closed
immersion.
We have an exact sequence of triangulated categories
$$0\to D_{X-U}(X)\xrightarrow{i_*} D(X)\xrightarrow{j^*}D(U)\to 0$$
and adjoint pairs $(i_*,i^!)$ and $(j^*,j_*)$. In particular,
$D_{X-U}(X)$ is a Bousfield subcategory of $D(X)$.
Furthermore, $j_*$ has finite cohomological dimension
(\ie, there is an integer $N$ such that if
$C\in D(U)$ and $H^n(C)=0$ for $n>0$, then $H^n(j_*C)=0$ for $n\ge N$).
Consequently, $i^!$ has also finite cohomological dimension.

Given $U$ and $U'$ two quasi-compact open subschemes of $X$, then
$D_{X-U}(X)$ and $D_{X-U'}(X)$ intersect properly and
$D_{X-U}(X)\cap D_{X-U'}(X)=D_{X-(U\cup U')}(X)$. If
$U\cup U'=X$, then the restriction functor
$D_{X-U'}(X)\iso D_{U-U'\cap U}(U)$ is an equivalence.

Given $\CF$ a finite family of open subschemes of $X$, then
$\{D_{X-U}(X)\}_{U\in\CF}$ is a cocovering of $D(X)$ if and only if $\CF$
is a covering of $X$.

\subsubsection{}
Let us start with the study of the affine case.

\smallskip
The following Proposition makes \cite[Proposition 6.1]{BoNee}  more
precise.

\begin{prop}
\label{BoNee}
Let $A$ be a commutative ring and $f_1,\ldots,f_n\in A$. Let
$I$ be the ideal of $A$ generated by $f_1,\ldots,f_n$.
Let $X=\Spec A$ and $Z=\Spec A/I$.

Let $K(f_1,\ldots,f_n)=\bigotimes_i (0\to \CO_X\xrightarrow{f_i}\CO_X\to 0)$ 
(with non zero terms in degrees $-n,\ldots,0$).

Then,
\begin{itemize}
\item Let $C\in D_Z(X)$ such that $H^0(C)\not=0$. Then, 
$\Hom_{D(X)}(K(f_1,\ldots,f_n),C)\not=0$. Given $\phi\in R\Gamma^0(C)$,
there are
integers $d_1,\ldots,d_n>0$ such that $\phi$ is in the image
of the canonical map $\Hom_{D(X)}(K(f_1^{d_1},\ldots,f_n^{d_n}),C)\to
\Hom_{D(X)}(\CO_X,C)=R\Gamma^0(C)$.
\item $K(f_1,\ldots,f_n)$ is a compact object of $D(X)$ that is a generator for
$D_Z(X)$.
\end{itemize}
\end{prop}

\begin{proof}
It is clear that $K(f_1,\ldots,f_m)$ is compact and supported by $Z$. Also, the first
statement of the Proposition implies the second one.

We have a distinguished triangle
$$K(f_1,\ldots,f_{m-1})\xrightarrow{f_m} K(f_1,\ldots,f_{m-1})\to
K(f_1,\ldots,f_m)\rightsquigarrow$$
giving an exact sequence
$$\Hom(K(f_1,\ldots,f_m),C)\to \Hom(K(f_1,\ldots,f_{m-1}),C)
\xrightarrow{f_m}\Hom(K(f_1,\ldots,f_{m-1}),C).$$

We prove the first assertion by induction on $m$. 
Since $\Hom(K(f_1,\ldots,f_{m-1}),C)$ is supported by $Z$ and non zero,
it follows that the kernel of the multiplication by $f_m$ is not zero,
hence $\Hom(K(f_1,\ldots,f_m),C)\not=0$.
By induction, there exists $d_1,\ldots,d_{m-1}>0$ and $\phi_{m-1}\in
\Hom(K(f_1^{d_1},\ldots,f_{m-1}^{d_{m-1}}),C)$ with image
$\phi\in R\Gamma^0(C)$.
There is $d_m>0$ such that $f_m^{d_m}\phi_{m-1}=0$. Then, there
is $\phi_m\in \Hom(K(f_1^{d_1},\ldots,f_m^{d_m}),C)$ with image $\phi_{m-1}$.
Now, $\phi_m$ has image $\phi\in R\Gamma^0(C)$.
\end{proof}

\begin{lemma}
\label{affineapprox}
Let $X=\Spec A$ be an affine scheme and $Z$ a closed subscheme defined
by $f_1=\cdots=f_n=0$. Let $\CK$ be the smallest additive subcategory
of $D_Z(X)$ containing the objects $K(f_1^{d_1},\ldots,f_n^{d_n})$
for $d_1,\ldots,d_n>0$.

Let $a\le b$ be two integers.
Let $C\in D_Z^{\le b}(X)$ with $C^i$ a vector bundle for $i\ge a$.
Then, there is $P\in \CK[-b]\ast\CK[-b+1]\ast\cdots\ast\CK[-a]$ and
$f:P\to C$ such that $H^i(\cone(f))=0$ for $i\ge a$.
\end{lemma}

\begin{proof}
We prove the Lemma by induction on $b-a$. By assumption,
$H^b(X)$ has finite type.
It follows from Proposition \ref{BoNee} that there is 
$K_1\in \CK$ and $f_1:K_1[-b]\to X$ such that $H^b(f_1)$ is surjective.
Let $C'=\cone(f_1)$. By induction, there is
$L\in \CK[-b+1]\ast\cdots\ast\CK[-a]$ and $g:L\to C'$ such that
$H^i(\cone(g))=0$ for $i\ge a$. Let 
$P$ be the cocone of the composition
$L\to C'\to K_1[-b+1]$. There is $f:P\to C$ making the following
diagram commutative
$$\xymatrix{
K_1[-b]\ar[r] & C\ar[r] & C'\ar[r] & K_1[-b+1] \\
K_1[-b]\ar[r]\ar@{=}[u] & P\ar[r]\ar@{.>}[u] & L\ar[r]\ar[u] &
 K_1[-b+1]\ar@{=}[u]
}$$
Since $\cone(f)\simeq\cone(g)$, we are done.
\end{proof}

\subsubsection{}
The following result is classical, although no published proof seems to exist
(when $Z=X$, cf \cite[Corollary 2.3 and Proposition 2.5]{Nee3} for the
separated case and \cite[Theorem 3.1.1]{BoVdB} for the general case).
The general constructions of \S \ref{seccoverings} reduce immediately its
proof to the affine case.

\begin{thm}
\label{genscheme}
Let $X$ be a  quasi-compact and quasi-separated scheme.
The perfect complexes on $X$ are the compact objects of $D(X)$.

Let $Z$ be a closed subscheme of $X$
with $X-Z$ quasi-compact. Then,
$D_Z(X)$ is generated by an object of
$D_Z(X)\cap D(X)^c$.
\end{thm}

\begin{proof}
Theorem \ref{cocovering} shows that compactness is of local nature in
the following sense~: an object $C\in D(X)$
is compact if and only there is a finite covering $\CC$ of
$X$ by quasi-compact open subschemes
such that the restriction of $C$ to an intersection of open subschemes
in $\CC$ is compact. Perfectness is obviously also of local nature, in
that sense. Since $X$ is  quasi-compact and quasi-separated,
we can even assume that the open
subschemes in the coverings are affine. This shows that compact complexes
are perfect.

Let us prove that perfect complexes are compact. The discussion above
reduces the problem to proving that bounded complexes of vector bundles
are compact.
Corollary \ref{genDG} shows that a bounded complex of vector bundles over
an affine scheme is compact. The discussion above allows us to deduce
that the same remains true for quasi-compact separated schemes, and then
for  quasi-compact and quasi-separated schemes.

\smallskip
The scheme $X$ has a finite covering $\CC$ by affine open subschemes with
$U-(U\cap U')$ defined by a finite number of equations for any $U,U'\in\CC$.
Theorem \ref{cocovering} reduces then the second part of the Theorem to the
case where $X$ is affine. If $Z$ is defined by the equations
$f_1=\cdots=f_n=0$, then 
$\bigotimes_i (0\to \CO_X\xrightarrow{f_i}\CO_X\to 0)$ is a generator
of $D_Z(X)$ that is compact in $D(X)$ (Proposition \ref{BoNee}).
\end{proof}
 
Note that we deduce $D_Z(X)^c=D_Z(X)\cap D(X)^c$ (Theorem
\ref{locNeeman} (i)).

\subsubsection{}
Given $C,D\in D_Z(X)$, we denote by $\amp(C)$ (resp. $\amp_D(C)$)
the smallest interval of
$\BZ$ such that $H^i(C)=0$ for $i\not\in\amp(C)$ (resp. $\Hom(D,C[i])=0$
for $i\not\in\amp_D(C)$). Given $I$ an interval of $\BZ$ and $m\ge 0$,
we put $I\pm m=\{i+j\}_{i\in I,j\in\BZ\cap [-m,m]}$.

The following Proposition relates in a precise way boundedness of a
complex and cohomological local boundedness
(cf \cite[Lemma 3.3.8]{BoVdB} for bounded cohomology implies
cohomologically locally bounded).

\begin{prop}
\label{amplitude}
Let $X$ be a  quasi-compact and quasi-separated scheme
and $Z$ be a closed subscheme of $X$ with $X-Z$ quasi-compact. 
Let $C\in D_Z(X)$. Then,
$C$ is cohomologically locally bounded (resp. bounded above, resp.
bounded below)
if and only if $H^i(C)=0$ for $|i|\gg 0$ (resp. for $i\gg 0$, resp.
for $i\ll 0$).

More precisely, let $G\in D(X)^c\cap D_Z(X)$ be a generator for $D_Z(X)$.
Then, there is an integer $N$ such that for any $C\in D_Z(X)$, then
$\amp_G(C)\subset \amp(C)\pm N$ and
$\amp(C)\subset \amp_G(C)\pm N$.
\end{prop}
 
\begin{proof}
Let $G'\in D_Z(X)^c$. Then, there is an
integer $d$ such that $G'\in\langle G\rangle_d$.
As a consequence, there is an integer
$m$ such that for any $C\in D_Z(X)$, then
$\amp_{G'}(C)\subset \amp_G(C)\pm m$.
Note that this shows it is enough to prove the more precise statements for one 
$G$.

\smallskip
Let us first assume that $X$ is affine. We take $G$ as in Proposition
\ref{BoNee}. Let $C\in D_Z(X)$. If $\Hom(G,C)=0$, then, $H^0(C)=0$.
Conversely, If $H^i(C)=0$ for $-n\le i\le 0$, then $\Hom(G,C)=0$.
So, the Proposition follows.

\smallskip
We denote by $\tilde{\i}:Z\to X$ the closed immersion.
Let $m$ be an integer such that for $C\in D(X)$, we have
$\amp(\tilde{\i}^!C)\subset \amp(C)\pm m$.
Let $G_0$ be a compact generator of $D(X)$. Since $\tilde{\i}_*G_0$
is compact and $G$ is a classical generator of $D(X)^c$, there is
an integer $m'$ such that for $C\in D(X)$, we have
$\amp_{\tilde{\i}_*G}(C)\subset\amp_{G_0}(C)\pm m'$.

Let $D\in D(X)$.
We have $\Hom(G,\tilde{\i}^!D)\simeq \Hom(\tilde{\i}_*G,D)$, so
$\amp_G(\tilde{\i}^!D)=\amp_{\tilde{\i}_*G}(D)\subset \amp_{G_0}(D)\pm m'$.

Let $C\in D_Z(X)$. Then, 
$\amp_G(C)\subset \amp_{G_0}(\tilde{\i}_*C)\pm m'$.
Since $\amp(\tilde{\i}_*C)=\amp(C)$, it follows that it is enough to
prove the first inclusion of the Proposition in the case where $Z=\emptyset$.
By induction, the Mayer-Vietoris triangle (Proposition \ref{mayervietoris} (2))
reduces the proof to the affine case, which we already considered.

\smallskip
Let $U_1,\ldots,U_n$ be an affine open covering of $X$.
We have canonical equivalences
$D_{Z-U_r\cap Z}(X)\iso D_{Z-U_r\cap Z}(\bigcup_{s\not=r} U_s)$ and
$D_{U_r\cap Z}(X-(Z-(U_r\cap Z)))\iso D_{U_r\cap Z}(U_r)$. So, we have
an exact sequence of triangulated categories
$$0\to D_{Z-U_r\cap Z}(\bigcup_{s\not=r} U_s)\xrightarrow{i_*}
D_Z(X)\xrightarrow{j^*} D_{U_r\cap Z}(U_r)\to 0$$
and an exact triangle of functors
$i_*i^!\to \Id_{D_Z(X)}\to j_*j^*\rightsquigarrow.$

We now show the second inclusion by induction on $n$.
Let $H$ be a compact generator of $D_{Z-U_r\cap Z}(\bigcup_{s\not=r} U_s)$.
By induction, there is an integer $N_1$ such that for every
$C'\in D_{Z-U_r\cap Z}(\bigcup_{s\not=r} U_s)$, we have
$\amp(C')\subset \amp_H(C')\pm N_1$.
Given $C\in D_Z(X)$, we have
$\Hom(H,i^!C)\simeq \Hom(i_*H,C)$. There is an integer $N_2$ such that
for any $C\in D_Z(X)$, we have $\amp_{i_*H}(C)\subset\amp_G(C)\pm N_2$.
So, given $C\in D_Z(C)$, we have
$\amp(i^!C)\subset \amp_G(C)\pm N_1\pm N_2$.
The proof above shows that there is $N_3$ such that given any 
$D\in D_Z(X)$, we have
$\amp_G(D)\subset \amp(D)\pm N_3$. In particular, for any 
$C\in D_Z(X)$, we have $\amp_G(i_*i^!C)\subset\amp(i^!C)\pm N_3\subset
\amp_G(C)\pm N_1\pm N_2\pm N_3$, hence
$\amp_G(j_*j^*C)\subset \amp_G(C)\pm N_1\pm N_2\pm N_3\pm 1$.

The study of the affine case shows there is $N_4$ such that
for any $C\in D_Z(X)$, then $\amp(j^*C)\subset\amp_{j^*G}(j^*C)\pm N_4=
\amp_G(j_*j^*C)\pm N_4\subset \amp_G(C)\pm N_1\pm N_2\pm N_3\pm N_4\pm 1$.
There is an integer $N_5$ such that for any $D\in D_{U_r\cap Z}(U_r)$, we have
$\amp(j_*D)\subset\amp(D)\pm N_5$.
Since $\amp(C)\subset\amp(i_*i^!C)\cup \amp(j_*j^*C)$, we deduce that
$\amp(C)\subset\amp_G(C)\pm N_1\pm N_2\pm N_3\pm N_4\pm N_5\pm 1$.
\end{proof}

\subsubsection{}
An object $C\in D(X)$ is {\em pseudo-coherent} if for every $a\in\BZ$ and
every point $x$ of
$X$, there is an open subscheme $U$ of $X$ containing $x$, a bounded complex
$D$ of vector bundles on $U$ and $f\in\Hom_{D(U)}(D,C_{|U})$ such
that $H^i(\cone(f))=0$ for $i\ge a$ (cf \cite[\S I.2]{SGA6} or
\cite[\S 2.2]{ThTr}). Pseudo-coherent complexes form a thick
subcategory of $D^-(X)$.

The following Proposition gives a substitute for global resolutions of
pseudo-coherent complexes. Such resolutions exist for 
schemes with a family of ample line bundles (cf \cite[\S II]{SGA6} or
\cite[Proposition 2.3.1]{ThTr}). It shows that pseudo-coherence of $C$ is
a condition on the functor $\Hom(-,C)$ restricted to compact objects.

\begin{prop}
\label{pseudocoherent}
Let $X$ be a quasi-compact and quasi-separated scheme and $Z$ a closed
subscheme with $X-Z$ quasi-compact.
Let $C\in D_Z(X)$. The following conditions are equivalent
\begin{itemize}
\item[(i)] $C$ is pseudo-coherent
\item[(ii)] given $a\in\BZ$, there is $D\in D_Z(X)^c$ and $f:D\to C$ such that
$H^i(\cone(f))=0$ for $i\ge a$
\item[(iii)]
given any $G\in D_Z(X)^c$ and any $a\in\BZ$, there is $D\in D_Z(X)^c$ and
$f:D\to C$ such that $\Hom(G,\cone(f)[i])=0$ for $i\ge a$.
\end{itemize}
\end{prop}

\begin{proof}
We prove (i)$\Rightarrow$(ii) by induction on the minimal number of affine
open subschemes in a covering of $X$. Let $X=U\cup V$ where $U$ is an
open affine subscheme and $V$ is a an open subscheme that can be covered
by strictly less affine open subschemes than $X$. 
Let $n$ be the minimal number of defining equations of $Z\cap (X-V)$
as a closed subscheme of $U$.

Let $C\in D_Z(X)$
be pseudo-coherent and let $a\in\BZ$. Then, $C_{|V}$ is pseudo-coherent and by
induction there is $D_1\in D_{Z\cap V}(V)^c$ and $f_1:D_1\to C_{|V}$ such that
$H^i(\cone(f_1))=0$ for $i\ge a-n$. Replacing $D_1$ by
$D_1\oplus D_1[d]$ and $f_1$ by $(f_1,0)$ for $d\gg 0$ odd, we can
assume in addition that $[D_1]=0$. Then, Theorem \ref{locNeeman} shows
that $f_1$ lifts to $f'_1:D'_1\to C$ where $D'_1\in D_Z(X)^c$.
Let $C_1=\cone(f'_1)$. Let $C_2=\tau^{\ge a-n}C_1$, an object of
$D_{Z\cap (X-V)}(X)$.
Lemma \ref{affineapprox} shows there is $D_2\in D_{Z\cap (X-V)}(U)$ a
bounded complex of free $\CO_U$-modules of finite type with
$D_2^i=0$ for $i<a-n$ and a map $f_2:D_2\to C_{2|U}$ such that
$H^i(\cone(f_2))=0$ for $i\ge a$. Via the equivalence
$D_{Z\cap (X-V)}(X)\iso D_{Z\cap (X-V)}(U)$, this map
corresponds to $f'_2:D'_2\to C_2$ with $D'_2\in D_{Z\cap (X-V)}(X)^c$.
We have 
$$\Hom(D'_2,(\tau^{<a-n}C_1)[1])\simeq
\Hom(D_2,(\tau^{<a-n}C_1)[1]_{|U})=0$$
hence there is $f_3:D'_2\to C_1$ lifting $f'_2$. Let $C_3$ be its cone.
We have a distinguished triangle
$\tau^{<a-n}C_1\to C_3\to \cone(f'_2)\rightsquigarrow$, hence
$H^i(C_3)=0$ for $i\ge a$.
Let $D$ be the cocone of the composition
$C\to C_1\to C_3$.
The octahedral axiom shows that $D\in D_Z(X)^c$ and we are done.
$$
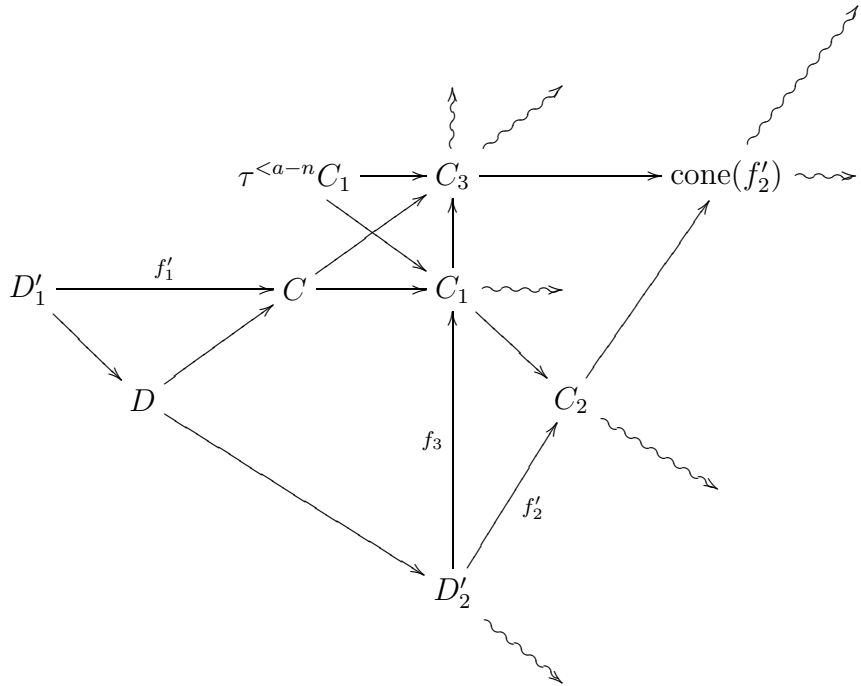
\xymatrix{
&&&&&&\\
&&&&&&\\
&&\tau^{<a-n}C_1\ar[r]\ar[dr]&C_3\ar[rr]\ar@{~>}[ur]\ar@{~>}[u]&&\cone(f'_2)
\ar@{~>}[r]\ar@{~>}[uur] & \\
D'_1\ar[rr]^-{f'_1}\ar[dr] &&C\ar[ur]\ar[r]&C_1\ar[u]\ar[dr]\ar@{~>}[r]&\\
&D\ar[ur]\ar[ddrr]&&&C_2\ar[uur]\ar@{~>}[dr] \\
&&&&&&\\
&&&D'_2\ar[uur]_-{f'_2}\ar[uuu]^{f_3}\ar@{~>}[dr]\\
&&&&&&\\
}$$

\smallskip
Since compact objects of $D_Z(X)$ are isomorphic, on affine open subschemes,
to bounded complexes of vector bundles, we have (ii)$\Rightarrow$(i).
The equivalence between (ii) and (iii) is given by Proposition \ref{amplitude}.
\end{proof}

We say that a noetherian scheme $X$ satisfies $(*)$ if given $G$ a compact
generator of $D(X)$ and given any $M\in X\mqcoh$,
there is $C\in \langle \tilde{G}\rangle_\infty$
and $f:C\to M$ such that $H^0(f)$ is surjective.

Note that if condition $(*)$ holds for one $G$, then it holds for all 
compact generators (cf Theorem \ref{equivgen} (3)).

\begin{lemma}
Let $X$ be an affine scheme or a quasi-projective scheme over a field.
Then, $X$ satisfies $(*)$.
\end{lemma}

\begin{proof}
The affine case is clear for $G=\CO_X$. The other case is solved by
Lemma \ref{powers} below.
\end{proof}

\begin{prop}
\label{clfpcoh}
Let $X$ be a noetherian scheme satisfying $(*)$. Then, 
the full subcategory of cohomologically locally finitely presented objects
of $D(X)$ is equivalent to $D^b_{\coh}(X)$.
\end{prop}

\begin{proof}
Let $G$ be a compact generator for $D(X)$.

Let $M\in D_{\coh}^b(X)$. Then, $M$ is cohomologically locally bounded
(Proposition \ref{amplitude}). Take $a\in\BZ$ such that
$\Hom(G,M[i])=0$ for $i<a$.
Consider $N$ as in Proposition \ref{amplitude}.
By Proposition \ref{pseudocoherent}, there is $C\in D(X)^c$ and
$f:C\to M$ such that $H^i(\cone(f))=0$ for $i\ge a-N$.
Then, $\Hom(G[i],f)$ is surjective for all $i$. It follows that
$M$ is cohomologically  locally finitely generated. So, every object
of $D_{\coh}^b(X)$ is cohomologically locally finitely presented
(Lemma \ref{carclfp}).

\smallskip
Let $C$ be a cohomologically locally finitely presented object.
Thanks to Proposition \ref{amplitude}, we know that $C$ has bounded cohomology.
Assume $C\not\in D^b_{\coh}(X)$ and take $i$ minimal such that
$H^i(C)$ is not coherent. Since $\tau^{<i}C\in D^b_{\coh}(X)$, it
follows from the first part of the Proposition that 
$\tau^{<i}C$ is cohomologically locally finitely presented, hence
$D=\tau^{\ge i}C$ is cohomologically locally finitely presented as well
(Proposition \ref{clfpthick}).
Condition $(*)$ shows that there is $E\in\langle
\tilde{G}\rangle_\infty$ and $f:E\to D$ such that
$H^i(f)$ is surjective. Lemma \ref{clfpfactor} shows that $f$ factors through 
a compact object $F$. Now, $H^i(F)$ is coherent, hence $H^i(D)$ is
coherent as well, a contradiction.
\end{proof}

\begin{rem}
Let $X$ be a quasi-compact quasi-separated scheme. Let us show that
given $M\in X\mqcoh$ of
finite type,
there is $C\in D(X)^c$ and $f:C\to M$ such that $H^0(f)$ is surjective.

Let $\CF$ be a finite covering of $X$ by affine open subschemes.
Given $U\in \CF$, there is a complex $C_U\in D(U)^c$ with $[C_U]=0$ and a map
$f_U:C_U\to \CF_{|U}$ such that $H^0(f)$ is onto. By Theorem
\ref{locNeeman}, there is $C(U)\in D(X)^c$, $\phi_U:C(U)_{|U}\iso C_U$ and
$f(U):C(U)\to \CF$ such that $f(U)_{|U}=f_U\phi_U$. Let
$C=\bigoplus_{U\in\CF}C(U)$ and $f=\sum f(U)$. Then, $C\in D(X)^c$
and $f$ is surjective.
\end{rem}

\begin{rem}
\label{resolproj}
Let $k$ be a field and $X$ a projective scheme over $k$. Let $A$ be a 
dg algebra such that $D(A)\simeq D(X)$. Then, $H^*(A)$
is finite dimensional. Given $C\in D(A)$ with $H^*(C)$ finite dimensional
and given $a\in\BZ$,
there is $D\in D(A)^c$ and $f:D\to C$ such that $H^i(\cone(f))=0$ for
$i>a$. This is a very strong condition on a dg algebra. For example,
the dg algebra $k[x]/x^2$ with $x$ in degree $1$ and differential zero
doesn't satisfy this condition.
\end{rem}

\subsection{Compact objects in bounded derived categories}
\label{seccompactbounded}
\begin{prop}
\label{compactsinj}
Let $\CA$ be an abelian category
with exact filtered colimits and a set $\CG$ of generators (\ie, a
Grothendieck category). Assume that
for any $G\in\CG$, the subobjects of $G$ are compact.

Then,
$(D^b(\CA))^c=\langle\CA^c\rangle_\infty$.
\end{prop}

\begin{proof}
An object $I$ of $\CA$ is injective if and only if for any $G\in\CG$ and any
subobject $G'$ of $G$, the canonical map
$\Hom_\CA(G,I)\to\Hom_\CA(G',I)$ is surjective
\cite[Proposition V.2.9]{Ste}.
Note that $G'$ is compact. It follows that a direct sum of injectives is
injective.

Let $M\in\CA^c$.
Let $\CF$ be a family of objects of $D^b(\CA)$.
Then, $\bigoplus_{F\in\CF}F$ exists in $D^b(\CA)$ if and only if
the direct sum, computed in $D(\CA)$, has bounded cohomology, \ie, if and
only if, there are integers $r$ and $s$ such that for any $F\in\CF$, we
have $H^i(F)=0$ for $i<r$ and for $i>s$.
Given $F\in\CF$, let $I_F$ be a complex of injectives quasi-isomorphic to $F$
with zero terms in degrees less than $r$. 
Since $\bigoplus_F I_F^j$ is injective, we have
$\Ext^i(M,\bigoplus_F I_F^j)=0$ for all $j$ and $i>0$.
Hence, 
\begin{align*}
\bigoplus_F \Hom_{D(\CA)}(M,F)&\iso
\bigoplus_F H^0\Hom^\bullet_\CA(M,I_F)\iso
H^0\bigoplus_F \Hom^\bullet_\CA(M,I_F)\iso
H^0\Hom^\bullet_\CA(M,\bigoplus_F I_F)\\
&\iso\Hom_{D(\CA)}(M,\bigoplus_F F).
\end{align*}
It follows that $M\in D^b(\CA)^c$.

\smallskip
Let $C\in D^b(\CA)^c$. We prove by induction on
$\max\{i|H^iC\not=0\}- \min\{i|H^iC\not=0\}$ that
$C\in \langle\CA^c\rangle_\infty$.

Take $i$ maximal such that $H^iC\not=0$. Then,
$\Hom_{D^b(\CA)}(C,M[-i])\iso \Hom_\CA(H^iC,M)$ for any $M\in\CA$.
It follows that $H^iC\in\CA^c$. As proven above, we deduce that
$H^iC[-i]\in D^b(\CA)^c$, hence $\tau^{\le i-1}C\in D^b(\CA)^c$.
By induction, $\tau^{\le i-1}C\in \langle\CA^c\rangle_\infty$ and we are done.
\end{proof}

\begin{cor}
\label{compactinb}
Let $A$ be a noetherian ring. Then,
$D^b(A\mMod)\iso D^b(A)^c$.

Let $X$ be a separated noetherian scheme. Then,
$D^b_{\coh}(X)\iso D^b(X)^c$.
\end{cor}

\begin{proof}
In the ring case, we take $\CG=\{A\}$.
In the geometric case, we take for $\CG$ the set of coherent sheaves,
cf \cite[Appendix B, \S 3]{ThTr}.
\end{proof}

\section{Dimension for derived categories of rings and schemes}
\label{secder}
\subsection{Resolution of the diagonal}
\label{secresoldiago}
Let $k$ be a field.
\subsubsection{}

\begin{lemma}
\label{boundpdim}
Let $A$ be a noetherian $k$-algebra such that
$\pdim_{A^\en}A<\infty$. Then,
$D^b(A)=\langle \widetilde{A}\rangle_{1+\pdim_{A^\en}A}$ and
$D^b(A\mMod)=\langle A\rangle_{1+\pdim_{A^\en}A}$. In particular,
$\dim D^b(A\mMod)\le \pdim_{A^\en}A$.
\end{lemma}

\begin{proof}
The discussion in \S \ref{devissage}, shows that
$D^b(A)=\langle\widetilde{A}\rangle_{1+\pdim_{A^\en}A}$.
Now, we have $D^b(A\mMod)\simeq D^b(A)^c$ (Corollary \ref{compactinb})
and the result follows from Corollary \ref{intercompact}.
\end{proof}

We say that a commutative $k$-algebra $A$
is {\em essentially of finite type} if it is the localization of a
commutative $k$-algebra of finite type over $k$.

Recall the following classical result~:
\begin{lemma}
\label{diago}
Let $A$ be a finite dimensional $k$-algebra
or a commutative $k$-algebra essentially of finite type.
Assume that given $V$ a simple $A$-module, then $Z(\End_A(V))$ is
a separable extension of $k$.
Then, $\pdim_{A^\en}A=\gldim A$.
\end{lemma}

\begin{proof}
Note that under the assumptions, $A^\en$ is noetherian.
In the commutative case,
$\gldim A=\sup\{\gldim A_\Gm\}_\Gm$ and
$\pdim_{A^\en}A=\sup\{\pdim_{(A_\Gm)^\en} A_\Gm\}_\Gm$
where $\Gm$ runs over the maximal ideals of $A$. It follows that it
is enough to prove the commutative case of the Lemma for $A$ local.

So, let us assume now $A$ is finite dimensional or is a commutative
local $k$-algebra essentially of finite type.

Let $0\to P^{-r}\to\cdots\to P^0\to A\to 0$ be a minimal projective
resolution of $A$ as an $A^\en$-module.
So, there is a simple $A^\en$-module $U$ with $\Ext^r_{A^\en}(A,U)\not=0$.
The simple module $U$ is isomorphic to a quotient of
$\Hom_k(S,T)$ for $S,T$ two simple $A$-modules. By assumption,
$\End_A(S)\otimes_k \End_A(T)^\circ$ is semi-simple, hence
$U$ is actually isomorphic to a direct summand of
$\Hom_k(S,T)$.

Then,
$$\Ext^r_A(T,S)\iso \Ext^r_{A^\en}(A,\Hom_k(T,S))\not=0,$$
hence, $r\le\gldim A$.

Now, given $N$ an $A$-module,
$0\to P^{-r}\otimes_A N\to\cdots\to P^0\otimes_A N\to N\to 0$
is a projective resolution of $N$, hence $r\ge\gldim A$, so
$r=\gldim A$.
\end{proof}

\begin{rem}
Note that this Lemma doesn't hold if
the residue fields of $A$ are not separable extensions of $k$.
Cf the case $A=k'$ a purely inseparable extension of $k$.
\end{rem}

Combining Lemmas \ref{boundpdim} and \ref{diago}, we get

\begin{prop}
\label{boundgldimk}
Let $A$ be a finite dimensional $k$-algebra
or a commutative $k$-algebra essentially of finite type.
Assume that given $V$ a simple $A$-module, then $Z(\End_A(V))$ is
a separable extension of $k$.

If $A$ has finite global dimension, then 
$D^b(A)=\langle \widetilde{A}\rangle_{1+\gldim A}$ and
$D^b(A\mMod)=\langle A\rangle_{1+\gldim A}$. In particular,
$\dim D^b(A\mMod)\le \gldim A$.
\end{prop}

\begin{rem}
The dimension of $D^b(A\mMod)$ can be strictly less than $\gldim A$
(this will be the case for example for a finite dimensional $k$-algebra $A$
which is not hereditary but which is derived equivalent to a hereditary
algebra).
This cannot happen if $A$ is a finitely generated commutative 
$k$-algebra, cf Proposition \ref{inegvar} below.
\end{rem}

\subsubsection{}
\label{diagovar}

Following \S \ref{devissage}, we have the following result (cf
\cite[\S 3.4]{BoVdB}).

\begin{prop}
\label{resoldiagvar}
Let $X$ be a separated noetherian scheme over $k$. Assume
there is a vector bundle $\CL$ on $X$ and
a resolution of the structure sheaf $\CO_\Delta$ of the diagonal in
$X\times X$
$$0\to \CF^{-r}\to\cdots\to\CF^0\to \CO_\Delta\to 0$$
with $\CF^i\in\add(\CL\boxtimes\CL)$.

Then, $D^b(X\mqcoh)=\langle\widetilde{\CL}\rangle_{1+r}$
and $D^b(X\mcoh)=\langle \CL\rangle_{1+r}$.
\end{prop}

\begin{proof}
Let $p_1,p_2:X\times X\to X$ be the first and second projections.
For $C\in D^b(X\mqcoh)$, we have
$C\simeq Rp_{1*}(\CO_\Delta\otimes^\BL p_2^* C)$.
It follows that
$C\in\langle \CL\otimes_k R\Gamma(\CL\otimes C)\rangle_{1+r}$, hence
$C\in\langle \widetilde{\CL}\rangle_{1+r}$.
Since $D^b(X\mqcoh)^c=D^b(X\mcoh)$ (Corollary \ref{compactinb}),
the second assertion follows from Corollary \ref{intercompact}.
\end{proof}

Note that the assumption of the Proposition forces $X$ to be smooth.

\begin{example}
\label{Pn}
Let $X=\BP^n_k$. Let us recall results of Beilinson
\cite{Bei}. The object $G=\CO\oplus\cdots\oplus\CO(n)$ is a classical
generator for $D^b(X\mcoh)$. We have $\Ext^i(G,G)=0$ for $i\not=0$.
Let $A=\End(G)$. We have $D^b(X\mcoh)\simeq D^b(A\mMod)$.
We have $\gldim A=n$, hence $D^b(A\mMod)=\langle A\rangle_{n+1}$
(Proposition \ref{boundgldimk}), so
$D^b(\BP^n\mcoh)=\langle\CO\oplus\cdots\oplus\CO(n)\rangle_{n+1}$.
Another way to see this is to use the resolution of the diagonal
$\Delta\subset X\times X$~:
$$0\to \CO(-n)\boxtimes \Omega^n(n)\to\cdots\to\CO(-1)\boxtimes\Omega^1(1)
\to\CO\boxtimes\CO\to \CO_{\Delta}\to 0.$$
By Proposition \ref{inegvar} below,
it follows that $\dim D^b(\BP^n\mcoh)=n$.
\end{example}

\begin{example}
In \cite{Kap}, Kapranov considers flag varieties (type $A$) and
smooth projective quadrics. For these varieties $X$, he
constructs explicit bounded resolutions of the diagonal
whose terms are direct sums of $\CL\boxtimes\CL'$, where $\CL$
and $\CL'$ are vector bundles. It turns out that these resolutions have
exactly $1+\dim X$ terms (this is the smallest possible number). By
Proposition \ref{inegvar}, it follows that $\dim D^b(X\mcoh)=\dim X$.

Starting from a smooth projective variety $X$, there is an ample
line bundle whose homogeneous coordinate ring is a Koszul algebra
\cite[Theorem 2]{Ba}. This provides a resolution of
$\CO_{\Delta}$ \cite[Theorem 3.2]{Kaw}.
Now, if the kernel of the $r$-th map of the resolution is
a direct sum of sheaves of the form $\CL\boxtimes\CL'$, where $\CL,\CL'$
are vector bundles, then $\dim D^b(X\mcoh)\le r$.
Note that this can work only if the class of $\CO_\Delta$ is in the image of the
product map $K_0(X)\times K_0(X)\to K_0(X\times X)$. The case of flag
varieties associated to reductive groups of type different from $A_n$ would
be interesting to study.
\end{example}

The following is our best result providing an upper bound for smooth
schemes.

\begin{prop}
\label{boundsmooth}
Let $X$ be a smooth quasi-projective scheme over $k$. Let $\CL$ be an ample
line bundle on $X$. Then, there is $r\ge 0$ such that
$D^b(X\mqcoh)=\langle \widetilde{G}\rangle_{2\dim X+1}$
and
$D^b(X\mcoh)=\langle G \rangle_{2\dim X+1}$
where $G=\CO\oplus\CL^{\otimes -1}\oplus\cdots\oplus \CL^{\otimes -r}$.
In particular,
$\dim D^b(X\mcoh)\le 2\dim X$.
\end{prop}

\begin{proof}
There is a resolution of the diagonal
$$\cdots\to C^{-i}\Rarr{d^{-i}}\cdots\to C^{0}\Rarr{d^0}\CO_{\Delta}\to 0$$
where $C^i\in\add(\{\CL^{-j}\boxtimes\CL^{-j}\}_{j\ge 0})$. 
Denote by $C$ the complex $\cdots\to C^{-i}\Rarr{d^{-i}}\cdots\to C^{0}\to 0$.
Let $n=\dim X$.
Truncating, we get an exact sequence
$$0\to C^{-2n-1}/\ker d^{-2n}\to C^{-2n}\to\cdots\to C^{-i}\Rarr{d^{-i}}\cdots
\to C^{0}\Rarr{d^0}\CO_{\Delta}\to 0$$
Since $X\times X$ is smooth of dimension $2n$, we have
$\Ext^{2n+1}(\CO_\Delta,C^{-2n-1}/\ker d^{-2n})=0$. So, the distinguished
triangle $C^{-2n-1}/\ker d^{-2n}[2n]\to \sigma^{\ge -2n}C\to \CO_\Delta
\rightsquigarrow$ splits, \ie,
$\CO_\Delta$ is a direct summand of the complex
$\sigma^{\ge -2n}C$. We conclude as in the proof of Proposition \ref{resoldiagvar}.
\end{proof}

\begin{rem}
We actually don't know any case of a smooth variety where
$\dim D^b(X\mcoh)>\dim X$. The first case to consider would be an elliptic
curve.
\end{rem}


\begin{rem}
Let $d$ be the largest integer such that $\Ext^d_{\CO_{X\times X}}
(\CO_{\Delta X},\CF)\not=0$ for some $\CF\in (X\times X)\mcoh$. Then, 
$\dim X\le d$. We don't know if the inequality can be strict.
\end{rem}

\subsubsection{}
For applications to finite dimensional algebras, we need to prove
certain results for the derived category of differential modules.
The theory of such derived categories mirrors that of the usual
derived category
of complexes of modules (forget the grading). We state here the constructions
and results needed in this paper.

\medskip
Let $A$ be a $k$-algebra. A differential $A$-module is an
$(A\otimes_k k[\eps]/(\eps^2))$-module. We view a differential
$A$-module as a pair $(M,d)$ where $M$ is an $A$-module and
$d\in\End_A(M)$ satisfying $d^2=0$ is given by the action of $\eps$.
The cohomology of a differential
$A$-module is the $A$-module $\ker d/\im d$.

The category of differential
$A$-modules has the structure of an exact category, where the exact sequences
are those exact sequences of $(A\otimes_k k[\eps]/(\eps^2))$-modules that
split by restriction to $A$. This is a Frobenius category and its
associated stable category is called the homotopy category of differential
$A$-modules.

A morphism of $A$-modules is a quasi-isomorphism if
the induced map on cohomology is an isomorphism.
We now define the derived category of differential
$A$-modules, denoted by $D\!\diff(A)$,
as the localization of the homotopy category of differential
$A$-modules in the class of quasi-isomorphisms. These triangulated categories
have a trivial shift functor.

We have a triangulated forgetful functor $D(A)\to D\!\diff(A)$.
Let $X,Y$ be two $A$-modules and $i\ge 0$. Then, the canonical map
$\Ext^i_A(X,Y)\iso \Hom_{D(A)}(X,Y[i])\to \Hom_{D\!\diff(A)}(X,Y)$
is injective and we have an isomorphism
$\prod_{n\ge 0}\Ext^n_A(X,Y)\iso \Hom_{D\!\diff(A)}(X,Y)$.

\subsubsection{}
\begin{lemma}
\label{HHnonzero}
Let $A$ be a $k$-algebra. Let $W$ be an $A$-module with $\pdim W\ge d$.
Then, there are $A^\en$-modules $M_0=A,M_1,\ldots,M_d$ which are projective
as left and as right $A$-modules, and elements
$\zeta_i\in\Ext^1_{A^\en}(M_i,M_{i+1})$ for $0\le i\le d-1$ such that
$(\zeta_{d-1}\cdots\zeta_0)\otimes_A \id_W$ is a non zero element of
$\Ext^d_A(W,M_d\otimes_A W)$.
\end{lemma}

\begin{proof}
Let
$\cdots\to C^{-2}\Rarr{d^{-2}} C^{-1}\Rarr{d^{-1}} C^0\Rarr{d^0} A\to 0$ be a
projective resolution of the $A^\en$-module $A$. Then,
$\cdots\to C^{-2}\otimes_A W\to C^{-1}\otimes_A W\to 
C^0\otimes_A W\to W\to 0$ is a projective resolution of $W$.
Let $\Omega^{-i}$ be the kernel of $d^{i+1}$ for $i\le -1$ and
$\Omega^0=A$. Let
$\zeta_i\in\Ext^1_{A^\en}(\Omega^i,\Omega^{i+1})$ given by the exact
sequence
$0\to \Omega^{i+1}\to C^{-i}\Rarr{d^{-i}}\Omega^i\to 0$.

Since
$\Ext^d_A(W,-)$ is not zero, it follows that the exact sequence
$$0\to \Omega^d\otimes_A W\to
C^{-d+1}\otimes_A W\to\cdots\to C^{-1}\otimes_A W\to 
C^0\otimes_A W\to W\to 0$$
gives a non zero element $\xi\in\Ext^d_A(W,\Omega^d\otimes_A W)$.
This element is equal to $(\zeta_{d-1}\cdots\zeta_0)\otimes_A \id_W$.
\end{proof}

The following result is our main tool to produce lower bounds for
the dimension.
\begin{lemma}
\label{borne}
Let $A$ be a $k$-algebra. Let $W$ be an $A$-module with $\pdim W\ge d$.
Let $\CT$ be $D(A)$ or $D\!\diff(A)$.
Then,
$W\not\in\langle\overline{A}\rangle_{\CT,d}$.
\end{lemma}

\begin{proof}
Assume $W\in\langle \overline{A}\rangle_{1+r}$ for some $r\ge 0$.
Let $W_{s-1}\to W_s\to V_s\rightsquigarrow$ be a family of
distinguished triangles, for $1\le s\le r$. We put $V_0=W_0$
and we assume $V_s\in\langle \overline{A}\rangle$ for
 $0\le s\le r$ and $W_r=W\oplus W'$ for some $W'$.

We use now Lemma \ref{HHnonzero}.
The element $\zeta_i$ induces a natural transformation
of functors $M_i\otimes_A-\to M_{i+1}[1]\otimes_A-$ from $\CT$ to
itself. Restricted to $\langle\overline{A}\rangle$, this transformation
is zero. It follows from Lemma \ref{compo0} that
$(\zeta_{d-1}\cdots\zeta_0)\otimes_A -$ vanishes on
$\langle \overline{A}\rangle_d$.
It follows that
$r\ge d$ (in case $\CT$ is the derived category of differential
$A$-modules, note that the canonical map $\Ext^d_A(W,M_d\otimes_A W)\to
\Hom_\CT(W,M_d\otimes_A W)$ is injective).
\end{proof}

We deduce the following crucial Proposition~:

\begin{prop}
\label{localregular}
Let $A$ be a commutative local noetherian $k$-algebra with maximal ideal $\Gm$.
Let $\CT$ be $D(A)$ or $D\!\diff(A)$.
Then, $A/\Gm\not\in \langle A\rangle_{\CT,\Krulldim A}$.
\end{prop}

\begin{proof}
We know that $\Krulldim A\le\gldim A=\pdim_A A/\Gm$ (cf for example
\cite[Theorem 41]{Ma}). The result follows now from Lemma \ref{borne}.
\end{proof}

\begin{rem}
Let $M,V\in D(A)$.
If $V\in\langle \overline{M}\rangle_{D(A),i}$, then
$F(V)\in\langle \overline{F(M)}\rangle_{D\!\diff(A),i}$, where
$F:D(A)\to D\!\diff(A)$ is the forgetful functor.
%
\end{rem}

From Lemma \ref{borne} and Propositions \ref{boundgldimk} and \ref{regularring},
we deduce

\begin{prop}
Let $A$ be a noetherian $k$-algebra of global dimension
$d\in\BN\cup\{\infty\}$. Assume $k$ is perfect.
Then, $d$ is the minimal integer $i$ such that
$A\mperf=\langle A\rangle_{i+1}$.
\end{prop}

We can now bound dimensions~:
\begin{prop}
\label{inegvar}
Let $X$ be a reduced separated scheme of finite type over $k$. Then, 
we have $\dim D^b(X\mcoh)\ge \dim X$.
\end{prop}

\begin{proof}
Let $M\in D^b(X\mcoh)$ such that $D^b(X\mcoh)=\langle M\rangle_{r+1}$.

Pick a closed point $x$ of $X$ with local ring $\CO_x$
of Krull dimension $\dim X$ such that $M_x\in\langle \CO_X\rangle$
(given $F$ a coherent sheaf over $X$, there is a dense open affine $U$
such that $F_{|U}$ is projective. Now, a complex with
projective cohomology splits).
Then, $k_x\in \langle \CO_x\rangle_{r+1}$. It follows from Proposition
\ref{localregular} that $r\ge \Krulldim \CO_x=\dim X$.
\end{proof}

From Propositions \ref{boundgldimk} and \ref{inegvar}, we deduce

\begin{thm}
\label{smoothaffine}
Let $X$ be a smooth affine scheme of finite type over $k$. Then,
$\dim D^b(X\mcoh)=\dim X$.
\end{thm}

\begin{rem}
\label{dual}
Let $A=k[x]/(x^2)$ be the algebra of dual numbers. The indecomposable
objects of $D^b(A\mMod)$ are $k[i]$ and $L_n[i]$ for $n\ge 1$ and
$i\in\BZ$, where $L_n$ is the cone of a non-zero map
$k\to k[n]$. It follows that $D^b(A\mMod)=\langle k\rangle_{2}$,
hence, $\dim D^b(A\mMod)=1$ (cf Proposition \ref{boundLoewy} below).

Note that the dimension of the category of perfect complexes
of $A$-modules is infinite by Proposition \ref{regularring} below.
Let us prove this directly. Given 
$C$ a perfect complex of $A$-modules, there is an integer $r$
such that $\Ext^i_{A^\en}(A,A)$ acts as $0$
on $\langle C\rangle$ for $i\ge r$. On the other hand, given 
$d$ an integer, then the canonical map
$-\otimes_A\id_{L_{rd+1}}:\Ext^1_{A^\en}(A,A))^{rd}\to
\Hom_{D^b(A)}(L_{rd+1},L_{rd+1}[rd])$ is not zero
(note that $L_{rd+1}$ is perfect).
So, $L_{rd+1}\not\in \langle C\rangle_d$
by Lemma \ref{compo0}.
\end{rem}

\begin{rem}
Let $k$ be a field and $A$ a finitely generated $k$-algebra. Can the
dimension of $D^b(A\mMod)$ be infinite ? We will show
that the dimension is finite if $A$ is finite dimensional
(Proposition \ref{boundLoewy}) or commutative and $k$ is perfect
(Theorem \ref{finitevar}).
\end{rem}

\subsection{Finite global dimension}
\label{secfingldim}
\subsubsection{}
We explain here a method of d\'evissage for derived categories of abelian
categories with finite global dimension.

\begin{lemma}
\label{equivseq}
Let $\CA$ be an abelian category and $C$ a complex of objects of $\CA$.
Assume $H^1C=\cdots=H^iC=0$ for some $i\ge 0$.
Let $0\to \ker d^0\Rarr{\alpha} L^0\Rarr{f^0}\cdots\Rarr{f^i}
L^{i+1}\Rarr{\beta} C^{i+1}/\im d^i\to 0$
be an exact sequence equivalent to
$0\to \ker d^0\to C^0\to\cdots\to C^{i+1}\to C^{i+1}/\im d^i\to 0$ (\ie,
giving the same element in $\Ext^{i+2}(C^{i+1}/\im d^i,\ker d^0)$).
Then, $C$ is quasi-isomorphic to the complex
$$\cdots\to C^{-2}\Rarr{d^{-2}} C^{-1}\Rarr{a} L^0\Rarr{f^0}\cdots\Rarr{f^i}
L^{i+1}\Rarr{b} C^{i+2}\Rarr{d^{i+2}}\cdots$$
where $a$ is the composite $C^{-1}\Rarr{d^{-1}}\ker d^0\Rarr{\alpha}L^0$ and
$b$ the composite $L^{i+1}\Rarr{\beta}C^{i+1}/\im d^i\Rarr{d^{i+1}}C^{i+2}$.
\end{lemma}

\begin{proof}
It is enough to consider the case of an elementary equivalence between
exact sequences.
Let 
$$\xymatrix{
0\ar[r] & \ker d^0 \ar[r]\ar@{=}[d] & L^0 \ar[r]\ar[d] & \cdots \ar[r] &
 L^{i+1}\ar[r]\ar[d] & C^{i+1}/\im d^i\ar[r]\ar@{=}[d] & 0 \\
0\ar[r] & \ker d^0 \ar[r] & C^0 \ar[r] & \cdots \ar[r] & C^{i+1}\ar[r] &
 C^{i+1}/\im d^i\ar[r] & 0 
}$$
be a commutative diagram, with the rows being exact sequences.
Then, there is a commutative diagram
$$\xymatrix{
\cdots\ar[r] & C^{-2}\ar[r] & C^{-1} \ar[rr]\ar@{=}[dd]\ar[dr] && L^0 \ar[r]\ar[dd] & \cdots \ar[r] &
 L^{i+1}\ar[r]\ar[dd] & C^{i+2}\ar[r]\ar@{=}[dd] & \cdots \\
&&& \ker d^0 \ar[ur] \\
\cdots\ar[r] & C^{-2}\ar[r] & C^{-1} \ar[rr] && C^0 \ar[r] & \cdots \ar[r] & C^{i+1}\ar[r] &
 C^{i+2}\ar[r] & \cdots
}$$
This induces a morphism of complexes from the first row to the last
row of the diagram and this is a quasi-isomorphism.
\end{proof}

\begin{lemma}
\label{splitcplx}
Let $\CA$ be an abelian category with finite global dimension $\le n$.
Let $C$ be a complex of objects of $\CA$.
Assume $H^iC=0$ if $n\not|\ i$. Then, $C$ is quasi-isomorphic to
$\bigoplus_i (H^{ni}C)[-ni]$.
\end{lemma}

\begin{proof}
Pick $i\in\BZ$. The sequence
$0\to \ker d^{ni}\to C^{ni}\to\cdots\to C^{n(i+1)}\to
C^{n(i+1)}/\im d^{n(i+1)-1}\to 0$
is exact. It defines an element of $\Ext_\CA^{n+1}(C^{n(i+1)}/\im d^{n(i+1)-1},
\ker d^{ni})$. This group is $0$ by assumption, hence the exact sequence
is equivalent to $0\to \ker d^{ni}\to\ker d^{ni}\Rarr{0}0\cdots 0\Rarr{0}
C^{n(i+1)}/\im d^{n(i+1)-1}\to C^{n(i+1)}/\im d^{n(i+1)-1}\to 0$.
Lemma \ref{equivseq} shows that $C$ is quasi-isomorphic to a complex
$D$ with $d_D^{ni}=\cdots=d_D^{n(i+1)-1}=0$.
Now, there is a morphism of complexes $(H^{n(i+1)}C)[-n(i+1)]\to D$ that
induces an isomorphism on $H^{n(i+1)}$.
So, for every $i$, there is a map $\rho_i$ in $D(\CA)$ from
$(H^{ni}C)[-ni]$ to $C$ that induces an isomorphism on $H^{ni}$.
Let $\rho=\sum_i \rho_i:\bigoplus_i (H^{ni}C)[-ni]\to C$. This
is a quasi-isomorphism.
\end{proof}

\begin{prop}
\label{filtgldim}
Let $\CA$ be an abelian category with finite global dimension $\le n$ with
$n\ge 1$.
Let $C$ be a complex of objects of $\CA$.
Then, there is a distinguished triangle in $D(\CA)$
$$\bigoplus_i D_i\to C\to \bigoplus_i E_i\rightsquigarrow$$
where $D_i=\sigma^{\ge ni+1}\tau^{\le n(i+1)-1}C$ is a complex with
zero terms outside $[ni+1,\ldots,n(i+1)-1]$ and
$E_i$ is a complex concentrated in degree $ni$.
\end{prop}

\begin{proof}
Let $i\in\BZ$. Let $f_i$ be the composition of the canonical
maps $\tau^{\le n(i+1)-1}C\to C$ with the canonical map
$\sigma^{\ge ni+1}\tau^{\le n(i+1)-1}C\to \tau^{\le n(i+1)-1}C$. 
Then, $H^r(f_i)$ is an isomorphism for $ni+2\le r\le n(i+1)-1$ and is
surjective for $r=ni+1$.
Let
$D=\bigoplus_i \sigma^{\ge ni+1}\tau^{\le n(i+1)-1}C$ and
$f=\sum_i f_i:D\to C$. Let $E$ be the cone of $f$.
We have an exact sequence
\begin{multline*}
\cdots\to H^{ni-2}D\iso H^{ni-2}C\to H^{ni-2}E\to
H^{ni-1}D\iso H^{ni-1}C\to H^{ni-1}E\to H^{ni}D\to\\
\to H^{ni}C\to H^{ni}E\to
H^{ni+1}D\twoheadrightarrow H^{ni+1}C\to H^{ni+1}E\to
H^{ni+2}D\iso H^{ni+2}C\to\cdots
\end{multline*}
Since $H^{ni}D=0$ for all $i$, 
we deduce that
$H^rE=0$ if $n\not|\ r$. The Proposition follows now from Lemma \ref{splitcplx}.
\end{proof}

\begin{rem}
Note there is a dual version to Proposition \ref{filtgldim} obtained by passing to 
the opposite category $\CA^\opp$.
\end{rem}

\subsubsection{}
\begin{prop}
\label{boundgldim}
Let $A$ be a ring with finite global dimension.
Then, $D^b(A)=\langle \widetilde{A}\rangle_{2+2\gldim A}$.

If $A$ is noetherian, then
$D^b(A\mMod)=\langle A\rangle_{2+2\gldim A}$ and
$\dim D^b(A\mMod)\le 1+2\gldim A$.
\end{prop}

\begin{proof}
Put $n=\gldim A$. Let $C\in D^b(A)$. Up to quasi-isomorphism,
we can assume $C$ is a bounded complex of projective $A$-modules.
We now use Proposition \ref{filtgldim}. An
$A$-module $M$ has a projective resolution of length $n+1$, hence
$M\in\langle \widetilde{A}\rangle_{n+1}$. So,
$\bigoplus_i E_i\in \langle \widetilde{A}\rangle_{n+1}$.
Similarly, we have $\bigoplus_i D_i\in \langle \widetilde{A}\rangle_{n+1}$,
hence $C\in \langle \widetilde{A}\rangle_{2+2n}$.

The second part of the Lemma follows from Corollaries 
\ref{compactinb} and \ref{intercompact}.
\end{proof}

The following characterization of regular algebras is due to
Van den Bergh in the noetherian case. It characterizes regularity as
a property of $D(A)$ as a triangulated category.

\begin{prop}
\label{regularring}
Let $A$ be a ring. Then, the following conditions are equivalent
\begin{itemize}
\item[(i)]
$A$ is regular, \ie, $\gldim A<\infty$
\item[(ii)]
$K^b(A\mProj)\iso D^b(A)$
\item[(iii)]
there is $G\in D(A)^c$ and $d\in\BN$
such that $\langle\widetilde{D(A)^c}\rangle_\infty=
\langle\tilde{G}\rangle_d$
\end{itemize}
If $A$ is noetherian, these conditions are equivalent to the following
\begin{itemize}
\item[(i')]
every finitely generated $A$-module has finite projective dimension
\item[(ii')]
$D^b(A\mMod)=A\mperf$
\item[(iii')]
$A\mperf$ is strongly finitely generated.
\end{itemize}
\end{prop}

\begin{proof}
The equivalence between the first two assertions is clear, since
$D^b(A)$ is classically generated by the 
$L[i]$, where $L$ runs over the $A$-modules and
$i$ over $\BZ$.

Put $D(A)^f=\langle\widetilde{D(A)^c}\rangle_\infty$. Note that
the canonical functor $K^b(A\mProj)\iso D(A)^f$ is an equivalence.
Let $C\in D(A)$.
As in Proposition \ref{catalg}, one shows that $C\in D^b(A)$ if and only
if $\Hom(-,C)_{|D(A)^f}$ is locally finitely presented.

Assume (iii). By Theorem \ref{Brown2}, we have
$C\in D(A)^f$ if and only if  $\Hom(-,C)_{D(A)^f}$ is locally finitely
presented. So, $D^b(A)=D(A)^f$ and (ii) holds.

Finally, (i)$\Rightarrow$(iii) follows from Proposition \ref{boundgldim}.

The proof for the remaining assertions is similar.
\end{proof}

\begin{rem}
For finite dimensional or commutative algebras over a perfect field, we obtained
in Proposition \ref{boundgldimk} the better bound
$\dim D^b(A\mMod)\le \gldim A$.
We don't know whether such a bound holds under the assumption of
Proposition \ref{boundgldim}.
\end{rem}

The construction of Proposition \ref{boundgldim} is not
optimal when $A$ is hereditary, since the $D_i$'s in Proposition
\ref{filtgldim} are then zero, \ie, every object of $D^b(A)$ is
isomorphic to a direct sums of complexes concentrated in one degree. We
get then the following result.

\begin{prop}
\label{dimhered}
Let $A$ be a hereditary ring. Then,
$D^b(A)=\langle \widetilde{A}\rangle_2$.

Assume now $A$ is noetherian. Then,
$D^b(A\mMod)=\langle A\rangle_2$.
\end{prop}


\begin{rem}
Proposition \ref{dimhered} generalizes easily to quasi-hereditary algebras.
Let $\CC$ be a highest weight category over a field $k$ with weight poset $\Lambda$
(\ie, the category of finitely generated modules over a quasi-hereditary algebra).
Then, there is a decomposition $D^b(\CC)=\CI_1\diamond\cdots\diamond\CI_d$
such that $\CI_i\simeq D^b(k^{n_i}\mMod)$ for some $n_i$ and where
$d$ is the maximal $i$ such that there
is $\lambda_1<\cdots<\lambda_i\in\Lambda$ \cite[Theorem 3.9]{CPS}.
It follows from Lemma \ref{decomposition} that $\dim D^b(\CC)<d$.
\end{rem}

\begin{rem}
It would interesting to classify algebraic triangulated categories of
dimension $1$. Which differential graded / finite dimensional
algebras can have such a derived
category ? This relates to work on quasi-tilted algebras.
\end{rem}

\subsubsection{}

The following Lemma is related to the non-commutative \cite[Lemma 4.2.4]{BoVdB}.

\begin{lemma}
\label{powers}
Let $X$ be a quasi-projective scheme over a field and $\CL$
an ample sheaf. Then, there are $r,l\ge 0$ such that for any 
$n\in\BZ$, we have
$\CL^{\otimes n}\in
\add(\{G[i]\}_{|i|\le r})^{\ast l}$,
where $G=\CL^{\otimes -r}\oplus\CL^{\otimes -r+1}\oplus\cdots
\oplus\CL^{\otimes r}$.
If $X$ is regular, then we can take $l=1+\dim X$.
\end{lemma}

\begin{proof}
Pick $s>0$ such that $\CL^{\otimes s}$ is very ample and let
$i:X\to \BP^N$ be a corresponding immersion (\ie, $\CL^{\otimes s}\simeq
i^*\CO(1)$).
Beilinson's resolution of the diagonal (cf example \ref{Pn}) shows that
for every $i<0$, there is an exact sequence
of vector bundles on $\BP^N$
$$0\to \CO(i)\to \CO\otimes V_0\to\CO(1)\otimes V_1\to\cdots\to
\CO(N)\otimes V_N\to 0$$
where $V_0,\ldots,V_N$ are finite dimensional vector spaces.
By restriction to $X$, we obtain an exact sequence
$$0\to \CL^{\otimes si}\Rarr{f^{-1}} \CO\otimes V_0\Rarr{f^0}
\CL^{\otimes s}\otimes V_1\Rarr{f^1} \cdots\Rarr{f^{N-1}}
\CL^{\otimes sN}\otimes V_N\to 0.$$
We get a similar exact sequence for $i>0$ by dualizing.
This shows the first part of the Lemma with $l=N+1$.

Assume now $X$ is regular of dimension $d$. Then,
$\Ext^{d+1}(M,\CL^{\otimes si})=0$, where $M=\coker f^{d-1}$.
Consequently, $\CL^{\otimes si}$ is a direct summand of the complex
$$0\to \CO\otimes V_0\Rarr{f^0}
\CL^{\otimes s}\otimes V_1\Rarr{f^1} \cdots\Rarr{f^{d-1}}
\CL^{\otimes sd}\otimes V_d\to 0.$$
Dualizing, we see that, for $i>0$, then
$\CL^{\otimes si}$ is a direct summand of a complex
$$0\to \CL^{\otimes -sd}\otimes V_d\to\cdots\to
\CL^{\otimes -s}\otimes V_1\to \CO\otimes V_0\to 0.$$
The Lemma follows.
\end{proof}

\begin{prop}
\label{boundscheme}
Let $X$ be a regular quasi-projective scheme over a field and
$\CL$ an ample sheaf.
Then, 
$D^b(X\mqcoh)=\langle \widetilde{G}\rangle_{2(1+\dim X)^2}$ and
$D^b(X\mcoh)=\langle G\rangle_{2(1+\dim X)^2}$ for
some $r>0$, where $G=\CL^{\otimes -r}\oplus \cdots\oplus \CL^{\otimes r}$.
In particular, $\dim D^b(X\mcoh)\le 2(1+\dim X)^2-1$.
\end{prop}

\begin{proof}
By Lemma \ref{powers}, there is $r>0$ such that
$\overline{\add}(\{\CL^{\otimes i}\}_{i\in\BZ})\subset
\langle \widetilde{G}\rangle_{1+\dim X}$ for all $i$, where
$G=\CL^{\otimes -r}\oplus\cdots\oplus\CL^{\otimes r}$.
Let $C\in D^b(X\mqcoh)$. Up to isomorphism, we can assume $C$ is a
bounded complex
with terms in $\overline{\add}(\{\CL^{\otimes i}\}_{i\in\BZ})$,
because $X$ is regular.  Now, proceeding as in the proof of Proposition
\ref{boundgldim}, we get
$C\in\langle \overline{\add}(\{\CL^{\otimes i}\}_{i\in\BZ})\rangle_{2+2\dim X}$.
\end{proof}

In the case of a curve, we have a slightly better (though probably
not optimal) result.

\begin{prop}
Let $X$ be a regular quasi-projective curve over a field. Then,
$\dim D^b(X\mcoh)\le 3$.
\end{prop}

\begin{lemma}
Let $X$ be a separated scheme of finite type over $k$ and
$U$ an open subscheme of $X$.
We have $\dim D^b(U\mcoh)\le \dim D^b(X\mcoh)$.
\end{lemma}

\begin{proof}
Lemma \ref{quotient} gives the result, via
the exact sequence
$0\to D^b_{X-U}(X\mcoh)\to D^b(X\mcoh)\to D^b(U\mcoh)\to 0$,
\end{proof}

\begin{prop}
\label{regularvar}
Let $X$ be a quasi-projective scheme over $k$. Then, the following
assertions are equivalent
\begin{itemize}
\item[(i)] $X$ is regular
\item[(ii)] every object of $D^b(X\mqcoh)$ is isomorphic to a bounded
complex of locally free sheaves
\item[(iii)] $D^b(X\mcoh)=X\mperf$
\item[(iv)] $X\mperf$ is strongly finitely generated
\end{itemize}
\end{prop}

\begin{proof}
It is clear that (ii)$\Rightarrow$(i) and (iii)$\Rightarrow$(i).

By Proposition \ref{boundscheme}, we have (i)$\Rightarrow$(ii)--(iv).

Assume (iv).
Since $X\mperf$ is strongly finitely generated, it follows from
Lemmas \ref{densedim} and \ref{quotient} that
$U\mperf$ is strongly finitely generated for any affine open $U$ of $X$
because the restriction functor $X\mperf\to U\mperf$ has dense image
(Theorem \ref{locNeeman}).
So, $U$ is regular by Proposition \ref{regularring}, hence $X$ is regular. 
So, (iv)$\Rightarrow$(i).
\end{proof}

\subsection{Nilpotent ideals}
\label{secnil}
\begin{lemma}
\label{nilpideal}
Let $A$ be a noetherian ring and $I$ a nilpotent (two-sided) ideal of $A$
with $I^r=0$. Let $M\in D^b((A/I)\mMod)$ such that
$D^b((A/I)\mMod)=\langle M\rangle_n$.
Then, $D^b(A\mMod)=\langle M\rangle_{rn}$.

In particular,
$\dim D^b(A\mMod)\le r (1+\dim D^b((A/I)\mMod))-1$.
\end{lemma}

\begin{proof}
Let $C$ be a bounded complex of finitely generated $A$-modules.
We have a filtration
$0=I^rC\subset I^{r-1}C\subset\cdots\subset IC\subset C$
whose successive quotients are bounded complexes of finitely
generated $(A/I)$-modules and the Lemma follows.
\end{proof}

We have a geometric version as well.

\begin{lemma}
\label{nilgeo}
Let $X$ be a separated noetherian scheme, $\CI$ a nilpotent ideal sheaf
with $\CI^r=0$ and $i:Z\to X$ the corresponding closed immersion.
Let $M\in D^b(Z\mcoh)$  such that
$D^b(Z\mcoh)=\langle M\rangle_n$.
Then, $D^b(X\mcoh)=\langle i_*M\rangle_{rn}$.
Similarly, for
$M\in D^b(Z\mqcoh)$  such that
$D^b(Z\mqcoh)=\langle \widetilde{M}\rangle_n$, then
$D^b(X\mqcoh)=\langle \widetilde{i_*M}\rangle_{rn}$.

In particular,
$\dim D^b(X\mcoh)\le r (1+\dim D^b(Z\mcoh))-1$.
\end{lemma}

For an artinian ring $A$, the Loewy length $\Ll(A)$
of $A$ is the smallest
integer $i$ such that $J(A)^i=0$, where $J(A)$ is the Jacobson radical
of $A$.

From Lemma \ref{nilpideal}, we deduce
\begin{prop}
\label{boundLoewy}
Let $A$ be an artinian ring.
Then, $D^b(A\mMod)=\langle A/J(A)\rangle_{\Ll(A)}$.
In particular, $\dim D^b(A\mMod)\le \Ll(A)-1$.
\end{prop}

\subsection{Finiteness for derived categories of coherent sheaves}
\label{schemes}
 
Let $k$ be a field.

\subsubsection{}

The following Theorem is due to Kontsevich, Bondal and Van den Bergh for $X$ non
singular \cite[Theorem 3.1.4]{BoVdB}.

\begin{thm}
\label{finitevar}
Let $X$ be a separated scheme of finite type over a perfect field $k$. Then,
there is
$E\in D^b(X\mcoh)$ and $d\in\BN$ such that
$$D(X\mqcoh)=\langle\bar{E}\rangle_d,\
D^b(X\mqcoh)=\langle\widetilde{E}\rangle_d\text{ and }
D^b(X\mcoh)=\langle E\rangle_d.$$

In particular, $\dim D^b(X\mcoh)<\infty$.
\end{thm}

Let us explain how the Theorem will be proved. It is
enough to consider the case where $X$ is reduced. Then, the structure
sheaf of the diagonal is a direct
summand of a perfect complex up to a complex supported on $Z\times X$, where
$Z$ is a closed subscheme with smooth dense complement. We conclude
by induction by applying the Theorem to $Z$.

\smallskip
Let us start with two Lemmas.

\begin{lemma}
\label{projdiago}
Let $A$ and $B$ be two finitely generated commutative $k$-algebras, where
$k$ is perfect.
Let $M$
be a finitely generated $(B\otimes A)$-module and
$\cdots\to P^{-1}\Rarr{d^{-2}} P^0\Rarr{d^{-1}} M\Rarr{d^0} 0$ be an exact
complex with $P^i$ finitely generated and projective.

If $M$ is flat as an $A$-module and $B$ is regular of dimension $n$,
then $\ker d^{-n}$ is a
projective $(B\otimes A)$-module.
\end{lemma}

\begin{proof}
Let $i\ge 1$, $\Gm$ a maximal ideal of $A$ and $\Gn$ a
maximal ideal of $B$. We have
$$\Tor^{B\otimes A}_i(\ker d^{-n},B/\Gn\otimes A/\Gm)\simeq
\Tor^{B\otimes A}_{n+i}(M,B/\Gn\otimes A/\Gm)\simeq
\Tor^B_{n+i}(M\otimes_A A/\Gm,B/\Gn)=0$$
since $B$ is regular with dimension $n$.
It follows that $\ker d^{-n}$ is projective (cf \cite[\S 18, Lemma 5]{Ma}).
\end{proof}

\begin{lemma}
\label{closedsupport}
Let $X$ be a separated noetherian scheme and $Z$ a closed subscheme of $X$,
given by the ideal sheaf $\CI$ of $\CO_X$. For $n\ge 1$, let
$Z_n$ be the closed subscheme of $X$ with ideal sheaf $\CI^n$ and
$i_n:Z_n\to X$ the corresponding immersion.

Then, given $C\in D^b_Z(X\mcoh)$, there is $n\ge 1$ and $C_n\in D^b(Z_n\mcoh)$
such that $C\simeq i_{n*}C_n$.
\end{lemma}

\begin{proof}
Let $\CF$ be a coherent sheaf on $X$ supported by $Z$. Then,
$\CI^n\CF=0$ for some $n$ and it follows that $\CF\iso i_{n*}(i_n^*\CF)$.
More generally, a bounded complex of coherent sheaves on $X$ that are
supported by $Z$ is isomorphic to the image under $i_{n*}$ of a bounded
complex of coherent sheaves on $Z_n$ for some $n$.

Let $\CF$ be a coherent sheaf on $X$. Let $\CF_Z$ be the subsheaf of $\CF$
of sections supported by $Z$. By Artin-Rees' Theorem \cite[\S 11.C Theorem
15]{Ma}, there is an integer $r$ such that 
$(\CI^m\CF)\cap \CF_Z=\CI^{m-r} (\CI^r\cap\CF_Z)$ for $m\ge r$.
Since $\CF_Z$ is a coherent sheaf supported by $Z$, there is an integer
$d$ such that $\CI^d\CF_Z=0$. So, 
$(\CI^{r+d}\CF)\cap\CF_Z=0$. It follows that the canonical
map $\CF_Z\to \CF/(\CI^{r+d}\CF)$ is injective.

\smallskip
We prove now the Lemma by induction on the number of terms of
$C$ that are not supported by $Z$.

Let $C=0\to C^r\Rarr{d^r}\cdots\Rarr{d^{s-1}} C^s\to 0$ be a complex of
coherent sheaves
on $X$ with cohomology supported by $Z$ and take $i$ minimal such that
$C^i$ is not supported by $Z$.

Since $C^{i-1}$ and $H^i(C)$ are supported by $Z$, it follows that
$\ker d^i$ is supported by $Z$. So, there is an integer $n$ such
that the canonical map $\ker d^i\to C^i/(\CI^n C^i)$ is injective.
Let $R$ be the subcomplex of $C$ with non zero terms
$R^i=\CI^n C^i$ and $R^{i+1}=d^i(\CI^n C^i)$ --- a complex homotopy
equivalent to $0$.
Let $D=C/R$. Then, the canonical map $C\to D$ is a quasi-isomorphism.
By induction, $D$ is quasi-isomorphic to a complex of coherent sheaves
on $Z_n$ for some $n$ and the Lemma follows.
\end{proof}

\begin{proof}[Proof of Theorem \ref{finitevar}]
We have $D^b(X\mqcoh)^c=D^b(X\mcoh)$ (Corollary \ref{compactinb}). So, 
the assertion about $D^b(X\mcoh)$ follows immediately
from the one about $D^b(X\mqcoh)$ by Corollary \ref{intercompact}.
We give the proof only for the case $D^b(X\mqcoh)$, the case
of $D(X\mqcoh)$ is similar and easier.
By Lemma \ref{nilgeo}, it is enough to prove the Theorem for $X$ reduced.

\smallskip
Assume $X$ is reduced and let $d$ be its dimension.
We now prove the Theorem by induction on $d$ (the case $d=0$ is trivial).

\smallskip
Let  $U$ be a smooth dense open subscheme of $X$.
The structure sheaf $\CO_{\Delta U}$ of the diagonal
$\Delta U$ in $U\times X$ is a perfect complex
by Lemma \ref{projdiago}.
By Thomason and Trobaugh's localization Theorem 
(Theorem \ref{locNeeman}),
there is a perfect complex $C$ on $X\times X$ and
a morphism $f:C\to \CO_{\Delta X}\oplus \CO_{\Delta X}[1]$ whose restriction to
$U\times X$ is an isomorphism.
Let $G$ be a compact generator for $D(X\mqcoh)$. Then,
$G\boxtimes G$ is a compact generator for
$D((X\times X)\mqcoh)$ \cite[Lemma 3.4.1]{BoVdB}. So,
there is $r$ such that $C\in \langle G\boxtimes G\rangle_r$ by
Theorem \ref{equivgen} (3).

Let $D$ be the cone of $f$. Then, $H^*(D)$ is supported by
$Z\times X$, where $Z=X-U$.
It follows that there is a closed subscheme $Z'$ of $X$ with underlying
closed subspace $Z$, a bounded
complex $D'$ of coherent $\CO_{Z'\times X}$-modules
and an isomorphism $(i\times \id)_* D'\iso D$ in $D^b((X\times X)\mcoh)$, where
$i:Z'\to X$ is the closed immersion (Lemma \ref{closedsupport}).
By induction, there is $M\in D^b(Z'\mcoh)$ and an integer $l$ such that
$D^b(Z'\mqcoh)=\langle\widetilde{M}\rangle_l$.

\smallskip
Let $p_1$ and $p_2$ be the first and second projections $X\times X\to X$
and $\pi:Z'\times X\to Z'$ be the first projection.
Let $\CF\in D^b(X\mqcoh)$.
We have a distinguished triangle
$$Rp_{1*}(C\otimes^\BL p_2^*\CF)\to \CF\oplus \CF[1]\to 
Rp_{1*}(D\otimes^\BL p_2^*\CF)\rightsquigarrow.$$
Since $C$ is perfect, we have
$C\otimes^\BL p_2^*\CF\in D^b((X\otimes X)\mqcoh)$, hence
$Rp_{1*}(C\otimes^\BL p_2^*\CF)$ has bounded cohomology.
It follows that $Rp_{1*}(D\otimes^\BL p_2^*\CF)$ has bounded cohomology as
well.
We have
$$Rp_{1*}(D\otimes^\BL p_2^*\CF)\simeq
Rp_{1*}(i\times \id)_*(D'\otimes^\BL \BL (i\times \id)^*p_2^*\CF))\simeq
i_* R\pi_* (D'\otimes^\BL (\CO_{Z'}\boxtimes \CF))$$
Note that $R\pi_* (D'\otimes^\BL (\CO_{Z'}\boxtimes \CF))$ is
an element of $D^b(Z'\mqcoh)$.
So,
$Rp_{1*}(D\otimes^\BL p_2^*\CF)\in \langle\widetilde{i_*M}\rangle_l$.

\smallskip
We have $(G\boxtimes G)\otimes^\BL p_2^*\CF\simeq
G\boxtimes (G\otimes^\BL \CF)$, hence
$Rp_{1*}((G\boxtimes G)\otimes^\BL p_2^*\CF)\simeq
G\otimes R\Gamma(G\otimes^\BL \CF)
\in \langle\widetilde{G}\rangle$ (note this has bounded cohomology).
So, $Rp_{1*}(C\otimes^\BL p_2^*\CF)\in \langle\widetilde{G}\rangle_r$.

Finally, $\CF\in \langle\widetilde{i_*M\oplus G}\rangle_{l+r}$ and we are
done.
\end{proof}

\begin{rem}
In Theorem \ref{finitevar}, one can require $E$ to be a sheaf (consider
$\bigoplus_i H^i(E)$).
\end{rem}

\begin{rem}
Note that when $X$ is smooth, then the proof shows the stronger functorial
result  as in \S \ref{devissage} ---
this is Kontsevich's result. This stronger property does not
hold in general for singular $X$, cf the case $X=\Spec k[x]/x^2$.
\end{rem}

\begin{rem}
Theorem \ref{finitevar} does not extend to the derived categories
$D_Z(X\mcoh)$. For example, $D_{\{0\}}^b(\BA^1_k\mcoh)$ is not
strongly finitely generated.
\end{rem}

\begin{rem}
Note that the proof works under the weaker assumption that $X$ is
a separated scheme of finite type over $k$ and the residue fields
at closed points are separable extensions of $k$.

We don't know how to bound the dimension of $D^b(X\mcoh)$ for singular $X$.
When $X$ is zero dimensional, then $\dim D^b(X\mcoh)=0$ if and only if
$X$ is smooth.

We don't know whether the inequality $\dim D^b((X\times Y)\mcoh)\le
\dim D^b(X\mcoh)+\dim D^b(Y\mcoh)$ holds for $X,Y$ separated
schemes of finite type over a perfect field.

Last but not least, we don't know a single case where $X$ is smooth
and $\dim D^b(X\mcoh)>\dim X$. For example, we don't know whether
$\dim D^b(X\mcoh)=1$ or $2$ for $X$ an elliptic curve over an
algebraically closed field.
\end{rem}

We now deduce that stable derived categories are strongly finitely generated
as well.
\begin{cor}
Let $X$ be a separated scheme of finite type over a perfect field $k$
and $\CT=D^b(X\mcoh)/X\mperf$. Then, $\dim\CT<\infty$.

Assume $X$ is Gorenstein, has enough locally free sheaves
and its singular locus is complete. Then
$\CT$ is $\Ext$-finite, hence every locally finite
cohomological functor is representable.
\end{cor}

\begin{proof}
The first statement is an immediate consequence of Theorem \ref{finitevar}
and Lemma \ref{quotient}.
The fact that $\CT$ is $\Ext$-finite is \cite[Corollary 2.24 and its proof]{Or}
and the representability statement is Corollary \ref{repExtfinite}.
\end{proof}

\subsubsection{}
Let $X$ be a projective scheme over a field $k$.
Given $C\in X\mperf$ and $D\in D^b(X\mcoh)$, then
$\dim\bigoplus_{i\in\BZ}\Hom(C,D[i])<\infty$.

\smallskip
The following result is given by \cite[Theorem A.1]{BoVdB}.

\begin{lemma}
\label{swap1}
An object $D\in D(X)$ is in $D^b(X\mcoh)$ if and only if
for all $C\in X\mperf$, we have $\dim\bigoplus_{i\in\BZ}\Hom(C,D[i])<\infty$.
\end{lemma}

\begin{proof}
The first implication has been recalled before.

Let $D\in D(X)$ such that $\Hom(-,D)_{|X\mperf}$ is locally finite.
Then, $\Hom(-,D)_{|X\mperf}$ is locally finitely presented (Proposition
\ref{locfin}), hence $D\in D^b(X\mcoh)$ (Proposition \ref{clfpcoh}).
\end{proof}

\begin{prop}
There is a fully faithful functor $S:X\mperf\to D^b(X\mcoh)$ and bifunctorial
isomorphisms
$$\Hom(C,D)^*\iso \Hom(D,S(C))$$
for $C\in X\mperf$ and $D\in D(X)$.
\end{prop}

\begin{proof}
The category $D(X)$ is cocomplete and has a compact generator
(Theorem \ref{genscheme}). Corollary \ref{Serre} shows the existence
of a functor $S:X\mperf\to D(X)$. 

By Lemma \ref{swap1},
if $C\in X\mperf$, then $S(C)\in D^b(X\mcoh)$.
\end{proof}

We can now prove a ``dual version'' of Lemma \ref{swap1}~:

\begin{lemma}
\label{swap2}
An object $C\in D(X)$ is in $X\mperf$ if and only if
for all $D\in D^b(X\mcoh)$, we have
$\dim\bigoplus_{i\in\BZ}\Hom(C,D[i])<\infty$.
\end{lemma}

\begin{proof}
The first implication has been recalled before.

Let $C\in D(X)$ such that
for all $D\in D^b(X\mcoh)$, we have
$\dim\bigoplus_{i\in\BZ}\Hom(C,D[i])<\infty$.
Let $D'\in X\mperf$. Then,
$$\Hom(D',C[i])^*\iso \Hom(C,S(D')[i]),$$
hence $\dim\bigoplus_{i\in\BZ}\Hom(D',C[i])<\infty$. It follows from
Lemma \ref{swap1} that $C\in D^b(X\mcoh)$.

Let $x$ be a closed point of $X$. We have 
$\dim\bigoplus_i\Hom(C,\CO_{\{x\}}[i])=
\dim\bigoplus_i \Hom(C_x,\CO_{\{x\}}[i])<\infty$.
This shows that $C_x$ is a perfect complex of $\CO_x$-modules.
Since $C\in D^b(X\mcoh)$, we deduce that $C$ is perfect.
\end{proof}

The following result was conjectured by Bondal --- the first statement
is \cite[Theorem A.1]{BoVdB}.

\begin{cor}
\label{swap}
Let $X$ be a projective scheme over a perfect field $k$.
\begin{itemize}
\item[(i)]
Every locally finite cohomological functor $(X\mperf)^{\opp}\to k\mMod$
is representable by an object of $D^b(X\mcoh)$.
\item[(ii)]
Every locally finite cohomological functor $D^b(X\mcoh)\to k\mMod$
is representable by an object of $X\mperf$.
\end{itemize}
\end{cor}

\begin{proof}
By Remark \ref{CKN}, a finite cohomological functor $(X\mperf)^{\opp}\to k\mMod$
is representable by an object of $D(X)$ and Lemma
\ref{swap1} says that the object must be in $D^b(X\mcoh)$. This shows (i).

By Theorem \ref{finitevar}, $D^b(X\mcoh)^\opp$ is strongly finitely
generated. So, Proposition \ref{locfin} and Corollary \ref{strongKaroubian}
show that every locally finite cohomological functor $D^b(X\mcoh)\to k\mMod$
is representable by an object of $D^b(X\mcoh)$ and Lemma
\ref{swap2} says that the object must be in $X\mperf$. This shows (ii).
\end{proof}

\begin{rem}
Similar results should hold for $X$ quasi-projective, with
$D^b(X\mcoh)$ replaced by its full subcategory of objects with compact support.
\end{rem}

\section{Applications to finite dimensional algebras}
\label{secfd}
\subsection{Auslander's representation dimension}
\subsubsection{}
Let $\CA$ be an abelian category.

\begin{defi}
\label{defAuslander}
The (Auslander) representation dimension $\repdim\CA$ is the
smallest integer $i\ge 2$ such that there is an object $M\in\CA$ with
the property that given any $L\in\CA$,
\begin{itemize}
\item[(a)]
there is an exact sequence
$$0\to M^{-i+2}\to M^{-i+3}\to\cdots\to M^0\to L\to 0$$
with $M^j\in\add(M)$ such that the sequence 
$$0\to \Hom(M,M^{-i+2})\to \Hom(M,M^{-i+3})\to\cdots\to \Hom(M,M^0)\to 
\Hom(M,L)\to 0$$
is exact
\item[(b)]
there is an exact sequence
$$0\to L\to {M'}^0\to {M'}^1\to\cdots\to {M'}^{i-2}\to 0$$
with ${M'}^j\in\add(M)$ such that the sequence 
$$0\to \Hom({M'}^{i-2},M)\to\cdots\to \Hom({M'}^1,M)\to\Hom({M'}^0,M)\to 
\Hom(L,M)\to 0$$
is exact.
\end{itemize}
\end{defi}

An object $M$ that realizes the minimal $i$ is called an {\em Auslander
generator}.

\smallskip
Note that either condition (a) or (b) implies that
$\gldim \End_\CA(M)\le i$, and
$\gldim \End_\CA(M)=i$ if $M$ is an Auslander generator
(cf \eg\ \cite[Lemma 2.1]{ErHoIySc}).
Note also that if condition (a) (resp. (b)) hold
for every $L$ in a dense subcategory $\CI$ of $\CA$, then, it holds for 
every object of $\CA$.

\smallskip
Note that $\repdim\CA=2$ if and only if $\CA$ has only finitely many
isomorphism classes of indecomposable objects.
Note also that $\repdim\CA=\repdim\CA^\opp$.

\subsubsection{}
Take $\CA=A\mMod$, where $A$ is a finite dimensional algebra over a field.
Then, we write $\repdim(A)$ for $\repdim(A\mMod)$.

Let $M\in\CA$ and $i\ge 2$. If $M$ satisfies (a) of Definition \ref{defAuslander},
then it contains a projective generator as a direct summand
(take $L=A$). More generally, the following are equivalent
\begin{itemize}
\item
$M$ satisfies (a) of Definition \ref{defAuslander}
and $M$ contains an injective cogenerator as a direct summand 
\item
$M$ satisfies (b) of Definition \ref{defAuslander}
and $M$ contains a projective generator as a direct summand 
\item
$M$ satisfies (a) and (b) of Definition \ref{defAuslander}.
\end{itemize}
So, the definition of representation dimension given here
coincides with Auslander's original definition (cf \cite[\S III.3]{Au} and
\cite[Lemma 2.1]{ErHoIySc}) when $A$ is not semi-simple. When $A$ is semi-simple,
Auslander assigns the representation
dimension $0$ whereas we define it to be $2$ here.
Iyama has shown \cite{Iy} that the representation dimension of
a finite dimensional algebra is finite.

Various classes of  algebras with representation dimension $3$
have been found~: algebras with radical square zero
\cite[\S III.5, Proposition p.56]{Au}, hereditary
algebras \cite[\S III.5, Proposition p.58]{Au} and more generally
stably hereditary algebras \cite[Theorem 3.5]{Xi},
special biserial algebras \cite{ErHoIySc}, local algebras
of quaternion type \cite{Ho}.

\subsubsection{}
\label{weak}
One can weaken the requirements in the definition of the representation
dimension as follows~:

\begin{defi}
\label{defweak}
The weak (resp. left weak, resp. right weak) representation dimension of $\CA$,
denoted by $\wrepdim(\CA)$ (resp. $\lwrepdim(\CA)$, resp. $\rwrepdim(\CA)$)
is the smallest integer $i\ge 2$ such that there is an object $M\in\CA$ with
the property that given any $L\in\CA$, 
there is a bounded complex $C=0\to C^r\to\cdots\to C^s\to 0$ of $\add(M)$ with 
\begin{itemize}
\item
$L$ isomorphic to a direct summand of $H^0(C)$
\item
$H^d(C)=0$ for $d\not=0$ and
\item
$s-r\le i-2$
(resp. and $C^d=0$ for $d>0$, resp. and $C^d=0$ for $d<0$).
\end{itemize}
\end{defi}

Note that $\lwrepdim \CA=\rwrepdim \CA^\opp$, $\wrepdim\CA=\wrepdim\CA^\opp$,
$$\inf\{\lwrepdim\CA,\rwrepdim \CA\}\ge \wrepdim\CA\textrm{ and }
\repdim\CA\ge \sup\{\lwrepdim\CA,\rwrepdim \CA\}.$$

\begin{prop}
\label{inegrepdim}
We have
$\dim D^b(A\mMod)\le \repdim(A)$
and $\dim \left(D^b(A\mMod)/A\mperf\right)\le \wrepdim(A)-2$.
\end{prop}

\begin{proof}
Let $M$ be an Auslander generator for $A\mMod$.
Let $C$ be a bounded complex of objects of $\add(M)$. Let
$L\in A\mMod$ and let $D$ be a bounded complex of $\add(M)$ together
with a map $f:D\to L$ such that $\Hom(M,f)$ is a quasi-isomorphism
(in particular, $f$ is a quasi-isomorphism). Then,
$\Hom^\bullet(C,f):\Hom^\bullet(C,D)\to\Hom^\bullet(C,L)$ is a
quasi-isomorphism, hence
$\Hom(C,f):\Hom_{K^b(A)}(C,D)\to\Hom_{K^b(A)}(C,L)$ is an isomorphism.

It follows by induction that every bounded complex of $A\mMod$
is quasi-isomorphic to a bounded complex of $\add(M)$, \ie,
the canonical functor $K^b(\add(M))\to D^b(A\mMod)$ is essentially surjective.
We have equivalences $K^b(\End(M)\mproj)\iso D^b(\End(M))$
and $K^b(\End(M)\mproj)\iso K^b(\add(M))$ and
$\dim D^b(\End(M))\le \gldim\End(M)$ by Proposition \ref{boundgldimk}. So,
$\dim D^b(A\mMod)\le \repdim(A)$ by Lemma \ref{quotient}.

\smallskip
Let $n=\wrepdim A$. There is $N\in A\mMod$ with the property that given
$L\in A\mMod$,
there is a bounded complex of $\add(N)$
$C=0\to C^r\to\cdots \to C^s\to 0$
with $H^i(C)=0$ for $i\not=0$, $L$ is a direct summand of $H^0(C)$ and
$s-r\le n-2$.
Then $L\in \langle C^s\rangle\diamond\cdots\diamond\langle C^r\rangle$.
Every object of $D^b(A\mMod)/A\mperf$ is isomorphic to an
object $L[r]$ for some $L\in A\mMod$ and $r\in\BZ$. Consequently,
$D^b(A\mMod)/A\mperf=\langle M\rangle_{n-1}$.
\end{proof}

\smallskip
In order to obtain lower bounds for the representation dimension
of certain algebras, we will actually construct lower bounds for
the weak representation dimension.

\begin{rem}
Note that the representation dimension as well as the invariants
of Definition \ref{defweak} are not invariant by derived equivalence
(consider for instance a derived equivalence between an algebra with finite
representation type and an algebra with infinite representation type).
\end{rem}

\begin{rem}
All the definitions given here for abelian categories make sense for
exact categories.
\end{rem}

\subsection{Stable categories of self-injective algebras}
Let $k$ be a field.
\subsubsection{}
For $A$ a self-injective finite dimensional $k$-algebra,
we denote by $A\mstab$ the stable category of
$A$. This is the quotient of the
additive category $A\mMod$ by the additive subcategory $A\mproj$.
The canonical functor $A\mMod\to D^b(A\mMod)$ induces an equivalence
$A\mstab\iso D^b(A\mMod)/A\mperf$ (\cite[Example 2.3]{KeVo} and
\cite[Theorem 2.1]{Ri}).
This provides $A\mstab$ with a structure of triangulated category.
Recall that $\Ll(A)$ denotes the Loewy length of $A$ (cf \S \ref{secnil}).

\begin{prop}
\label{Loewy}
Let $A$ be a non-semisimple self-injective algebra. Then,
$$\Ll(A)\ge \repdim A\ge \wrepdim A\ge 2+\dim A\mstab.$$
\end{prop}

\begin{proof}
The first inequality is \cite[\S III.5, Proposition p.55]{Au} (use
$M=A\oplus A/J(A)\oplus A/J(A)^2\oplus\cdots$).
The second inequality is trivial (cf \S \ref{weak}).
The last inequality is given by Proposition \ref{inegrepdim}.
\end{proof}
\subsubsection{}
We study here self-injective algebras with representation dimension $3$.

Recently, various properties have been found for
algebras of representation dimension $3$ (cf for example \cite{IgTo}).
Here is a result in this direction concerning self-injective algebras.

Note that $\repdim A=2$ if and only if $\dim A\mstab=0$. Consequently,
if $\repdim A=3$, then $\dim A\mstab=1$ (cf Proposition \ref{Loewy}).

Given $M$ an $A$-module, we denote by $\Omega M$ the kernel of a surjective
map from a projective cover of $M$ to $M$ and by $\Omega^{-1}M$ the 
cokernel of an injective map from $M$ to an injective hull of $M$.

\begin{lemma}
\label{simple2}
Let $A$ be a self-injective $k$-algebra and $C=0\to C^0\to C^1\to 0$ an
indecomposable complex of finitely generated $A$-modules with
$H^0(C)=0$ and $H^1(C)=S$ simple. Then,
\begin{itemize}
\item $C^0$ and $C^1$ have no non-zero projective direct summand
\item or $C^1$ is projective indecomposable and $C^0\simeq\Omega S$.
\end{itemize}
\end{lemma}

\begin{proof}
If $C^0$ has a non-zero projective summand $L$, then $L$ is
injective and the restriction of $d=d_C^0$ to $L$ is a split injection. In
particular, $C$ has a direct summand isomorphic to $0\to L\Rarr{\id}L\to 0$,
which is impossible.

Assume now that $C^1$ has a submodule $N$ such that $C^1/N$ is
projective indecomposable. If $N\not\subseteq \im d$, then there is
$P\subseteq\im d$ such that $C^1=N\oplus P$. So, $C$ has
a direct summand isomorphic to $0\to P\Rarr{\id}P\to 0$~: this
is impossible. So, $C^0\iso \im d=N\oplus N'$ with $N'\simeq\Omega S$.
The indecomposability of $C$ gives $N=0$.
\end{proof}

\begin{lemma}
Let $A$ be a self-injective $k$-algebra with $\repdim A=3$ and 
$M$ an Auslander generator. Assume $\Omega^{-1}M$ has no simple direct summand.
Then, the number of simple $A$-modules (up to isomorphism) is less than or
equal to the number of isomorphism classes of non projective
indecomposable summands of $M$.
\end{lemma}

\begin{proof}
Since $\repdim A=3$,
for every simple $A$-module $S$, there is
an exact sequence
$0\to M_1\to M_0\to S\to 0$ with $M_0$ and $M_1$ in $\add(M)$.
By Lemma \ref{simple2},
we can assume
that $M_0$ and $M_1$ have
no non-zero projective direct summands. Then, we have
$[S]=[M_0]-[M_1]$ in $K_0(A\mMod)$. It follows that the non-projective
indecomposable summands of $M$ generate $K_0(A\mMod)$.
\end{proof}

Let $A$ and $B$ be two self-injective algebras. A stable equivalence
of Morita type between $A$ and $B$ is the data of a finite dimensional $(A,B)$-bimodule
$X$, projective as an $A$-module and as a right $B$-module,
and of a finite dimensional $(B,A)$-bimodule $Y$,  
projective as a $B$-module and as a right $A$-module, such that
$$X\otimes_B Y\simeq A\oplus\text{ projective as }(A,A)-\text{bimodules}$$
$$Y\otimes_A X\simeq B\oplus\text{ projective as }(B,B)-\text{bimodules}.$$

Stable equivalences of Morita type preserve the representation dimension
\cite[Theorem 4.1]{Xi}~:

\begin{prop}
\label{auslanderstable}
Let $A$ and $B$ be two self-injective $k$-algebras and
$X$ be an $(A,B)$-bimodule inducing a stable equivalence between $A$ and $B$.

Let $M$ be an Auslander generator for $B$.
Then, $X\otimes_B M$ is an Auslander generator for $A$.

In particular, $\repdim A=\repdim B$.
\end{prop}

\begin{proof}
Given $M$ containing a progenerator as a direct summand, the property
(a) (resp. (b)) for $L$ in Definition \ref{defAuslander}
is equivalent to the same property for $L\oplus P$,
where $P$ is some fixed projective module.

Let $Y$ be a $(B,A)$-bimodule inverse to $X$.
Let $V$ be an $A$-module. Then,
$X\otimes_B Y\otimes_A V\simeq V\oplus P$ with $P$ projective. Starting with an exact
sequence
resolving $Y\otimes_A V$ as in (a) or (b) of Definition \ref{defAuslander},
we get one for $V\oplus P$ by applying $X\otimes_B-$
Now, applying $\Hom_B(X\otimes_B M,-)$ to that new exact sequence gives
the same result as applying $\Hom_A(Y\otimes_A X\otimes_B M,-)$ to the original
exact sequence. Since $Y\otimes_A X\otimes_B M\simeq M\oplus\text{ projective}$,
we indeed get an exact sequence.
\end{proof}

The following Proposition gives a bound for the number of non-projective
simple modules of a self-injective algebra which is stably equivalent
(\`a la Morita) to a given self-injective algebra with representation
dimension $3$.

\begin{prop}
Let $A$ be a self-injective $k$-algebra with $\repdim A=3$ and 
$M$ an Auslander generator. Let $B$ be a self-injective $k$-algebra.
Assume there is a stable equivalence of Morita type
between $A$ and $B$.

Then, the number of simple non-projective $B$-modules (up to isomorphism)
is less than or equal to twice the number 
of isomorphism classes of indecomposable summands of $M$.
\end{prop}

\begin{proof}
Replacing $B$ by a direct factor, one can assume $B$ has no
simple projective module.

Let $Y$ be a $(B,A)$-bimodule inducing a stable equivalence
and $N=Y\otimes_A M$. Then, $N$ is an Auslander generator for $B$ and
$\repdim B=3$ (Proposition \ref{auslanderstable}).
Let $R$ be the subgroup of $K_0(B\mMod)$
generated by the classes of the non-projective indecomposable
summands of $N$. Note that the rank of $R$ is at most the number of
isomorphism classes of indecomposable non-projective summands of $M$.

Let $S$ be a simple $B$-module with $[S]\not\in R$.
There is an exact sequence
$0\to N_1\to N_0\to S\to 0$ with $N_0$ and $N_1$ in $\add(N)$ and by
Lemma \ref{simple2},
$N_0$ is a projective cover of $S$ and $N_1\simeq \Omega S$. In
particular, $\Omega S$ is a direct summand of $N$.

So, the number of simple $B$-modules with $[S]\not\in R$ is
at most the number of isomorphism classes of
indecomposable non-projective summands of $N$.
\end{proof}

\begin{rem}
This Proposition, which was the starting point of this paper, led us to
investigate the existence of self-injective algebras
with representation dimension greater than $3$. This Proposition is
related to the problem of the equality of the number of simple non projective
modules for two stably equivalent algebras.
\end{rem}

\subsubsection{}

The following Theorem gives the first known examples of algebras with
representation dimension $>3$.

\begin{thm}
\label{exterior}
Let $n\ge 1$ be an integer. Then,
$\dim \Lambda(k^n)\mstab=\repdim \Lambda(k^n)-2=n-1$ and
$\dim D^b(\Lambda(k^n)\mMod)=n$.
\end{thm}

\begin{proof}
Put $A=\Lambda(k^n)$ and $B=k[x_1,\ldots,x_n]$.

Let us recall a version of Koszul duality 
\cite[\S 10.5, Lemma ``The `exterior' case'']{Ke}.
We have an equivalence of triangulated categories $R\Hom^\bullet(k,-)$
between $D^b(A\mMod)$ and $\CT$, the subcategory of
the derived category of differential graded $B$-modules classically
generated by $B$ ($\CT$ is also the subcategory of compact objects by
Corollary \ref{genDG}).
Note that this is a special case of \S \ref{dgend}, using the fact that
$R\End^\bullet(k)$ is a dg algebra quasi-isomorphic to its cohomology algebra
$B$.
This equivalence sends $A$ to $k$, so it
induces an equivalence of triangulated categories between $A\mstab$ and
$\CT/\CI$, where $\CI$ is the thick subcategory of $\CT$ classically
generated by $k$.
Denote by $F:\CT\to\CT/\CI$ the quotient functor.

Let $M\in\CT$ such that $\CT/\CI=\langle F(M)\rangle_{\CT/\CI,r+1}$.
Up to isomorphism,
we can assume $M$ is finitely generated and projective as a $B$-module.
Let $\CF$ be the sheaf over $\BP^{n-1}$ corresponding to the graded
$B$-module $M$. The differential on $M$ gives a map
$d:\CF\to \CF(1)$. Let $\CG=\ker d(1)/\im d$.
Pick $x$ a closed point of $\BP^{n-1}$ such that
$\CG_x$ is a projective $\CO_x$-module. Then,
there is a projective $\CO_x$-module $R$ such that
$\ker d_x=\im d_x\oplus R$.
We have an exact sequence $0\to R\to \CF_x\to \CF_x/R\to 0$ of
differential $\CO_x$-modules. Since $\CF_x/R$ is acyclic, it follows that
$R\to \CF_x$ is an isomorphism in $D\!\diff(\CO_x)$, 
the derived category of differential $\CO_x$-modules.
Let $I(x)$ be the prime ideal of $B$ corresponding to the line $x$ of $\BA^n$.
Note that the differential graded $B$-module $B/I(x)$ (the differential
is $0$) is in $\CT$.
So, $F(B/I(x))\in \langle F(M)\rangle_{\CT/\CI,r+1}$, hence
$k_x\in \langle \CF_x\rangle_{D\!\diff(\CO_x),r+1}$, hence
$k_x\in \langle \CO_x\rangle_{D\!\diff(\CO_x),r+1}$.
By Proposition \ref{localregular}, we get
$r\ge n-1$. Hence, $\dim A\mstab\ge n-1=\Ll(A)-2$. Now, Proposition \ref{Loewy}
gives the conclusion.

\smallskip
The proof of the inequality $\dim D^b(\Lambda(k^n)\mMod)\ge n$ is
similar (and easier). Proposition \ref{boundLoewy} gives the inequality
$\dim D^b(\Lambda(k^n)\mMod)\le n$.
\end{proof}

\subsubsection{}

We assume here that $k$ is a field of characteristic $p>0$.

\begin{prop}
\label{defect}
Let $G$ be a finite group and $B$ a block of $kG$.
Let $D$ be a defect group of $B$. Then,
$\dim B\mstab=\dim (kD)\mstab$ and
$\dim D^b(B\mMod)=\dim D^b((kD)\mMod)$
\end{prop}

\begin{proof}
Recall that a defect group $D$ of $B$ is a (smallest) subgroup such that
the identity functor of $B\mMod$ is a direct summand of
$\Ind_D^G\Res_D^G$.
Since $kD$ is a direct summand of $B$
as a $(kD,kD)$-bimodule, the Proposition follows from
Lemma \ref{quotient}.
\end{proof}

Given $P$ a finite $p$-group and $Q$ a maximal subgroup of
$P$, we denote by $\beta_Q\in H^2(P,\BZ/p)$ the class
of the exact sequence
$$0\to \BZ/p\to \Ind_Q^P\BZ/p\xrightarrow{x-1}\Ind_Q^P\BZ/p\to \BZ/p\to 0$$
where $x\in P-Q$.

The following Proposition gives a recursive bound for the dimension
of the stable category.

\begin{prop}
Let $G$ be a finite group and $B$ a block of $kG$.
Let $D$ be a defect group of $B$. Let
$D_1,\ldots,D_n$ be a family of maximal subgroups of $D$ such that
$\beta_{D_1}\cdots \beta_{D_n}=0$ (such a family exists and one
can assume $n\le \frac{p+1}{p^2}|D:\Phi(D)|$).

Then, $\dim B\mstab<2\sum_i(1+\dim kD_i\mstab)$ and
$\dim D^b(B\mMod)<2\sum_i(1+\dim D^b(kD_i\mMod))$.
\end{prop}

\begin{proof}
By Proposition \ref{defect}, it is enough to consider the
case where $G=D$. Then, \cite[Lemma 3.9]{Ca} asserts that
there is a $kG$-module $M$ which has $k$ as a direct summand and
has a filtration
$0=M_0\subset M_1\subset\cdots\subset M_{2n}=M$
with $M_{2i-1}/M_{2i-2}\simeq \Ind_{D_i}^G\Omega^{t_i}k$ and
$M_{2i}/M_{2i-1}\simeq \Ind_{D_i}^G\Omega^{t'_i}k$ for some
integers $t_i,t'_i$. We conclude as in \S \ref{devissage}.

The existence of the family is Serre's Theorem on product of Bockstein's,
cf \eg\ \cite[Theorem 7.4.3]{Ben}.
\end{proof}

\begin{thm}
\label{2group}
Let $G$ be a finite group, $B$ a block of $kG$ over a field $k$
of characteristic $2$. Let $D$ be a defect group of $B$.
Then, $\repdim B\ge 2+\dim B\mstab >r$ and
$\dim D^b(B\mMod)\ge r$, where $r$ is the $2$-rank of $D$.
\end{thm}

\begin{proof}
The first inequality is given by Proposition \ref{Loewy}.
By Proposition \ref{defect}, it suffices to prove the Theorem for
$G=D$ and $B=kD$.
Let $P$ be an elementary abelian $2$-subgroup of $D$ with rank the
$2$-rank of $D$. Then,
$\dim kP\mstab\le \dim kD\mstab$ by Lemma \ref{quotient}.
Now, $kP\simeq \Lambda (k^r)$ and the Theorem follows from Theorem
\ref{exterior}.
The derived category assertion has a similar proof.
\end{proof}

Let us recall a conjecture of D.~Benson~:

\begin{conj}[Benson]
\label{conjBenson}
Let $G$ be a finite group, $B$ a block of $kG$ over a field
$k$ of characteristic $p$. Then,
$\Ll(B)> p$-$\mathrm{rank}(D)$.
\end{conj}

From Theorem \ref{2group} and Proposition \ref{Loewy}, we deduce~:

\begin{thm}
\label{Benson2}
Benson's conjecture \ref{conjBenson} holds for $p=2$.
\end{thm}

\end{document}